\documentstyle[amsfonts,amsmath,theorem,epsfig,subfigure,twoside]{report}
\def\shadowbox{\hbox{\rule[-0.0ex]{0.1ex}{1.2ex}%
\hspace{-0.1ex}\rule[-0.0ex]{1.2ex}{0.1ex}%
\hspace{0.0ex}\rule[-0.0ex]{0.1ex}{1.2ex}\hspace{-1.3ex}%
\rule[1.15ex]{1.25ex}{0.1ex}\hspace{-0.0ex}\rule[-0.25ex]{0.3ex}{1.1ex}%
\hspace{-1.2ex}\rule[-0.25ex]{1.1ex}{0.25ex}}}
\def\qed{\ifmmode \hbox{\hfill\shadowbox}
     \else \hphantom{x}\hfill\shadowbox \fi}

\newenvironment{myproof}{\noindent{\bf Proof.}\rm}{\par\medskip}

\newtheorem{theorem}{Theorem}[section]
\newtheorem{lemma}[theorem]{Lemma}
\newtheorem{corollary}[theorem]{Corollary}

\newtheorem{conjecture}[theorem]{Conjecture}

\theorembodyfont{\rmfamily}
\newtheorem{definition}[theorem]{Definition}
\theoremheaderfont{\sc\bf}

\newtheorem{example}[theorem]{Example}
\def\remark{{\noindent \bf Remark:\hspace{0.5em}}}

\def \R{{\mathbb R}}
\def \C{{\mathbb C}}
\def \Q{{\mathbb Q}}
\def \N{{\mathbb N}}
\def \Z{{\mathbb Z}}
\def \T{{\mathbb T}}
\def \H{{\mathbb H}}
\def \A{{\mathcal A}}
\def \B{{\mathcal B}}
\def \D{{\mathcal D}}
\def \F{{\mathcal F}}
\def \S{{\mathcal S}}
\def \O{{\Omega}}
\def \Hi{{\mathcal H}}
\def \Hchan{{\mathbf H}}
\def \Ll{{\mathcal L}}
\def \Lone{{L^1(\R^d)}}
\def \Ltwo{{L^2(\R^d)}}
\def \Lp{{L^{ ^{\scriptstyle p}}(\R^d)}}
\def \Lq{{L^{ ^{\scriptstyle q}}(\R^d)}}
\def \Linf{{L^{ ^{\scriptstyle \infty}}(\R^d)}}
\def \lq{{\ell^{ ^{\scriptstyle q}}(\Z^d)}}
\def \lp{{\ell^{ ^{\scriptstyle p}}(\Z^d)}}
\def \ll{{\ell^{ ^{\scriptstyle 2}}(\Z^d)}}

\def \a{{\alpha}}
\def \b{{\beta}}
\def \d{{\delta}}
\def \e{{\epsilon}}
\def \g{{\gamma}}
\def \G{{\Gamma}}
\def \k{{\kappa}}
\def \p{{\psi}}
\def \t{{\tau}}
\def \w{{\omega}}
\def \l{{\ell}}
\def \la{{\lambda}}
\def \La{{\Lambda}}
\def \m{{\mu}}
\def \n{{\nu}}
\def \r{{\rho}}
\def \th{{\theta}}
\def \e{{\epsilon}}
\def \s{{\sigma}}
\def \x{{\xi}}
\def \z{{\zeta}}
\def \La{{\Lambda}}
\def \ph{{\varphi}}
\def \Ph{{\Phi}}

\def \be{\begin{equation}}
\def \bea{\begin{eqnarray}}
\def \bpm{\begin{bmatrix}}
\def \bbm{\begin{bmatrix}}
\def \bc{\begin{center}}
\def \ec{\end{center}}

\def \nn{{\nonumber}}
\def \np{{\newpage}}

\def\Lsp{{\boldsymbol L}}
\def\lsp{{\boldsymbol\ell}}
\def\LtR{{\Lsp^2(\R)}}
\def\ltZ{{\lsp^2(\Z)}}

\def\ord{{\cal O}}

\def\cond{\mbox{cond\/}}
\def\condR{{\cond(R_{\ggam,T,F})}}

\def\Dil{{\cal D}_{\alpha}}
\def\Chirp{{\cal C}_{\beta}}

\def\SQI{{S^{-\frac{1}{2}}}}

\def\eps{\varepsilon}
\def\supp{{supp}}
\def \Ps{{\Psi}}
\def \Ph{{\Phi}}

\newcommand{\ggam}{{g_{\sigma}}}

\newcommand{\psikl}{\ph_{k,\l}}
\newcommand{\phikl}{\phi_{k,\l}}
\newcommand{\psigam}{{\ph_{\sigma}}}

\newcommand{\Dtf}{\Delta \tau_{f}}
\newcommand{\Dnf}{\Delta \nu_{f}}
\newcommand \wh{\widehat}
\newcommand \wt{\widetilde}

\begin{document}
\pagenumbering{roman}
\pagestyle{plain}
\renewcommand{\baselinestretch}{1.4}\large
\begin{center}
{\Large \bf Banach Algebras of Integral Operators, Off-Diagonal Decay, and Applications in Wireless Communications}
\vskip 0.2in
\centerline{By}
\vskip 0.1in
\centerline{\large  SCOTT BEAVER}
\centerline{\large B.S. Engineering Physics (Lehigh University) 1990}
\centerline{\large B.S. Mathematics (Auburn University at Montgomery) 1994}
\centerline{\large M.S. Mathematics (University of Arizona) 1997}
\vskip 0.3in
\centerline{\large  DISSERTATION}
\vskip 0.1in
\centerline{\large Submitted in partial satisfaction of the requirements for the
degree of}
\vskip 0.08in
\centerline{\large  DOCTOR OF PHILOSOPHY}
\vskip 0.08in
\centerline{\large in}
\vskip 0.08in
\centerline{\large MATHEMATICS}
\vskip 0.08in
\centerline{\large in the}
\vskip 0.08in
\centerline{\large  OFFICE OF GRADUATE STUDIES}
\vskip 0.08in
\centerline{\large of the}
\vskip 0.08in
\centerline{\large  UNIVERSITY OF CALIFORNIA}
\vskip 0.08in
\centerline{\large  DAVIS}
\end{center}
\vskip 0.08in
\centerline{\large Approved:}
\vskip 0.18in
\centerline{\underbar{\hskip 2.5in}}
\vskip 0.14in
\centerline{\underbar{\hskip 2.5in}}
\vskip 0.14in
\centerline{\underbar{\hskip 2.5in}}
\vskip 0.14in
\centerline{\large Committee in Charge}
\vskip 0.08in
\centerline{\large 2004}

\np
\renewcommand{\baselinestretch}{1.05}
\large
\tableofcontents


\newpage
\large
{\Large \bf ACKNOWLEDGEMENTS} \\

Permit me to take a moment to give thanks to those who have given me aid and comfort during my years in graduate school.  First and foremost, I thank my wife Cheryl who supported me both emotionally and financially, and who was a constant inspiration and a source of mathematical wealth.  Additional inspiration has been provided by my daughter Elizabeth, who at the time of this writing is almost fourteen months old.  I am also deeply indebted to my advisor Thomas Strohmer, for whose patience, guidance, and seemingly endless supply of problems on which to work I am most grateful.  I thank my parents Frederick and Barbara Beaver, my other family members Robin, Matthias, and Joshua Limmeroth, and Glorian and Gene DeLorme for their love and support.  I express my most heartfelt thanks to Mark Gorden, my great friend in Davis who gave me a place to stay when I made the many trips from Albuquerque to Davis, but more importantly, for whose forthright manner, insights into life, and terrific conversations I am the richer.  I also wish to thank two great friends and fraternity brothers Todd Radano and Kevin Parmenter, who were there with a kick in the tail when I needed it most.  I thank professor Bruno Nachtergaele, who, while in the capacity of Chair of the Graduate Group in Mathematics, always provided wise counsel.  I must also thank my friends Mike Becker, Mark Banaszek, Angel Pineda, Hung Quan, Bill Bowser, Sherry Smith, and Danny Carden.  Finally, I must thank my graduate program advisor in the Department of Mathematics, Celia Davis, whose generous assistance alleviated many of the difficulties which came with four years of graduate school at a university 1100 miles away from my home.     \\

This work was supported in part by VIGRE Grant DMS0135345.


\newpage 
\begin{center}
\underline{\bf \Large Abstract} 
\end{center}

\medskip

\large

One may care to study operators and their inverses under conditions for which polynomial decay is too slow, and exponential decay is too fast or not exactly preserved.  In this paper I establish that a broad class of Banach *-algebras of infinite integral operators, defined by the property that the kernels of the elements of the algebras possess subexponential off-diagonal decay, is inverse closed in $\B(\Ltwo).$  This means that each subalgebra of $\B(\Ltwo)$ under consideration contains the inverse of each of its invertible elements.\\

A dividend of the above work concerning symmetry of Banach algebras is demonstrated en route to proving the above facts about decay of kernels of integral operator inverses.  A *-algebra is symmetric if the (algebraic) spectrum of positive elements is positive.  Historically it has often been quite difficult to verify that a given involutive algebra is symmetric, but using techniques discovered by Gr\"ochenig and Leinert I will show that the algebras under consideration are symmetric. \\

In the second part of this dissertation, I present the results of the IEEE Transactions on Communications paper written jointly by Thomas Strohmer and me.  We develop a comprehensive framework for the design of orthogonal frequency-division multiplexing (OFDM) systems, using techniques from Gabor frame theory, sphere-packing theory, representation theory, and the Heisenberg group.  We begin with a Gaussian pulse, which has optimal time-frequency localization (TFL), and orthogonalize it (using a clever method first discovered by the chemist L\"owdin) with respect to the time-frequency shifts of an OFDM system.  This process retains the optimal amount of the TFL of the Gaussian, vastly improving the robustness of the standard OFDM system to dispersion caused by multipath propagation and the Doppler effect.  Our next contribution to the OFDM algorithm is the introduction of Lattice-OFDM, wherein we admit more general lattices than the usual separable lattices.  This provides even more defense against dispersion, and  we present a theorem giving the OFDM system designer the ability to optimally adapt the lattice and pulse to the channel conditions.

\np
\pagestyle{myheadings} 
\pagenumbering{arabic}
\markright{  \rm \normalsize CHAPTER 1. \hspace{0.5cm}
 INTRODUCTION }
\large 
\chapter{Introduction}
\thispagestyle{myheadings}
\section{Banach algebra results}
This dissertation is a study of off-diagonal decay properties of kernels of invertible integral operators and their inverses.  Previous strides in the area have been made by Gel'fand, Gr\"ochenig, Leinert, Jaffard, and Gohberg, among others (see Section \ref{prevresults}).  I extend the results in \cite{GL03} on the off-diagonal decay of infinite matrices and their inverses, and the symmetry of infinite matrix algebras, to the appropriate class of distributions on $\R^{2d}.$

\begin{figure}[htb]
\begin{center}
\subfigure{
\epsfig{file=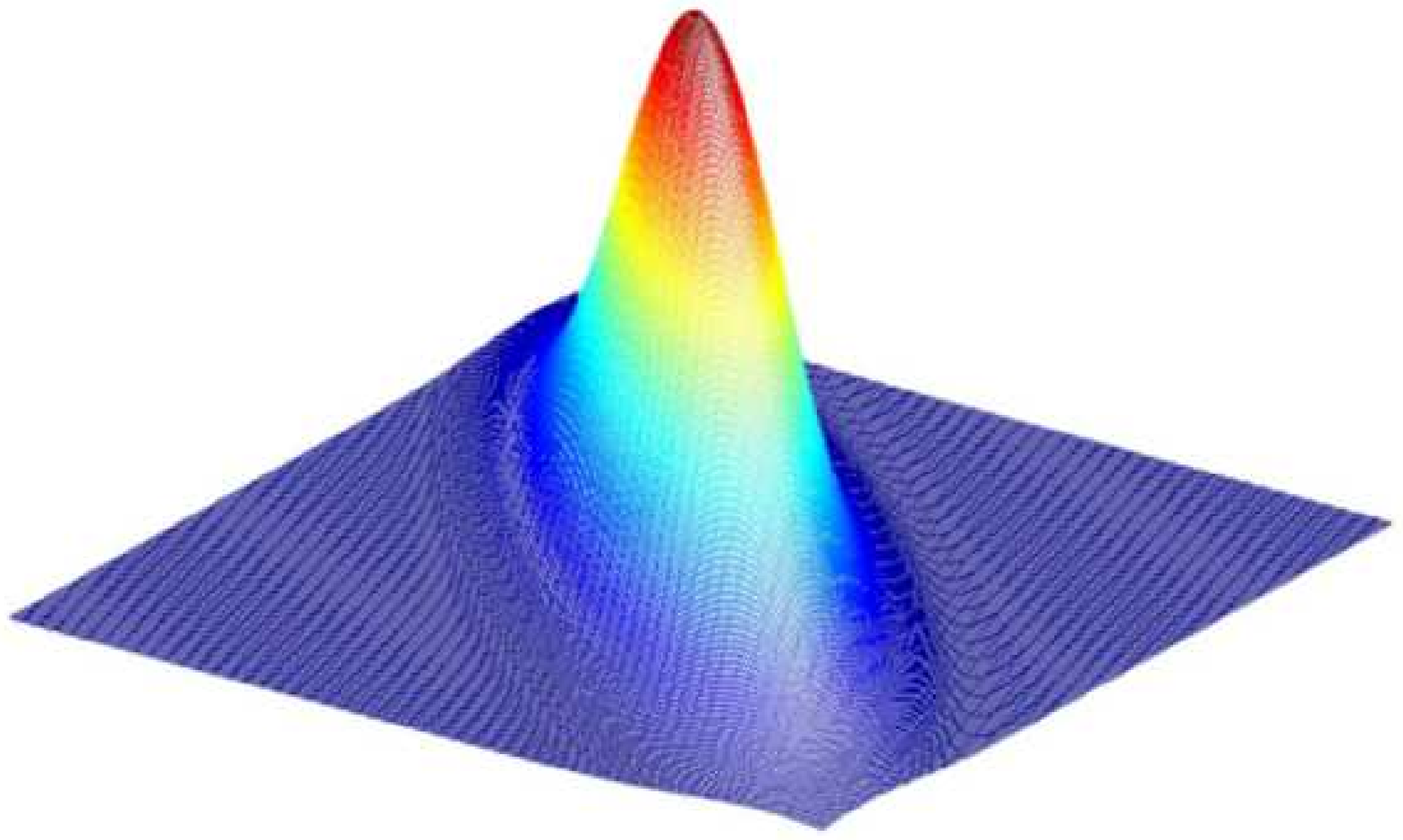,width=60mm}}
\subfigure{
\epsfig{file=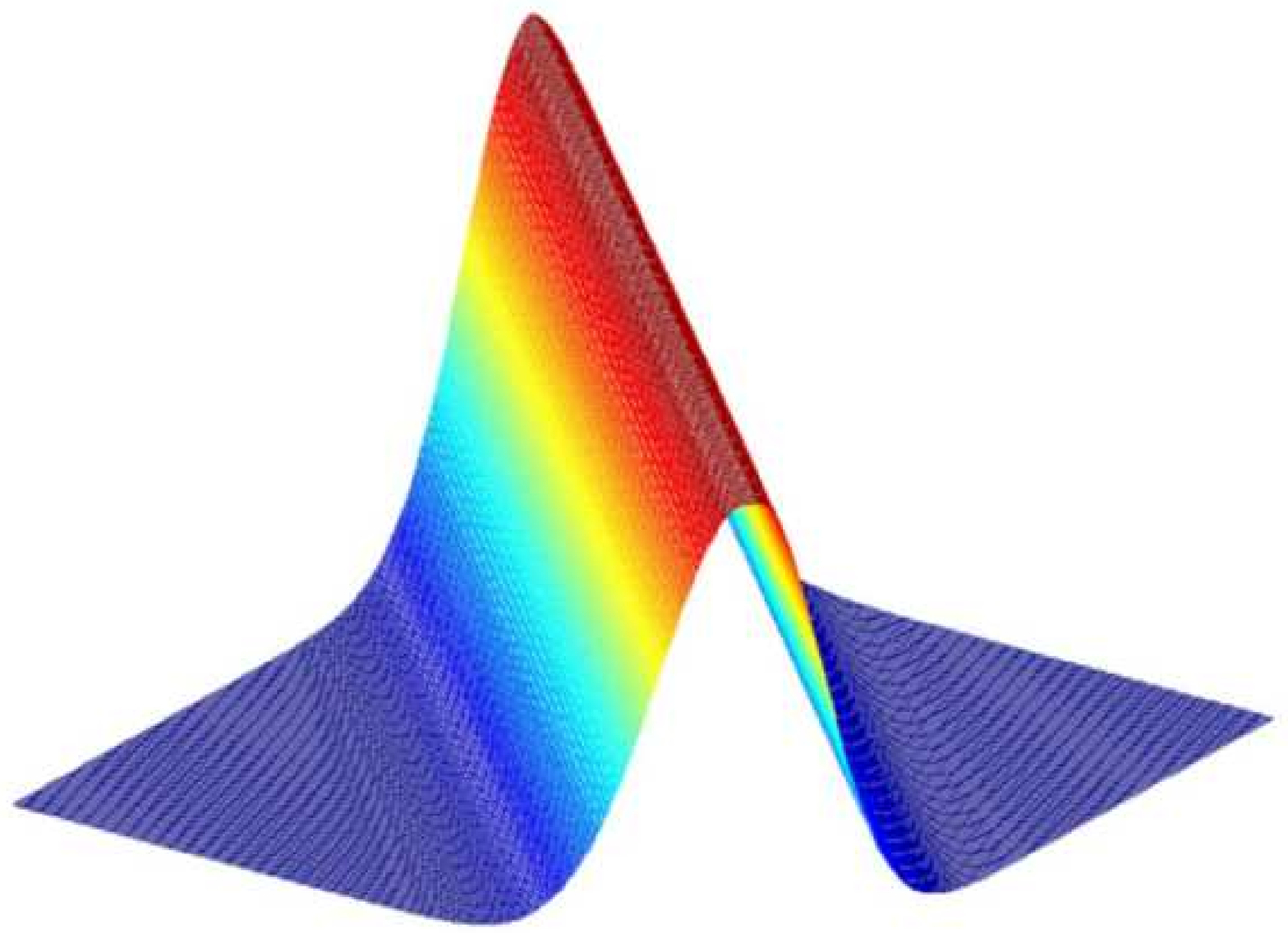,width=60mm}}
\caption{{Kernels, of integral operators, exhibiting off-diagonal decay}}
\end{center}
\end{figure}

One is often compelled to study operators and their inverses under conditions for which polynomial decay is too slow, and exponential decay is too fast or not exactly preserved \cite{GL03}.  My result establishes that a broad class of involutive Banach algebras (Banach *-algebras) of infinite integral operators, defined by the property that the kernels of the algebra elements possess subexponential off-diagonal decay, is inverse closed in $\B(\Ltwo).$  This means that each subalgebra of $\B(\Ltwo)$ under consideration contains the inverse of each of its elements.\\

An additional result, concerning symmetry of Banach algebras, is established en route while proving the above facts about decay of kernels of integral operator inverses.  An involutive algebra is {\em symmetric} if the spectrum of positive elements is positive.  Historically it has often been quite difficult to verify that a given involutive algebra is symmetric.  Losert recently solved a 30-year-old conjecture that the group algebra of a compactly generated, locally compact group of polynomial growth is symmetric.  The result of Gr\"ochenig and Leinert includes the establishment of symmetry for a large class of bi-infinite matrix algebras.  Extending their result about the inverse closedness of Banach algebras of bi-infinite matrices to inverse closedness of unital Banach algebras of infinite integral operators permits me to conclude that the algebras I consider are symmetric as well.  

\subsection{Applications:}

1.  Continuous-time blurring by convolution.  It is recognized that implementation will require discretization, quantization, and truncation, but a thorough understanding of the continuous theoretical model is often crucial to avoid setting up the model in a non-optimal way, of which there are examples from non-uniform sampling theory.\\

2.  Wireless channel scattering functions.  In practical terms a channel scattering function is often the sum of a delta ``function'' and a compact operator\footnote{This is the case when the channel is described by {\em Ricean fading}, where there is a line-of-sight component in the received signal.}, or the sum of a delta and a low-rank operator, each of which is invertible.\\

3.  Any integral equation whose kernel has appropriate decay properties but is not singular.

\subsection{Previous results in the field}\label{prevresults}
To illustrate the context of my main result, I briefly discuss previous results in the area, all of which regard operators on $\l^2(\Z^d)$:\\

1.  Commutative Cases:

\indent (a) {\em Wiener's Lemma} (Wiener, illustrated in e.g., \cite{RS00}) : $\l^1(\Z^d)$ is an algebra in which the multiplication is convolution.  For $a\in\l^1(\Z^d)$ consider the convolution operator $T_a$ acting on $\l^2(\Z^d),$ the action of which is given by $T_ac=a*c$.  Wiener's Lemma in this case is: if $a\in\l^1(\Z^d)$ and if $T_a$ is invertible as an operator on $\l^2(\Z^d),$ then $a$ is invertible in $\l^1(\Z^d)$ and hence the inverse operator is of the form $T_a^{-1}=T_b$ for some $b\in\l^1(\Z^d).$\\

\indent (b) {\em Weighted Wiener's Lemma} (implicitly Gel'fand, Raikov, Shilov \cite{GRS}, explicitly Strohmer \cite{Str00}):  Let $A=\{a_{k\l}\}$ be a Hermitian positive definite Laurent matrix, generated by the $0^{th}$ row denoted $a_k,$ with inverse $A^{-1}=\{b_{k\l}\},$ generated by the $0^{th}$ row denoted $b_k,$ and let $v(k)$ be a positive submultiplicative ($v(k+\l)\leq v(k)v(\l)$) weight function such that $\sum_{k\in\Z}|a_k|v(k)<\infty.$  If $v$ satisfies the condition
\be
\underset{n\rightarrow\infty}{\lim}\frac{1}{\sqrt[n]{v(-n)}}=1\qquad\mbox{and}\qquad \underset{n\rightarrow\infty}{\lim}\sqrt[n]{v(n)}=1\nn
\end{equation}
then $\sum_{k\in\Z}|b_k|v(k)<\infty.$\\

\indent 2.  Non-commutative cases:

\indent (c) {\em Non-commutative Wiener's Lemma} (Sj\"ostrand \cite{Sjo95}):  Let $C=C_{k,\l}\,,\,k,\l\in\Z^d,$ be a matrix with $|C_{k,\l}|\leq a(k-\l)$ for some $a\in\l^1(\Z^d),$ i.e. $C$ is dominated by the Laurent function generated by $a$ (see Section \ref{intops}).  If $C$ has a bounded inverse $D$ on $\l^2(\Z^d),$ then there exists $b\in\l^1(\Z^d)$ depending only on $a$ such that $|D_{k,\l}|\leq b(k-\l).$\\

\indent (d) {\em Exponential off-diagonal decay of infinite matrices} (Jaffard \cite{Jaf90}, Gohberg \cite{Goh90}):  Suppose $A$ is bounded on $\ll,$ with bounded inverse.  If either $a_{k\l}=0$ for $|k-\l|>M$ or if $|a_{k\l}|={\mathcal O}(e^{-\a|k-\l|}),\a>0,$ then $B=A^{-1}$ satisfies $|a_{k\l}|={\mathcal O}(e^{-\b|k-\l|})$ for some $\b>0.$\\

\indent (e) {\em Polynomial off-diagonal decay of infinite matrices} (Jaffard \cite{Jaf90}):  If $A$ is boundedly invertible on $\l^2(\Z^d)$ and $|a_{k\l}|={\mathcal O}(|k-\l|^{-s})$ for some $s>d,$ then $B=A^{-1}$ has the {\it same} polynomial-type off-diagonal decay $|b_{k\l}|={\mathcal O}(|k-\l|^{-s}).$\\

\indent (f) {\em Sub-exponential off-diagonal decay of infinite matrices} (Gr\"ochenig and Leinert \cite{GL03}):  This case includes (c).  Let $\odot$ denote the Hadamard product and let $\widetilde{v}_{k\l}:=v(|k-\l|).$  If $A$ is boundedly invertible on $\l^2(\Z^d)$ and the row- and column-sums of $|A\odot\,\widetilde{v}|$ are bounded, then $B=A^{-1}$ has the same off-diagonal decay as $A$ (i.e., the row- and column-sums of $|B\odot\,\widetilde{v}|$ are bounded).  Or, if $|a_{k\l}|\leq C\,v^{-1}(|k-\l|)$ for some $C,$ then $|b_{k\l}|\leq C'v^{-1}(|k-\l|).$\\

The contribution to the field provided by the first part of this dissertation is the extension of the results of Gr\"ochenig and Leinert to the continuous case.  That is to say, I will show that if $K$ is an integral operator on $\Ltwo,$ and the kernel $k$ of $K$ has sub-exponential off-diagonal decay measured either pointwise or additively, then the kernel of $K^{-1}$ has the same type of decay, respectively.

\section{Improving the standard OFDM communications scheme}In the second part of this dissertation, beginning with Chapter \ref{OBFDMinterpol}, we study orthogonal frequency-division multiplexing and introduce some extensions of the traditional OFDM scheme.  In this paper, we develop a comprehensive framework for the design of OFDM systems by introducing new pulse-shaping and lattice design methods.\\

Transmission over wireless channels can lead to two kinds of dispersion, namely intersymbol interference (ISI) and interchannel interference (ICI). ISI is caused (in part) by time dispersion due to multipath propagation and ICI results from frequency dispersion due to the Doppler effect.  Two considerations affect the dispersion robustness of a given OFDM system, namely the time-frequency (TF) localization of the pulse shape, and the time-frequency lattice from which the data are transmitted.\\

In Section \ref{contLowd}, we discuss pulse shaping.  We begin with the Gaussian due to its optimal time-frequency localization, and orthogonalize via a clever method first discovered by the chemist L\"owdin.  This yields a pulse which retains as much of the excellent TF localization of the Gaussian as possible under orthogonalization, a fact which we prove by appealing to Gabor analysis techniques.  Numerical experiments, the results of which are presented in section \ref{psofdmcpofdm}, illustrate the advantage of the L\"owdin-orthogonalized Gaussian over the standard square-wave in the presence of the Doppler effect. \\

OFDM is a special case of Biorthogonal Frequency-Division Multiplexing or BFDM.  In BFDM in general, one uses a different set of pulses for demodulation than the pulses used for transmission.  In Section \ref{ss:bfdm}, we briefly discuss BFDM, in which one sacrifices the orthogonality requirement but in return can be afforded improved time-frequency localization of the transmitter pulse (one can use a Gaussian, in fact).  This may improve robustness to multipath and Doppler dispersion, since the dispersion should affect only the transmitter pulse.  We then present an extension of the BFDM framework which we call extended BFDM, or EBFDM, which is a continuous interpolation between OFDM and BFDM.  Of all BFDM systems, OFDM, which uses the same pulses for transmission and demodulation, is optimally robust against additive white gaussian noise (AWGN).  EBFDM should have the attribute of permitting the communications system designer to fine-tune the balance between robustness to ICI/ISI and robustness to AWGN.  Inconclusive numerical studies are presented in Section \ref{pbofdmcpofdm}.\\

In Chapter \ref{GenPulseDesign}, we introduce Lattice-OFDM (LOFDM) in which the lattice parameters are based on sphere-packing considerations (thus maximizing spectral efficiency, or alternatively, maximizing robustness to ICI/ISI) and optimally adapted to channel conditions described by he channel scattering function.  We complete the pulse-shaping by applying a sequence of unitary operators, derived from the metaplectic representation of the symplectic group (see Appendix \ref{repSvNHg}), which optimally adapt the pulse to the lattice and channel.  Our framework affords the OFDM system designer an optimal approach to combatting dispersion.  Numerical results are presented to demonstrate the superiority of LOFDM in doubly dispersive channels. 

\np
\pagestyle{myheadings} 
\markright{  \rm \normalsize CHAPTER 2. \hspace{0.5cm}
 PRELIMINARIES AND NOTATION}
\large 
\chapter{Preliminaries and Notation}
\thispagestyle{myheadings}

\section{Basics}
Let $\Z,$ $\Q,$ $\R,$ and $\C$ denote the integers, the rational numbers, the real numbers and the complex numbers, respectively, and set $i=\sqrt{-1}.$  Let $y=a+bi\in\C;$ we denote by $\overline{y}$ the complex conjugate $\overline{y}=a-bi.$   \\

Denote by $\wt{\O}$ the set of all functions defined on $\R^{2d}.$  \\

``There exists'' will be deonted by the symbol $\,\exists,$ and ``for all'' by $\,\forall.$   The set with no elements, called the {\em empty set}, is denoted by $\emptyset.$\\

A {\em vector space} or {\em linear space} $X$ over a scalar field $\R$ is a set of elements on which are defined operations of elementwise addition $(+)$ and (commutative) scalar multiplication with the following properties:\\

\indent (a) the operation $+$ is commutative, i.e., for all $x,y\in X,$ we have $x+y=y+x;$\\
\indent (b) the operation $+$ is associative, i.e., for all $x,y,z\in X$ we have $(x+y)+z=x+(y+z);$\\
\indent (c) there is a zero vector such that for all $x\in X$ we have $x+0=x;$\\
\indent (d) for all $x\in$ there exists a unique element $-x$ such that $x+(-x)=0;$\\
\indent (e) since $\R$ is a field, it contains the element $1,$ for which $1x=x$ for all $x\in X;$\\
\indent (f) for all $x,y\in X$ and $\la,\mu\in\R$ we have $(\la+\mu)x=\la x+\mu x,\, (\la\mu)x=\la(\mu x),$ and $\la(x+y)=\la x+\la y.$\\

A {\em norm} on a linear space $X$ is a function $\|\cdot\|\rightarrow\R$ satisfying the properties\\

\indent (a) $\|x\|\geq 0$ for all $x\in X;$\\
\indent (b) $\|\la x\|=|\la|\|x\|$ for all $\la\in\C$ and $x\in X;$\\
\indent (c) $\|x+y\|\leq\|x\|+\|y\|$ for all $x,y\in X;$\\
\indent (d) $\|x\|=0\,\Leftrightarrow\,x=0.$\\

If a function $f$ is absolutely $p^{th}$-power integrable, we denote the $p$-norm of $f$ by $\|f\|_p,$ where 
\be
\|f\|_p=\left(\int_{\R^d}|f(x)|^p\,dx\right)^\frac{1}{p}.
\end{equation}
Note that for two functions $f_1$ and $f_2$ we can have $\|f_1-f_2\|_p=0$ even if there exists $x$ such that $f_1(x)\neq f_2(x).$  Thus is generated an equivalence class structure.  The space of equivalence classes of functions $f$ for which $\|f\|_p<\infty$ is denoted $\Lp.$\footnote{We shall follow the usual convention in denoting equivalence classes by $f,$ but referring to $f$ as a ``function.''}\\

The {\em essential supremum} of a nonnegative function $f,$ denoted $\underset{x\in\R^d}{\mbox{ess }\sup}\, f(x),$ is defined by
\be
\underset{x\in\R^d}{\mbox{ess }\sup}\,\,f(x)=\inf\,\bigl\{M\,|\,m\{x\,|\,f(x)>M\}=0\bigr\}
\end{equation}
where $m(\G)$ denotes the measure of a set $\G.$  Again, the essential supremum does not capture pointwise differences between two functions, and generates equivalence classes.  The space of equivalence classes of functions $f$ for which
\be
\|f\|_\infty :=\underset{x\in\R^d}{\mbox{ess }\sup}\,|f(x)|\,<\infty
\end{equation}
is denoted $\Linf.$  Not all interesting functions are in $\Lp$ or $\Linf.$ The space of infinitely differentiable functions on $\R^d$ is denoted $C^\infty(\R^d),$ and $C_0 :=\{f\,|\,f\,\text{ continuous and }\,f(x)\rightarrow 0 \,\text{ as }\,\|x\|\rightarrow\infty\}.$

\section{Banach and Hilbert spaces}\label{BanHilSp}

A normed linear space $(X,\|\cdot\|)$ is a linear space $X$ equipped with a norm $\|\cdot\|.$  A {\em Cauchy sequence} $\{x_n\}$ is a sequence in which for every $\e>0$ there exists an $N\in\N$ such that for all $n,m>N,$ we have $\|x_n-x_m\|<\e.$  A linear space $X$ is {\em complete} if every Cauchy sequence in $X$ converges to an element in $X.$  A {\em Banach space} is a complete normed linear space.  Familiar examples of Banach spaces are the $\Lp,$ with $1\leq p\leq\infty.$  
\begin{lemma}
Let $\B_1$ and $\B_2$ are Banach spaces with the same metric $d.$  Then $\B_1\cap\B_2$ is a Banach space.\\  
\begin{myproof} Let $\B_1,\B_2$ be Banach spaces with metric $d$ such that $\B_1\cap\B_2\neq\emptyset.$  We need to show that every Cauchy sequence in $\B_1\cap\B_2$ converges to an element $x$ in $\B_1\cap\B_2.$  Let $\{x_n\}_{n\in\N}\subset\B_1\cap\B_2.$  Then since $\{x_n\}\subset\B_1,\,\{x_n\}$ converges to an element $x$ in $\B_1.$  By virtue of its membership in $\B_2,\,\{x_n\}$ converges to an element $x'$ in $\B_2.$  If $d(x,x')=\e>0,$ then $\exists\,\e'$ and $N$ such that either $d(x_N,x)>\e'$ or $d(x_N,x)>\e'$ for all $n>N,$ contradicting the Cauchy property of $\{x_n\}.$  Thus $x'=x.$  We conclude that $\{x_n\}$ converges to a unique element in $\B_1\cap\B_2,$ and so $\B_1\cap\B_2$ is complete in the metric $d,$ and as such is a Banach space with norm inherited from $\B_1$ and $\B_2.$ 
\end{myproof}
\end{lemma}

An {\em inner product} on a linear space $X$ is a map $\langle\cdot,\cdot\rangle\,:\,X\times X\rightarrow\C$ with the properties that for all $x,y,z\in X$ and $\la,\mu\in\C,$ we have\\

\indent (a) $\langle x\,,\,\la y+\mu z\rangle=\la\langle x\,,\,y\rangle+\mu\langle x\,,\,z\rangle$\\
\indent (b) $\langle y\,,\,x\rangle=\overline{\langle x\,,\,y\rangle}$\\
\indent (c) $\langle x\,,\,x\rangle\geq 0$\\
\indent (d) $\langle x\,,\,x\rangle=0\Leftrightarrow x=0.$\\

\noindent We call a linear space equipped with an inner product an inner product space.  If $X$ is an inner product space, we can define a norm on $X$ by $\|x\|=\sqrt{\langle x\,,\,x\rangle}.$  A {\em Hilbert space} is a inner product space which is complete under the norm derived from the inner product.  Famous examples of Hilbert spaces are $\Ltwo$ and $\ell^2(\Z^d),$ the space of square-summable sequences of complex $d$-tuples.

\section{Algebras}

An {\it algebra} is a vector space $\A$ on which the following properties are satisfied:
\bea
&&a(bc)=(ab)c\label{algebra1}\\
&&a(b+c)=ab+ac\label{algebra2}\\
&&(a+b)c=ac+bc\label{algebra3}\\
&&\la(ab)=(\la a)b=a(\la b)\label{algebra4}
\end{eqnarray}
for each $a,b,c\in\A$ and scalar $\la.$\\

An {\em involutive algebra}, is an algebra $\A$ on which there exists a unary operation (*) such that
\bea
a^{**}&=&a\\
(a+b)^*&=&a^*+b^*\\
(\la a)^*&=&\la^*a^*\\
(ab)^*&=&b^*a^*,
\end{eqnarray}
for all $a\in\A.$\\

We call $\A$ a {\em unital algebra} if $\A$ contains the multiplicative identity.\\

A {\it normed algebra} is an algebra $\A$ on which is defined a norm $\|\cdot\|$ which is {\it submultiplicative}, i.e.,
\be
\|ab\|\leq\|a\|\|b\|\nn
\end{equation}
for all $a,b\in\A.$\\

A {\it Banach algebra} is a normed algebra which is also a Banach space, i.e., complete in the norm $\|\cdot\|.$\\

\section{Fourier considerations}
The Fourier transform of $f\in\Lone$ denoted by $\hat{f}$ is 
\be
\hat{f}(\w) = \int_{\R^d}f(t) e^{-2\pi i \w\cdot t} dt.
\label{fourier}
\end{equation}
We also write ${\cal F}f$ instead of ${\hat f}.$  If $\hat{f}\in\Lone,$ then we have the inversion formula
\be
f(t) = \int_{\R^d}\hat{f}(\w) e^{-2\pi i \w\cdot t} d\w.
\label{fourierinv}
\end{equation}
The Fourier transform can be extended to a unitary operator on $\Ltwo;$ from here on the notation $\F$ shall represent this unitary operator where applicable.\\

Fix a function $g\in \Lq,\,q\geq 1,\,g\neq 0,$ and let $p=\frac{q}{q-1}.$  Then the {\em short-time Fourier transform} (STFT) of a function $f\in \Lp$ is given by
\be
V_gf(x,\w)=\int_{\R^d}f(t)\overline{g(t-x)}e^{-2\pi i \w\cdot t}\,dt,\quad x,\w\in\R^d.\label{stft}
\end{equation}

The {\em ambiguity function} of $f \in \Ltwo$ is defined to be
\begin{equation}
\label{ambi}
Af(x,\omega) = \int_{\R^d}
f\big(t+\frac{x}{2}\big) \overline{f\big(t-\frac{x}{2}\big)} 
e^{-2\pi i \w\cdot t} dt,
\end{equation}
where $\overline{f(t)}$ denotes the conjugate of $f(t)$.
Let $p\geq 1\,,\,q=\frac{p}{p-1}.$  The {\em cross-ambiguity function} of $f\in \Lp, g\in \Lq$ is
\begin{equation}
\label{crossambi}
A(f,g)(x,\omega) = \int_{\R^d}
f\big(t+\frac{x}{2}\big) \overline{g\big(t-\frac{x}{2}\big)}
e^{-2\pi i \w\cdot t} dt = e^{\pi i \w\cdot x}V_gf(x,\w).
\end{equation}

\section{Tempered distributions}\label{tempereddistributions}
Let $\Z_+=\{n\in\Z\,\,|\,\,n\geq 0\}$ denote the nonnegative integers.  A {\em multi-index} $\a=(\a_1,\ldots,\a_d)\in\Z^d_+$ is a $d$-tuple of nonnegative integers.  Let $x=(x_1,\ldots,x_d)\in\R^d$ and $\a\in\Z^d_+,$ and define
\bea
D^\a&=&\left(\frac{\partial}{\partial x_1}\right)^{\a_1}\cdots\left(\frac{\partial}{\partial x_d}\right)^{\a_d},\nn\\
x^\a&=&\prod^d_{j=1}x^{\a_j}_j\nn\\
|x|&=&\sqrt{x^2_1+\cdots+x^2_d}\label{normabs}
\end{eqnarray}
We will use the notation in \eqref{normabs} for the finite-dimensional 2-norm throughout the sequel.  The {\em Schwartz space} $\mathcal{S}$ consists of all $\ph\in\C^\infty(\R^d)$ such that $\underset{x\in\R^d}{\sup}\left|x^\a D^{\,\b}\ph(x)\right|$ is finite for every pair of multi-indices $\a,\b\in\Z^d_+.$  Elements of $\mathcal{S}$ thus decay faster than the inverse of any polynomial as $|x|\rightarrow\infty.$\\

The space of continuous linear functionals on $\mathcal{S},$ or the topological dual space of $\mathcal{S},$ is denoted $\mathcal{S}^*.$  Elements of $\mathcal{S}^*$ are called {\em tempered distributions}.  The flagship example of a tempered distribution is the Dirac delta ``function'' defined by
\be
\d(\ph)=\ph(0).\nn
\end{equation}

The most intuitive way to view $\d$ is as a measure on $\mathcal{S},$ whereby
\be
\int_{\R^d}\ph(x)\d(x)\,dx=\ph(0).\nn
\end{equation}
We can extend the domain of the point measure $\d(x)$ to any function defined on $\R^d$ (including the constant function $1$) enabling us to say $\int_{\R^d}\d(x)\,dx=1.$  Note that in the sense of Lebesgue, $\int_{\R^d}\d(x)\,dx=0.$  The Fourier transform of a tempered distribution $\ph$ can be defined by $\widehat{\ph}(f)=\ph(\tilde{f})$ where $\tilde{f}=\F^{-1}(f), f\in\S.$  Thus 
\be
\widehat{\d}(f)=\int_{\R^d}\d(x)\tilde{f}(x)\,dx=\tilde{f}(0)=\int_{\R^d}\frac{1}{\sqrt{2\pi}}f(\w)\,d\w\label{dist1}
\end{equation}
which implies that, in the sense of measure, $\widehat{\d}(\w)$ can be thought of as the constant $\frac{1}{\sqrt{2\pi}}.$  Thus $\widehat{\d}(\w)$ is a multiple of the usual Lebesgue measure, and $\d(x)$ can be thought of as the inverse Fourier transform of the Lebesgue measure, up to a constant.

\section{Integral operators}\label{intops}
Let $h$ be a member of some function space $X$ on $\R^d.$  An {\em integral operator} on $X$ is an operator $K$ defined by
\be
K(h)(x):=\int_{\R^d}k(x,y)\,h(y)\,dy\label{integraloperator}
\end{equation}
when $Kx\in X,$ i.e. $K:X\rightarrow X.$  We call $k(x,y),\,x,y\in\R^d$ the {\em kernel} of $K.$  In this paper we shall remove the distinction between $K$ and $k$ since $k$ defines $K$ uniquely, and will refer to the kernel and the operator by the same symbol.  That is, in the sequel $k(h)=K(h).$\\

Let $f:\R^d\rightarrow\R$ be any function.  The Laurent function on $\R^{2d}$ induced by $f$ is $f_L(x,y)=f(x-y).$  We say that $f$ {\em generates} $f_L.$  This is the bi-infinite analog of Toeplitz functions on $\R^{2d}_+,$ and the idea of a Laurent function is of course well-defined in discrete settings as well.  The Laurent function $f_L$ can be viewed as the kernel of a convolution operator, provided that the convolution is defined.  Let $g$ be a function defined on $\R^d.$  Then the convolution of $f$ with $g$ is defined to be
\be
f*g(t):=\int_{\R^d}f(x)g(t-x)\,dx=\int_{\R^d}f(t-x)g(x)\,dx,\label{convdef}
\end{equation}
provided the integral in \eqref{convdef} is defined.  If we view $f_L$ as an operator, we have $f*g=f_L(g).$

\section{Continuous analog of matrix multiplication}\label{contanalogmatmult}

We write explicitly the integral operator formed by composition of integral operators.  Let $1\leq p,q\leq\infty,\,\frac{1}{p}+\frac{1}{q}\!=\!1,\,f\!=\!f(t,y),\,g\!=\!g(x,t),$ and define the restriction operator $R$ as $R_{\Z^{2d}}(h)(x,y):=h|_{\Z^{2d}}(x,y).$  We further demand that \\

(1) each ``row'' in the matrix $A:=R_{\Z^{2d}}(f)(t,y)$ is absolutely $p$-summable, and \\

(2) each ``column'' in the matrix $B:=R_{\Z^{2d}}(g)(x,t)$ is absolutely $q$-summable \\

(making $R_{\Z^{2d}}(f)$ and $R_{\Z^{2d}}(g)$ bounded linear operators on $\lq$ and $\lp,$ respectively).  For $A,$ let the row index ($y$-axis for $f$) be $j$ and the column index ($t$-axis for $f$) be $\l,$ and for $B$ let the row index ($t$-axis for $g$) be $\l$ and the column index ($x$-axis for $g$) be $k.$  We have 
\be
(AB)_{jk}=\underset{\l\in\Z}{\sum}A_{j\l}B_{\l k}.
\end{equation}

It is clear that the analogous continuous multiplication, defined for the moment on the space of all functions $g:\R^{2d}\rightarrow\R$ (and for the moment ignoring convergence considerations), denoted as $\tilde{\O}(\R^{2d})$ which respects axes and orientation must be
\be
(f\star g)\,(x,y):=\int_{\R^d} f(t,y)g(x,t)\,dt\label{fstarg},
\end{equation}
where we shall assume that $x=0$ if $g\in\tilde{\O}(\R^d).$  To visualize the difference in indexing between this operation and matrix multiplication, consider $g\in\tilde{\O}(\R).$  In this case, the convention should be that $g$ is visualized vertically, with positive $y$-values at the top and negative $y$-values at the bottom.  Of course the operation $\star$ satisfies
\be
(f\star g)\star h=f\star (g\star h),\label{assoc}
\end{equation}
when both $f\star g$ and $g\star h$ exist.  A proof of associativity is given in Appendix \ref{starassoc}.\\

We shall later show that certain classes of tempered distributions in $\S^*$ act on the Hilbert space $\Ltwo.$  In this case, we have, for $h_1,\,h_2\in\Ltwo,$
\bea
\langle h_1,f\,h_2\rangle&=&\int_{\R^d}\int_{\R^d}h_1(y)\overline{f(t,y)\,h_2(t)}\,dt\,dy\nn\\
&=&\int_{\R^d}\int_{\R^d}h_1(y)\overline{f((t,y)}\,dy\,\overline{h_2(t)}\,dt\nn\\
&=&\int_{\R^d}\int_{\R^d}\overline{f(y,t)}\,h_1(y)\,dy\,\overline{h_2(t)}\,dt\nn\\
&=&\langle f^*h_1,\,h_2\rangle\nn
\end{eqnarray}
which implies that $f^*=(\overline{f})^T.$  Thus, the adjoint operation in any subspace of $\wt{\O}$ which acts on $\Ltwo$ is given by reflection about the line $t=y$ followed by conjugation, i.e.,
\be 
f^*(t,y)=\overline{f(y,t)}.\label{invol}
\end{equation}
In addition, we have (abusing notation slightly) $(f\star g)^*=g^*\star f^*,$ from
\bea
(f\star g)^*&=&\left(\int_{\R^d}f(t,y)g(x,t)\,dt\right)^*\nn\\
&=&\left(\overline{\int_{\R^d}f(t,y)g(x,t)\,dt}\right)^T\nn\\
&=&\left(\int_{\R^d}\overline{f(t,y)g(x,t)}\,dt\right)^T\nn\\
&=&\int_{\R^d}\overline{g(t,x)}\,\,\overline{f(y,t)}\,dt.\label{fgaster}
\end{eqnarray}

\section{$\star$-invertibility}\label{starinvert}

The left unit for the multiplication operation $\star$ is $\d(t-y);$ the right unit is $\d(x-t);$ clearly these are equal.  But $\star$ need not be a group operation.  By illustrating the action of a convolution operator $f_L$ on a function in $\Ltwo$ which ``oscillates to death'' in the sense of Lieb and Loss \cite{LiebLoss}, we will show that $f_L$ is not boundedly invertible.  Consider the Gaussian $f(t)=e^{-\pi |t|^2}\!,\,t\in\R^d,$ and form the Laurent operator $f_L$ on $\Ltwo$ with 
\be
f_L(t,y)=e^{-\pi|y-t|^2}.
\end{equation}
Note that $f_L$ is not Hilbert-Schmidt (the kernel $f_L\not\in L^2(\R^{2d});$ Hilbert-Schmidt operators are compact and hence not invertible).  However, we do have $\underset{\|\,g\|_2=1}{\inf}\|f_L\star g\,\|_2=0,$ which precludes invertibility by a bounded linear operator.  Indeed, consider $g_\w=(2^{\frac{d}{4}})f(t)e^{2\pi i\w\cdot t}$ which has $\|g_\w\|_2=1,\forall\,\w\in\R^d.$  We have 
\bea 
f_L\star g_\w(y)&=&\int_{\R^d}e^{-\pi|y-t|^2}e^{-\pi|t|^2}e^{2\pi i\w\cdot t}\,dt\nn\\
&=&\F(e^{-\pi|y-t|^2}e^{-\pi|t|^2})(y,\w)\nn
\end{eqnarray}
which is just the STFT of $f$ with $f$ as its own window.  Since $\|V_ff\|_{L^2(\R^{2d})}=\|f\|^2_{\Ltwo}<\infty,$ we infer that the slices $\|V_ff(y,\w=\w_0)\|_{\Ltwo}\rightarrow 0$ as $|\w_0|\rightarrow\infty$ which implies $\|f_L\star g_\w\|_\Ltwo\rightarrow 0$ as $|\w|\rightarrow\infty$ (this result is related to the Riemann-Lebesgue lemma).  Proposition 4.8 in \cite{Doug98} states that if $T$ is a bounded operator on a Hilbert space, then $T$ is (boundedly) invertible on the Hilbert space $\LtR$ if and only if $T$ is bounded below and has dense range; thus we have shown that $f_L$ is not invertible.\\

In the sequel we will drop the symbol $\star$ and simply denote $f\star g$ by $fg$ when the context of the multiplication is clear.\\

We next consider the existence of invertible integral operators, specifically invertible {\em convolution} operators.  Let $f,g$ be functions in function spaces yet to be determined.  We define the convolution of $f$ with $g$ to be
\be
(f*g)(t)=\int_{\R^d}f(x)g(t-x)\,dx .
\end{equation}
We wish to invert this operation on $g.$  Young's inequality says that if $f,g\in L^1(\R^{d})$ then $f*g\in L^1(\R^{d}).$  In this case we have 
\be 
\widehat{f*g}=\widehat{f}\,\cdot\,\widehat{g}
\end{equation}
and thus 
\be
g=\F^{-1}\biggl(\frac{\widehat{f*g}}{\widehat{f}}\biggr)=\F^{-1}\biggl(\frac{\widehat{f}\cdot\widehat{g}}{\widehat{f}}\biggr)=\F^{-1}(\widehat{g})\label{finvghat}
\end{equation}
showing that the convolution operation on $g$ is invertible (since $\widehat{f}\in C_0(\R^{d}),$ we can use a limiting process at the points where $\widehat{f}=0$).\\

But what about the operator $f^{-1}$ in general, operating on functions not necessarily formed by convolution?  When is it a bounded operator on, say, $\Lone$?$\,$\footnote{We could ask when $f^{-1}$ is a bounded operator on the Hilbert space $\Ltwo$ but for the purposes of the present argument, we will not complicate matters by appealing to the limit arguments necessary for the treatment of $\F$ and $\F^{-1}$ in $\Ltwo.$}  
By inspection of \eqref{finvghat} we have
\be
f^{-1}\,h=\F^{-1}\left(\,\widehat{h}/\widehat{f}\,\right).\nn
\end{equation}
If $\underset{\w}{\inf\,\,}|\widehat{f}(\w)|>0,$
\be
\frac{\widehat{h}(\w)}{\widehat{f}(\w)}\leq\frac{\widehat{h}(\w)}{\underset{\w}{\inf\,\,}\widehat{f}(\w)}=C\,\widehat{h}(\w)
\end{equation}
which is in the domain of $\F^{-1}$ if $\widehat{h}\in\Lone.$  Obviously if $\underset{\w\rightarrow+\infty}{\lim} \widehat{f}(\w)=0$ or $\underset{\w\rightarrow-\infty}{\lim} \widehat{f}(\w)=0,$ then $\frac{\widehat{h}(\w)}{\widehat{f}(\w)}$ is not bounded, hence not integrable and not in the domain of $\F^{-1}.$  It is also possible that, for a given $h,\,\,\frac{\widehat{h}(\w)}{\widehat{f}(\w)}\in\Lone$ with $\widehat{f}(\w)=0$ at a countable number of points $\w_n.$  Consider, e.g., an $\widehat{f}$ which is differentiable, and whose slope is sufficiently large, and sufficiently increasing as $n\rightarrow\pm\infty,$ at the zero crossings of $\widehat{f}(\w)$ so as to ensure that $\frac{\widehat{h}(\w)}{\widehat{f}(\w)}\in\Lone.$  For example, suppose that (loosely speaking)  away from the zero crossings, $\widehat{f}=1,$ and near the zero crossings $\frac{\widehat{h}(\w)}{\widehat{f}(\w)}$ had the same type of vertical decay behavior as $(x-\w_n)^{-\frac{1}{4n^2}}$ (for any $\e>0,\,\,\int^{\w_n+\e}_{\w_n-\e}(x-\w_n)^{-\frac{1}{4n^2}}\,dx$ exists and as $n\rightarrow\infty$ is decreasing sufficiently rapidly so that the sum of all the integrals about the countable number of singularities is finite).  But no such $f$ will satisfy $\frac{\widehat{h}(\w)}{\widehat{f}(\w)}\in\Lone$ for all $h\in\Lone.$\\

So what kind of function is $f$ with $\underset{\w}{\inf\,\,}\widehat{f}(\w)>0$?  The Riemann-Lebesgue lemma states that $h\in\Lone\Rightarrow\widehat{h}\in C_0(\R^d).$  Thus neither $f$ nor $\widehat{f}$ can be an $\Lone$ function.  Could $f$ be a member of another $\Lp$ space?  Certainly not.  Let $C=\underset{\w}{\inf\,\,}\widehat{f}(\w).$  Then $\F^{-1}\widehat{f}=\frac{C}{\sqrt{2\pi}}\,\d+f_2,$ where $f_2$ may be in some $\Lp$ space but $\d\not\in\Lp$ for any $p$ since $\d^{ ^{\scriptstyle p}}=\d$ is not a function in the Lebesgue sense (see Section \ref{tempereddistributions}).\\

Thus for a convolution operator to be invertible, it is necessary that the kernel of the operator possess as a summand a Dirac delta distribution.  I.e., 
\be
f_L(x,y)=f(x-y)\qquad\mbox{with}\qquad f=\tilde{f}+\d.\label{distcond}
\end{equation}
We will call $f$ the {\em generator} of the operator $f_L.$  In the discrete setting, the analogous condition is not sufficient; this is called the {\em Wiener-Pitt phenomenon}.  For any locally compact abelian group $G,$ define the {\em Wiener algebra} $W(G)$ as 
\be
W(G)=\left\{g\in\Lone\left|\right.\hat{g}\in\Lone\right\}.\label{wieneralg}
\end{equation}
In our setting, it is known that if the nondistributional part $\tilde{f}$ of the generator $f$ of the Laurent operator is in the Wiener algebra $W(G),$ then the Wiener-Pitt phenomenon does not occur and condition \eqref{distcond} is sufficient as well as necessary.  We assume in the sequel that if an operator $g_L$ under consideration is a convolution operator, its generator $g$ has nondistributional part $\tilde{g}\in W(\R^d)\subset\Lone.$  We can then safely say that $g_L$ acts on $\Ltwo.$

\section{Spaces of distributions}
Recall that for a convolution operator $g_L$ to be invertible, we must have $g_L(x,y)=g(x-y)+\d(x-y),$ where $g:\R^d\rightarrow\R.$  We extend that idea now to the kernel of any integral operator.  Formally adjoin to $\widetilde{\O}$ the unit element $e:=\d(x-y)$ to yield the class $\O,$ whose elements are of the form $f=\tilde{f}+\mu e,\,\tilde{f}\in\widetilde{\O}\,,\,\mu\in\R.$  We shall consider the class $\O_1,$ whose elements are of the form $f=\tilde{f}+\mu e+\sum_{k=1}^\infty\,c_k\d(x_k,y_k),\,\tilde{f}\in\widetilde{\O}\,,\,\mu\,,c_k\in\R,$ and where each hyperplane defined by constant $x$ or $y$ supports only a finite number of deltas.

\np
\pagestyle{myheadings}
\chapter{Classes of Integral Operators Classified by Off-Diagonal Decay of the Kernel}\label{decaydistributions}
\thispagestyle{myheadings} 
\markright{  \rm \normalsize CHAPTER 3. \hspace{0.5cm}
  INTEGRAL OPERATORS CLASSIFIED BY KERNEL DECAY }
In this chapter we will show that if the kernel of an integral operator possesses subexponential decay, whether measured by weighted ``pointwise'' essential supremum boundedness, or by essential supremum boundedness of the weighted partial integral, then the kernel of the inverse integral operator possesses the same type of decay, measured in the same manner as the original kernel.
\section{Weight functions}\label{weights}
We define the relevant classes of distributions under consideration by their decay relative to weight functions.  We measure the decay of an integral operator kernel in two ways: by boundedness of its weighted integral or pointwise essential supremum boundedness against weights.  Weights are used in order to quantify the convergence properties.  In addition, the quantitative nature of measurement against particular weights will be of some advantage when deriving numerical methods and assessing their performance.
\begin{definition}
An {\it admissible} weight $v$ is 
\indent\indent (a) {\it subexponential}: 
\begin{equation}
v(x)=e^{\k(|x|)}\,\,\mbox{with}\,\,\k \mbox{ concave and }\k(0)=0,\quad\mbox{and}\label{subexpdef}
\end{equation}
\indent\indent (b) satisfies the Gel'fand-Raikov-Shilov (GRS)-condition (cf. \cite{GRS})
\be
\underset{n\rightarrow\infty}{\lim}\,v(nx)^{1/n}\,=\,1 \quad\forall x\in\R^d.\label{grs}
\end{equation}
\end{definition}

We shall discuss only admissible weights, which lie in $L^p_{loc}(\R^d)\,,\,1\leq p\leq\infty,$ in what follows.  An example is $\t_s(x)=(1+|x|)^s,\,s>2d.$  But to further clarify what sort of function is an admissible weight, we state some of the implied properties of such weights.   Note that since $\k$ is concave and $\k(0)=0,\, v$ is {\it submultiplicative}, i.e.,
\be
v(x+y)\leq v(x)v(y).\label{submult}
\end{equation}
Although $\k(t), t\in\R$ is concave $v(x)=e^{\k(|x|)}$ is neither log-concave nor log-convex.  In fact, if $x\in\R^2,$ then every point on the graph of $v(x)$ resembles a saddle point (see Figure \ref{fig:notconcconv} in Appendix \ref{weights1}).  As a consequence of the concavity of $v$ and the fact that $\k$ is nonnegative, $\k$ must be nondecreasing (increasing if $v$ is strictly concave).  Finally, we note that since $\k$ is concave, $v$ must be continuous.  For a more detailed analysis of each of these properties and an illustration of the GRS-condition, see Appendices \ref{weights1} and \ref{grsappendix}, respectively. \\

\section{The distribution classes $\Ll^1_v$ and $\Ll_v$}
We define the relevant classes of distributions.
\begin{definition} The class ${\mathcal L}_{v}^{1}$ consists of all elements $f\in\O_1$ such that
\bea
S_{x,v}(\tilde{f})&:=&\underset{x\in\R^d}{\mbox{ess }\sup}\,\int_{\R^d}\,|\tilde{f}(x,y)|\,v(x\!-\!y)\,dy<\infty,\quad  \mbox{ and }\nn\\
S_{y,v}(\tilde{f})&:=&\underset{y\in\R^d}{\mbox{ess }\sup}\,\int_{\R^d}\,|\tilde{f}(x,y)|\,v(x\!-\!y)\,dx<\infty. \nn
\end{eqnarray}
Define a norm on $\Ll^1_v$ as
\be
\|f\|_{{\mathcal L}_{v}^{1}}:=\max \bigl\{S_{x,v}(\tilde{f})\,,\, S_{y,v}(\tilde{f})\bigr\}+|\mu|=\|\tilde{f}\|_{\Ll^1_v}+|\mu|.\label{defineLv1}
\end{equation}
\end{definition}
Write ${\mathcal L}^1,S_x(f),S_y(f)$ in the case of $v\equiv 1.$  Note that $f\in\Ll^1_v$ is quite similar to a Hille-Tamarkin operator (cf. e.g., \cite{Lab01}), the difference being that integrability of members of $\Ll^1_v$ is measured against weights, and that elements of $\Ll^1_v$ may contain the identity and/or point measures.  Alternatively, $\Ll^1_v$ can be (nearly) viewed in the context of the $L^P(\R^{2d})$ spaces of mixed-norm (cf.\cite{BePa61}).\footnote{Note the capital $P,$ generally signifying $P=\{p_i\},i=1,\ldots n$ with $\sum_{i=1}^np_i=1;$ in our case $n=2.$}  Obviously the above relationships are stronger when $v\equiv 1.$
\begin{definition}
The class ${\mathcal L}_{v}$ consists of all elements $f\in\O_1$ such that
\be
\underset{x,y\in\R^d}{\mbox{ess }\sup}\,|\tilde{f}(x,y)|\,v(x\!-\!y)<\infty.\nn
\end{equation}
Define a norm on $\Ll_v$ as
\be
\|f\|_{\Ll_v}:=\underset{x,y\in\R^d}{\mbox{ess }\sup}\,|\tilde{f}(x,y)|\,v(x\!-\!y)+|\mu|=\|\tilde{f}\|_{\Ll_v}+|\mu|.\label{defineLv}
\end{equation}
\end{definition}
\begin{figure}[htb]
\begin{center}
\subfigure{
\epsfig{file=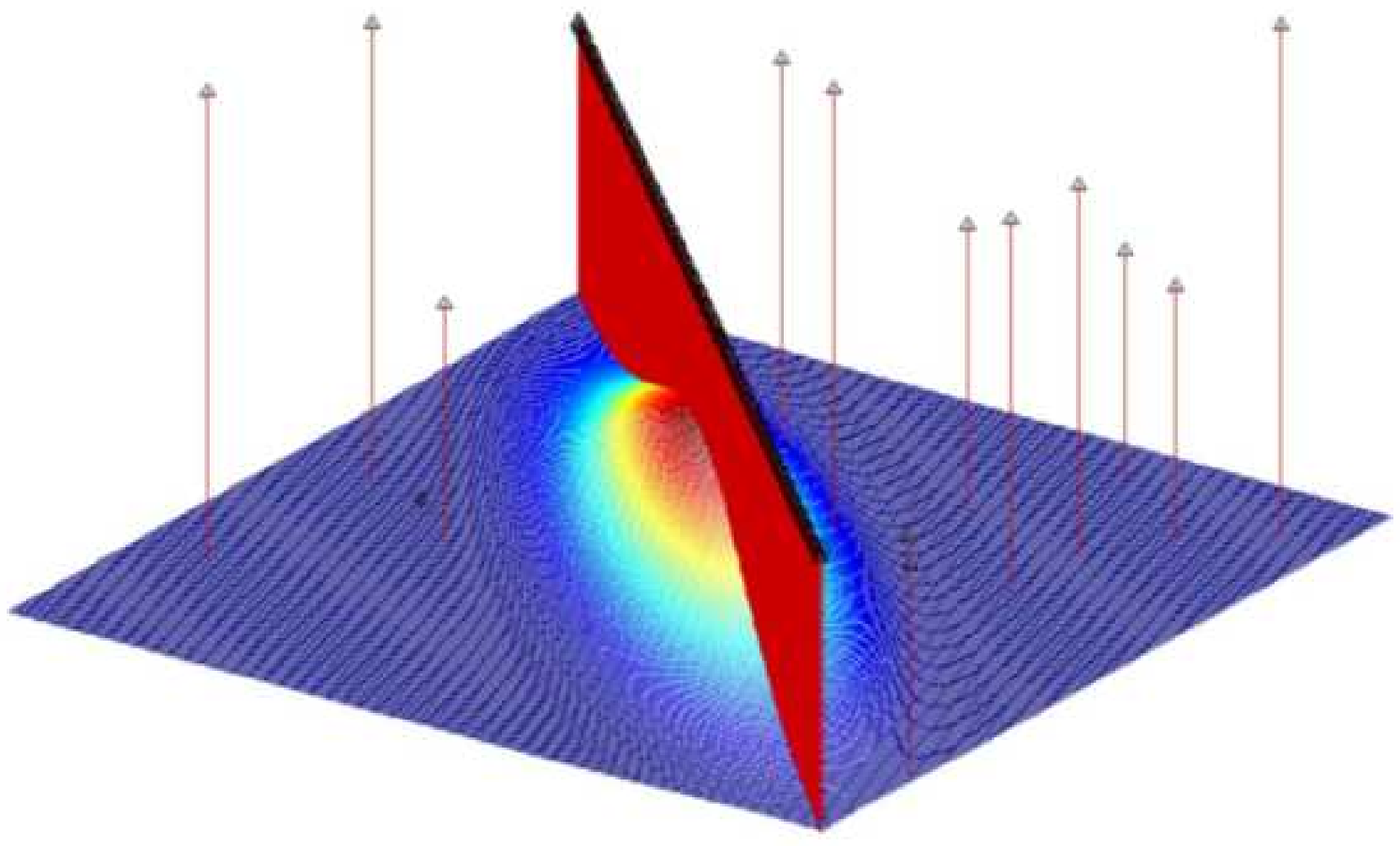,width=70mm}}
\subfigure{
\epsfig{file=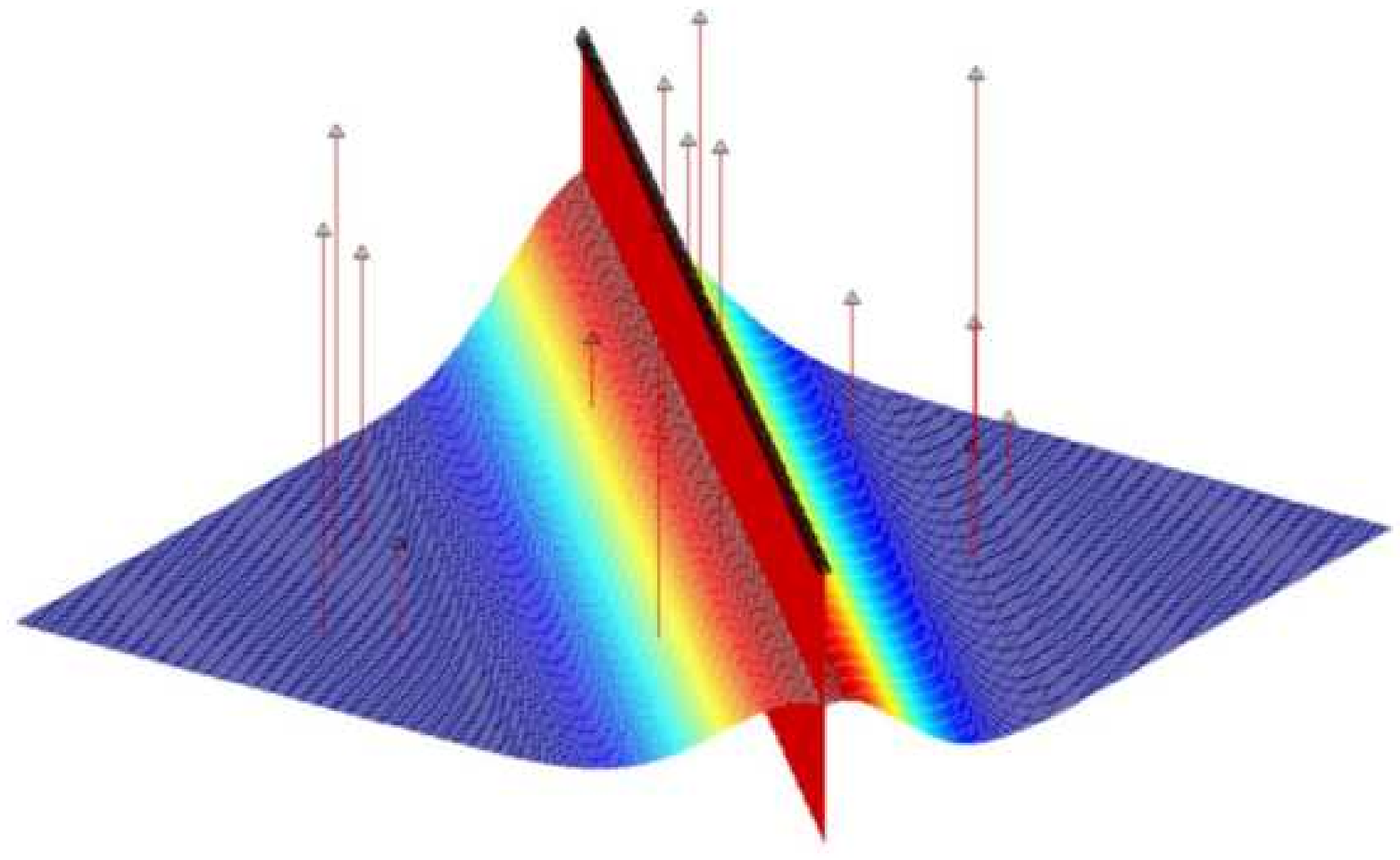,width=70mm}}
\caption{{Depictions of truncated elements of $\Ll^1_v$ or $\Ll_v,$ derived from compact and noncompact operators}}
\end{center}
\end{figure}

In terms of Hadamard multiplication, if we let $v_L(x,y)=v(x-y),$ we have 
\be
\Ll^1_v=\left\{f\in\O_1\,\Bigl|\Bigr.\,\underset{y\in\R^d}{\mbox{ess }\sup}\int_{\R^d}\!(\tilde{f}\odot v_L)\,dx\!<\!\infty \mbox{  and }\underset{x\in\R^d}{\mbox{ess }\sup}\int_{\R^d}\!(\tilde{f}\odot v_L)\,dy\!<\!\infty\right\}\nn
\end{equation}
and
\be
\Ll_v=\left\{f\in\O_1\,\Bigl|\Bigr.\,\underset{(x,y)\in\R^{2d}}{\mbox{ess }\sup}\,(\tilde{f}\odot v_L)(x,y)<\infty\right\}.\nn
\end{equation}

The inclusion of a infinite number of delta-distributions, with only a finite number supported on any given hyperplane defined by constant $x$ or $y,$ in members of $\Ll^1_v$ and $\Ll_v$ is permitted since the deltas make only finite contributions to the integrals on sets of measure zero, thus preserving the equivalence class structure of $\Ltwo$ and inducing equivalence class structures on both $\Ll^1_v$ and $\Ll_v.$  We shall adopt the same convention as one does for the $\Lp$-spaces, namely that although $f$ in, say, $\Ll^1_v$ represents an equivalence class, we will take liberties with the terminology and say ``$f$ is a operator in $\Ll^1_v.$''\\

It is also interesting to note that neither $\Ll^1_v$ nor $\Ll_v$ are separable (Appendix \ref{nonsep}).

\subsection{The unital Banach algebra $\Ll^1_v$}
\begin{lemma}  If $v$ is admissible, then $\Ll_v^1$ is a unital Banach $*$-algebra of bounded operators acting on $\Ltwo.$

\begin{myproof} Considering first the unweighted case, if $f\in\Ll^1,$ then $f$ is simultaneously bounded on $\Linf$ ($\int_{\R^d} |\tilde{f}(x,y)|\,dx$ is bounded) and $\Lone$ ($\int_{\R^d} |\tilde{f}(x,y)|\,dy$ is bounded).  Thus $f$ is bounded on $\Lp,\,1\leq p\leq\infty,$ and in particular $\Ltwo,$ by Schur's test (see Appendix \ref{schurstest}).  So $\Ll^1$ is a set of operators bounded on $\Ltwo,$ but is it an algebra?  Yes, the algebra axioms ( \eqref{algebra1} - \eqref{algebra4}) are easily verified.  In fact $\Ll^1$ is a *-algebra.  Also, $\Ll^1$ is a Banach space, since $\Ll^1$ coincides with $L^{1,\infty}_1$ of \cite{Gro01}, which was shown to be a Banach space in, e.g., \cite{BePa61}.  It remains to show that $\Ll^1$ is a {\em normed algebra}; i.e., that $\|\cdot\|_{\Ll^1}$ is submultiplicative.  Let $f,g\in\Ll^1,$ with $f=\tilde{f}+\mu_1 e,\,g=\tilde{g}+\mu_2 e.$  The definition of $\|\cdot\|_{\Ll^1_v}$ yields
\be
\|fg\|_{\Ll^1}=\max\left\{\|\tilde{f}\tilde{g}\|_{op,\infty}\,,\,\|\tilde{f}\tilde{g}\|_{op,1}\right\}+|\mu_1||\mu_2|,\nn
\end{equation}
and we have
\bea
\|\tilde{f}\tilde{g}\|_{op,\infty}&=&\underset{\|h\|_\infty=1}{\sup}\|\tilde{f}\tilde{g}h\|_\infty\nn\\
&=&\underset{h\in\Linf}{\sup}\frac{\|\tilde{f}\tilde{g}h\|_\infty}{\|h\|_\infty}\nn\\
&=&\underset{h\in\Linf}{\sup}\frac{\|\tilde{f}\tilde{g}h\|_\infty}{\|h\|_\infty}\frac{\|\tilde{g}h\|_\infty}{\|\tilde{g}h\|_\infty}\nn\\
&=&\underset{h\in\Linf}{\sup}\frac{\|\tilde{f}\tilde{g}h\|_\infty}{\|\tilde{g}h\|_\infty}\frac{\|\tilde{g}h\|_\infty}{\|h\|_\infty}\nn
\end{eqnarray}
\vspace{-.1in}
\bea
&\leq&\underset{h\in\Linf}{\sup}\frac{\|\tilde{f}\tilde{g}h\|_\infty}{\|\tilde{g}h\|_\infty}\underset{h\in\Linf}{\sup}\frac{\|\tilde{g}h\|_\infty}{\|h\|_\infty}\nn\\
&\leq&\underset{h'\in\Linf}{\sup}\frac{\|\tilde{f}h'\|_\infty}{\|h'\|_\infty}\|\tilde{g}\|_{op,\infty}\nn\\
&=&\|\tilde{f}\|_{op,\infty}\|\tilde{g}\|_{op,\infty}\nn
\end{eqnarray}
Thus the $\infty$-norm is submultiplicative.  The $1$-norm is also submultiplicative, proven in the same manner.  And because
\bea
&&\max\biggl\{\|\tilde{f}\|_{op,\infty}\|\tilde{g}\|_{op,\infty}\,,\,\|\tilde{f}\|_{op,1}\|\tilde{g}\|_{op,1}\biggr\}+|\mu_1||\mu_2|\nn\\
&\leq&\max\biggl\{\|\tilde{f}\|_{op,\infty}\,,\,\|\tilde{f}\|_{op,1}\biggr\}+|\mu_1|\nn\\
&&\qquad\cdot\max\biggl\{\|\tilde{g}\|_{op,\infty}\,,\,\|\tilde{g}\|_{op,1}\biggr\}+|\mu_2|\,,\nn
\end{eqnarray}
we conclude that $\|\cdot\|_{\Ll^1}$ is submultiplicative, which establishes that $\Ll^1$ is a normed algebra, and thus a Banach algebra (with unit by contruction).\\

For the weighted case, we begin by showing that it is reasonable for $f\in\Ll^1_v$ to act on $\Ltwo,$ by first demonstrating that $\tilde{f}$ acts on $\Ltwo,$ then extending the conclusion to all of $\Ll^1_v.$  We can no longer use Schur's Test.  Let $h\in \Lp,$ with $\|h\|_p=1$ (it is easy to see why no generality is lost under this assumption).  Then the action of $\tilde{f}\in\Ll^1_v$ on $h\in \Lp,\,1\leq p\leq\infty,$ is given by
\be
(\tilde{f}h)(y)=\int_{\R^d}\tilde{f}(t,y)h(t)\,dt.
\end{equation}
Let $g_L(t,y)$ be the Laurent operator over $\R^{2d}$ formed from translations of
\be
g_a(y)=\tilde{f}\left(\mbox{ arg}\left\{\underset{x}{\max}\Bigl\{\underset{x\in\R^d}{\mbox{ess }\sup}\int_{\R^d}|\tilde{f}(t,x)||h(t)|\,dt\,\Bigr\}\right\}\,,\,y\right).\nn
\end{equation}
I.e., we form the Laurent operator $g_L(t,y):=g_a(y-t)$ from the slice of $\tilde{f}$ which acts on $h$ to yield the largest $\underset{y\in\R^d}{\mbox{ess }\sup}\int_{\R^d}|\tilde{f}(t,y)||h(t)|\,dt.$  It is clear that $g_L$ inherits membership in $\Ll^1_v$ from $\tilde{f},$ and we can conclude that $g_a\in\Lone.$  Then we have, for $1\leq p<\infty,$
\begin{eqnarray}
\|\tilde{f}h\|_p^p&=&\int_{\R^d}|(\tilde{f}h)(y)|^p\,dy\nn\allowdisplaybreaks\\
&=&\int_{\R^d}\Biggl|\int_{\R^d} \tilde{f}(t,y)h(t)\,dt\Biggr|^p\,dy\nn\allowdisplaybreaks\\
&\leq&\int_{\R^d}\Biggl(\int_{\R^d} |\tilde{f}(t,y)|\,|h(t)|\,dt\Biggr)^p\,dy\nn\allowdisplaybreaks\\
&\leq&\int_{\R^d}\Biggl(\int_{\R^d} |g_L(t,y)|\,|h(t)|\,dt\Biggr)^p\,dy\label{g1}\allowdisplaybreaks\\
&=&\int_{\R^d}\Biggl(\int_{\R^d} |g_a(y-t)|\,|h(t)|\,dt\Biggr)^p\,dy\nn\allowdisplaybreaks\\
&=&\int_{\R^d}(\,|g_a|*|h|\,(y)\,)^p\,dy\nn\allowdisplaybreaks\\
&=&\|(|g_a|*|\,h|)\|_p^p\nn\\
&\leq&\|g_a\|^p_1\|h\|_p^p\label{young2}\allowdisplaybreaks\\
&=&\|\tilde{f}\|^p_{\Ll^1_v}\|h\|_p^p\,<\,\infty\label{bydefgnv1} ,
\end{eqnarray}
where \eqref{g1} follows from the definition of $g,$ \eqref{young2} follows from Young's inequality with $r=p$ and $q=1,$ and \eqref{bydefgnv1} follows from the definition of $g_a(y)$ and the fact that $v(k)\geq 1.$   Extending the analysis to general $f=\tilde{f}+\mu e\in\Ll^1_v$ where $\mu\neq 0$ we get
\bea
\|fh\|^p_p&=&\|(\tilde{f}+\mu e)h\|^p_p\nn\\
&\leq&\|\tilde{f}h\|^p_p+|\mu|^p\|h\|^p_p\nn\\
&\leq&\|\tilde{f}\|^p_{\Ll^1_v}\|h\|^p_p+|\mu|^p\|h\|^p_p\nn\\
&\leq&\|f\|^p_{\Ll^1_v}\|h\|^p_p\,\,\,\mbox{ (recall $p\geq 1$)}.
\end{eqnarray}
Thus for all $h\in\Lp,$ we have $\|fh\|_p\leq M\|h\|_p$ and we conclude that any $f\in \Ll^1_v$ is a bounded linear operator on $\Lp$ (and in particular $\Ltwo$). \\

For $\Ll^1_v,$ the algebra axioms are all easily verified.  To show submultiplicativity, we again begin by considering elements of $\Ll^1_v$ which do not include the unit as a summand.  We have
\be
\int_{\R^d}|\tilde{f}\tilde{g}|(x,y)\,v(x\!-\!y)\,dx=\int_{\R^d} \biggl|\int_{\R^d}\tilde{f}(t,y)\tilde{g}(x,t)\,dt\biggr|\,v(x\!-\!y)\,dx\nn
\end{equation}
\vspace{-.3in}
\bea
&\leq&\int_{\R^d}\int_{\R^d}|\tilde{f}(t,y)|\,|\tilde{g}(x,t)|\,v(t-y)\,v(x\!-\!t)\,dt\,dx\nn\\
&\leq&\int_{\R^d}\int_{\R^d}|\tilde{f}(t,y)|\,|\tilde{g}(x,t)|\,v(t-y)\,v(x\!-\!t)\,dx\,dt\nn\\
&\leq&\int_{\R^d}|\tilde{f}(t,y)|\,v(t-y)\,\left(\int_{\R^d}|\tilde{g}(x,t)|\,v(x\!-\!t)\,dx\right)\,dt\nn
\end{eqnarray}
\vspace{-.3in}
\bea
&\leq&\int_{\R^d}\left(\underset{t\in\R^d}{\mbox{ess }\sup}\,\int_{\R^d}|\tilde{g}(x,t)|\,v(x\!-\!t)\,dx\right)\,|\tilde{f}(t,y)|\,v(t\!-\!y)\,dt\nn\\
&\leq&\|\tilde{f}\|_{\Ll^1_v}\|\tilde{g}\|_{\Ll^1_v}\,,\nn
\end{eqnarray}
and likewise for 
\be
\int_{\R^d}|\tilde{f}\,\tilde{g}|(x,y)\,v(x\!-\!y)\,dy.\nn
\end{equation}
Taking appropriate suprema yields $\|\tilde{f}\tilde{g}\|_{\Ll^1_v}\leq\|\tilde{f}\|_{\Ll^1_v}\|\tilde{g}\|_{\Ll^1_v}.$  Extending the analysis to the treatment of elements $f=\tilde{f}+\mu_1 e\,,\,g=\tilde{g}+\mu_2 e\in\Ll_v,$ we have 
\bea
\|fg\|_{\Ll^1_v}&=&\|(\tilde{f}+\mu_1e)(\tilde{g}+\mu_2e)\|_{\Ll^1_v}\nn\\
&=&\|\tilde{f}\tilde{g}+\mu_1 \tilde{g}e+\mu_2e \tilde{f}+\mu_1\mu_2e\|_{\Ll^1_v}\nn\\
&\leq&\|\tilde{f}\tilde{g}\|_{\Ll^1_v}+|\mu_1|\|\tilde{g}\|_{\Ll^1_v}+|\mu_2|\|\tilde{f}\|_{\Ll^1_v}+|\mu_1||\mu_2|\nn\\
&\leq&\|\tilde{f}\|_{\Ll^1_v}\|\tilde{g}\|_{\Ll^1_v}+|\mu_1|\|\tilde{g}\|_{\Ll^1_v}+|\mu_2|\|\tilde{f}\|_{\Ll^1_v}+|\mu_1||\mu_2|\nn\\
&=&(\|\tilde{f}\|_{\Ll^1_v}+|\mu_1|)(\|\tilde{g}\|_{\Ll^1_v}+|\mu_2|)\nn\\
&=&\|f\|_{\Ll^1_v}\|g\|_{\Ll^1_v}\label{submultLv1}
\end{eqnarray}
showing that $\|\cdot\|_{\Ll^1_v}$ submultiplicative, and we can conclude that $\Ll^1_v$ is a Banach algebra.
\end{myproof}
\end{lemma}

\subsection{The unital algebra $\Ll_v$}
\begin{definition}
A function $g$ is called {\em subconvolutive} if $g$ satisfies $g^{-1}*g^{-1}\leq Cg^{-1},$ for some $C>0.$  See Appendix \ref{subcon} for an elementary example of an admissible weight which is subconvolutive. 
\end{definition} 
\begin{lemma}  If $\int_{\R^d} v^{-1}(x)\,dx\,<\infty,$ then ${\mathcal L}_{v}\subset {\mathcal L}^1 .$  If in addition $v$ is subconvolutive, then ${\mathcal L}_{v}$ is a unital $*$-algebra of bounded operators acting on $\Ltwo.$\\
\begin{myproof}
As before, let $f=\tilde{f}+\mu e,\,f\in\Ll_v,$ with $v$ assumed subconvolutive.  We first show $\Ll_v\subset\Ll^1.$  Since
\be
|\tilde{f}(x,y)|v(x-y)\leq\underset{x',y'}{\mbox{ ess sup }}|\tilde{f}(x',y')|v(x'-y')\nn
\end{equation}
almost everywhere, we conclude that
\be
|\tilde{f}(x,y)|\leq \|\tilde{f}\|_{{\mathcal L}_{v}}v^{-1}(x-y),\label{fltnormnv}
\end{equation}
except possibly on a set of measure zero.  Then if $\int_\R v^{-1}(x)\,dx\,<\infty,$ we must have 
\bea
\underset{y\in\R^d}{\mbox{ess sup}}\left\{\int_\R |f(x,y)|\,dx\right\}&=&\underset{y\in\R^d}{\mbox{ess sup}}\left\{\int_\R |\tilde{f}(x,y)+\mu e(x,y)|\,dx\right\}\nn\\
&\leq&\|\tilde{f}\|_{{\mathcal L}_{v}}\int_\R v^{-1}(x-y)\,dx+|\mu|\\
&=&\|\tilde{f}\|_{{\mathcal L}_{v}}\int_\R v^{-1}(x)\,dx+|\mu|\,<\infty,
\end{eqnarray}
and similarly for $\underset{x\in\R^d}{\mbox{ess sup}}\left\{\int_\R |f(x,y)|\,dy\right\}.$  Thus ${\mathcal L}_{v}\subset {\mathcal L}^1 .$  Next, we show that ${\mathcal L}_v$ is an algebra of bounded operators on $\Lp.$  We begin by demonstrating that arbitrary $f\in{\mathcal L}_{v}$ in fact acts on $\Lp.$  Let $h\in\Lp.$  Recall the (formal) action of $f\in{\mathcal L}_v$ on $h$ is given by
\be
(fh)(y)=\int_\R f(t,y)h(t)\,dt\nn.
\end{equation}
Then
\bea
\int_\R|(fh)(y)|^pdy\!&=&\!\int_\R\left|\int_\R f(t,y)h(t)\,dt\right|^pdy\nn\\
&\leq&\|f\|_{{\mathcal L}_{v}}^p\int_\R\left|\int_\R v^{-1}(t-y)h(t)\,dt\right|^pdy\label{eq2}\\
&=&\|f\|_{{\mathcal L}_{v}}^p\cdot\|v^{-1}*h\|_p^p\label{vinvconv}\\
&\leq&\|f\|_{{\mathcal L}_{v}}^p\cdot\|v^{-1}\|^p_1\cdot\|h\|_p^p\,<\,\infty\label{young1}
\end{eqnarray}
where \eqref{eq2} follows from the estimate \eqref{fltnormnv}, \eqref{vinvconv} is a consequence of the symmetry of $v,$ and \eqref{young1} follows from Young's inequality with $r=q$ and $p=1.$  Thus for all $h\in\Lp,\,\|fh\|_p<\infty$ and we conclude that any $f\in {\mathcal L}_v$ induces a bounded linear operator on $\Lp$ (in particular for $p=2$).  To show that ${\mathcal L}_v$ is an algebra, we need to show that the class ${\mathcal L}_v$ is closed under multiplication.  Without loss of generality we can assume the submultiplicativity constant $C\geq 1.$  We have
\be
|(\tilde{f}\tilde{g})(x,y)|=\left|\int_\R \tilde{f}(t,y)\tilde{g}(x,t)\,dt\right|\nn
\end{equation}
\vspace{-.2in}
\bea
&\leq&\|\tilde{f}\|_{{\mathcal L}_{v}}\|\tilde{g}\|_{{\mathcal L}_{v}}\!\left|\int_\R \!v^{-1}(t\!-\!y)v^{-1}(x\!-\!t)\,dt\right|\label{excmz1}\\
&\leq&\|\tilde{f}\|_{{\mathcal L}_{v}}\|\tilde{g}\|_{{\mathcal L}_{v}}\,C\,v^{-1}(x-y)\label{subcon1}
\end{eqnarray}
\vspace{-.2in}
\bea
\Rightarrow |(\tilde{f}\tilde{g})(x,y)|v(x-y)&\leq& C\|\tilde{f}\|_{{\mathcal L}_{v}}\|\tilde{g}\|_{{\mathcal L}_{v}}.\label{excmz2}\nn\\
\Rightarrow\qquad\qquad\|\tilde{f}\tilde{g}\|_{{\mathcal L}_v}&\leq& C\|\tilde{f}\|_{{\mathcal L}_{v}}\|\tilde{g}\|_{{\mathcal L}_{v}},\label{submultLvtilde}
\end{eqnarray}
where \eqref{subcon1} follows from the subconvolutivity of $v,$ and \eqref{excmz1} and \eqref{excmz2} hold except possibly on a set of measure zero.  Invoking \eqref{submultLvtilde} in a calculation nearly identical to \eqref{submultLv1} yields
\bea
\|fg\|_{\Ll_v}&=&\|(\tilde{f}+\mu_1\d)(\tilde{g}+\mu_2\d)\|_{\Ll_v}\nn\\
&=&\|\tilde{f}\tilde{g}+\mu_1 \tilde{g}\d+\mu_2\d \tilde{f}+\mu_1\mu_2\d\|_{\Ll_v}\nn\\
&\leq&\|\tilde{f}\tilde{g}\|_{\Ll^1_v}+|\mu_1|\|\tilde{g}\|_{\Ll^1_v}+|\mu_2|\|\tilde{f}\|_{\Ll_v}+|\mu_1||\mu_2|\nn\\
&\leq&C\,\|\tilde{f}\|_{\Ll^1_v}\|\tilde{g}\|_{\Ll^1_v}+|\mu_1|\|\tilde{g}\|_{\Ll^1_v}+|\mu_2|\|\tilde{f}\|_{\Ll_v}+|\mu_1||\mu_2|\nn\\
&\leq&C\,\|\tilde{f}\|_{\Ll^1_v}\|\tilde{g}\|_{\Ll^1_v}+C\,|\mu_1|\|\tilde{g}\|_{\Ll^1_v}+C\,|\mu_2|\|\tilde{f}\|_{\Ll_v}+C\,|\mu_1||\mu_2|\nn\\
&=&C\,(\|\tilde{f}\|_{\Ll_v}+|\mu_1|)(\|\tilde{g}\|_{\Ll_v}+|\mu_2|)\nn\\
&=&C\,\|f\|_{\Ll_v}\|g\|_{\Ll_v}\label{submultLv}
\end{eqnarray}
Thus $fg$ satisfies the definition of ${\mathcal L}_v.$  The class $\Ll_v$ possesses a unit by contruction, and thus ${\mathcal L}_v$ is a unital algebra of bounded linear operators on $\Ltwo.$
\end{myproof}
\end{lemma}

We see that under the norm $\|\cdot\|_{\Ll_v}$ and the operator multiplication $fg,\, \Ll_v$ is not strictly submultiplicative, so we seek a different norm under which $\Ll_v$ might satisfy submultiplicativity.  A reasonable candidate might be $\|\cdot\|'_{{\mathcal L}_v}$ defined on ${\mathcal L}_v$ by
\be
\|f\|'_{{\mathcal L}_v}=\underset{h\in {\mathcal L}_v\,,\,\|\,h\|_{{\mathcal L}_v}=1}{\mbox{ess }\sup}\|fh\|_{{\mathcal L}_v},
\end{equation}
which is the continuous analog to the norm $\|\cdot\|'_{\A_v}$ mentioned briefly in \cite{GL03}.  
Unfortunately, the two norms $\|\cdot\|_{\Ll_v}$ and $\|\cdot\|'_{\Ll_v}$ are not equivalent.   In \cite{GL03} the corresponding norms $\|\cdot\|'_{\A_v}$ and $\|\cdot\|_{\A_v}$ are equivalent but the relationship breaks down in the continuous case; consider the sequence $\{f_n\}_{n\in\N}$ where $f_n=v^{-1}(x)\chi_{[-\frac{1}{n},\frac{1}{n}]}(x).$  We have $\|f_n\|_{\Ll_v}=1\,\,\forall\,\,n,$ but $\underset{n\rightarrow\infty}{\lim}\|f_n\|'_{\Ll_v}\rightarrow 0.$\\ 

\subsubsection{Transforming $\Ll_v$ into a Banach space}
We define a new norm and space, and show that, under broad circumstances, under the new norm, the new space is equal to $\Ll_v.$\\

Fix $s>0,$ and let $u$ be an admissible weight on $\R^d.$  Let $\t_s(x)=(1+|x|)^s,$ (also admissible) and set
\be
v(x)=u(x)\t_s(x).
\end{equation}
Note that
\be
\t_s(x+y)\leq\k(\t_s(x)+\t_s(y))\label{polysubadd}
\end{equation}
for some $\k.$  To prove this we consider two cases (in each we assume $x,y>0$):\\   
\bea
&&x>y\quad\Rightarrow\quad\frac{(1+x+y)^s}{(1+x)^s+(1+y)^s}<\frac{(1+2x)^s}{(1+x)^s}<2^s;\quad\mbox{ and }\nn\\
&&x\leq y\quad\Rightarrow\quad\frac{(1+x+y)^s}{(1+x)^s+(1+y)^s}\leq\frac{(1+2y)^s}{(1+y)^s}<2^s;\nn
\end{eqnarray}
note that $\k>1.$  We also have 
\bea
\!\!\!\!\!\!\!v(x\!+\!y)&=&u(x\!+\!y)\t_s(x\!+\!y)\nn\\
&\leq&\k u(x)u(y)(\t_s(x)\!+\!\t_s(y))\nn\\
&=&\k[v(x)u(y)+u(x)v(y)].\label{wtineq2}
\end{eqnarray}
Define the Banach space $\B_{u,s}:=\Ll^1_u\cap\Ll_v$ with norm
\be
\|f\|_{\B_{u,s}}=\k\|\tilde{f}\|_{\Ll^1_u}+\|\tilde{f}\|_{\Ll_v}+|\mu|.\label{bnorm}
\end{equation}

\begin{lemma}
(a) The class $\B_{u,s}$ with norm \eqref{bnorm} is an involutive unital Banach algebra.\\
(b) If $s>d,$ then $\Ll_v\subset\Ll^1_u$ and thus $\B_{u,s}=\Ll_v$ with equivalent norms.

\begin{myproof}
(a) Suppose $f,g\in\B_{u,s}.$  We have
\bea
\|\tilde{f}\tilde{g}\|_{\Ll_v}&=&\underset{x,y\in\R^d}{\mbox{ess }\sup}\,\,\biggl|\int_{\R^d}\tilde{f}(t,y)\tilde{g}(x,t)\,dt\,\biggr|v(x-y)\nn\allowdisplaybreaks\\
&\leq&\underset{x,y\in\R^d}{\mbox{ess }\sup}\,\,\int_{\R^d}\k\,|\tilde{f}(t,y)|\,v(t-y)\,|\tilde{g}(x,t)|\,u(x-t)\,dt\nn\\
&&+\,\underset{x,y\in\R^d}{\mbox{ess }\sup}\,\,\int_{\R^d}\k\,|\tilde{f}(t,y)|\,u(t-y)\,|\tilde{g}(x,t)|\,v(x-t)\,dt\nn\\
&\leq&\k\,\|\tilde{f}\|_{\Ll_v}\underset{x\in\R^d}{\mbox{ess }\sup}\,\,\int_{\R^d}|\tilde{g}(x,t)|\,u(x-t)\,dt\nn\\
&&+\,\k\,\|\tilde{g}\|_{\Ll_v}\underset{y\in\R^d}{\mbox{ess }\sup}\,\,\int_{\R^d}|\tilde{f}(t,y)|\,u(t-y)\,dt\nn\allowdisplaybreaks\\
&\leq&\k\,(\|\tilde{f}\|_{\Ll_v}\|\tilde{g}\|_{\Ll^1_u}+\|\tilde{f}\|_{\Ll^1_u}\|\tilde{g}\|_{\Ll_v}).\label{submultaddav}
\end{eqnarray}
Thus it follows that
\bea
\|fg\|_{\Ll_v}&=&\|(\tilde{f}+\mu_1 e)(\tilde{g}+\mu_2 e)\|_{\Ll_v}\nn\allowdisplaybreaks\\
&=&\|\tilde{f}\tilde{g}+\mu_2\tilde{f}+\mu_1\tilde{g}+\mu_1\mu_2 e\|_{\Ll_v}\nn\allowdisplaybreaks\\
&\leq&\|\tilde{f}\tilde{g}\|_{\Ll_v}+|\mu_2|\|\tilde{f}\|_{\Ll_v}+|\mu_1|\|\tilde{g}\|_{\Ll_v}+|\mu_1||\mu_2|\nn\allowdisplaybreaks\\
&\leq&\k(\|\tilde{f}\|_{\Ll_v}\|\tilde{g}\|_{\Ll^1_u}+\|\tilde{f}\|_{\Ll^1_u}\|\tilde{g}\|_{\Ll_v})\nn\allowdisplaybreaks\\
&&\quad+|\mu_2|\|\tilde{f}\|_{\Ll_v}+|\mu_1|\|\tilde{g}\|_{\Ll_v}+|\mu_1\|\mu_2|\nn\allowdisplaybreaks\\
&\leq&\k(\|f\|_{\Ll_v}\|g\|_{\Ll^1_u}+\|f\|_{\Ll^1_u}\|g\|_{\Ll_v}).\label{cruxoflastpart}
\end{eqnarray}
(We actually do not need \eqref{cruxoflastpart} to show that $\|\cdot\|_{\B_{u,s}}$ is submultiplicative, but we shall need it in Section \ref{invclsdforsgd}).  We also find
\bea
\|\wt{(fg)}\|_{\Ll_v}&\leq&\|\tilde{f}\tilde{g}\|_{\Ll_v}+|\mu_2|\|\tilde{f}\|_{\Ll_v}+|\mu_1|\|\tilde{g}\|_{\Ll_v}\nn\allowdisplaybreaks\\
&\leq&\k(\|\tilde{f}\|_{\Ll_v}\|\tilde{g}\|_{\Ll^1_u}+\|\tilde{f}\|_{\Ll^1_u}\|\tilde{g}\|_{\Ll_v})\nn\allowdisplaybreaks\\
&&\quad+|\mu_2|\|\tilde{f}\|_{\Ll_v}+|\mu_1|\|\tilde{g}\|_{\Ll_v}.\label{cruxofnextpart}
\end{eqnarray}
To establish that $\|\cdot\|_{\B_{u,s}}$ is submultiplicative, we require that
\be
\|\wt{(fg)}\|_{\Ll^1_u}\leq\|\tilde{f}\|_{\Ll^1_u}\|\tilde{g}\|_{\Ll^1_u}+\mu_2\|\tilde{f}\|_{\Ll^1_u}+\mu_1\|\tilde{g}\|_{\Ll^1_u},\label{tilde1}
\end{equation}
which is easily proven using the definitions.  Using \eqref{cruxofnextpart} and \eqref{tilde1}, we have
\bea
\!\!\!\!\!\!\|fg\|_{\B_{u,s}}&=&\k\|\widetilde{(fg)}\|_{\Ll^1_u}+\|\widetilde{(fg)}\|_{\Ll_v}+|\mu_1||\mu_2|\nn\\
&\leq&\k\,(\|\tilde{f}\|_{\Ll^1_u}\|\tilde{g}\|_{\Ll^1_u}+\mu_2\|\tilde{f}\|_{\Ll^1_u}+\mu_1\|\tilde{g}\|_{\Ll^1_u})\nn\\
&+&\k(\|\tilde{f}\|_{\Ll_v}\|\tilde{g}\|_{\Ll^1_u}+\|\tilde{f}\|_{\Ll^1_u}\|\tilde{g}\|_{\Ll_v})+\mu_2\|\tilde{f}\|_{\Ll_v}+\mu_1\|\tilde{g}\|_{\Ll_v}\nn\\
&+&|\mu_1||\mu_2|.
\end{eqnarray} 
But 
\bea
\|f\|_{\B_{u,s}}\|g\|_{\B_{u,s}}\!\!\!\!&=&\!\!\k^2\|\tilde{f}\|_{\Ll^1_u}\|\tilde{g}\|_{\Ll^1_u}\nn\\
&&+\k(\|\tilde{f}\|_{\Ll_v}\|\tilde{g}\|_{\Ll^1_u}+|\mu_1|\|\tilde{g}\|_{\Ll^1_u})\nn\\
&&+\k(\|\tilde{f}\|_{\Ll^1_u}\|\tilde{g}\|_{\Ll_v}+|\mu_2|\|\tilde{f}\|_{\Ll^1_u})\nn\\
&&+\|\tilde{f}\|_{\Ll_v}\|\tilde{g}\|_{\Ll_v}+|\mu_2|\|\tilde{f}\|_{\Ll_v}\nn\\
&&+|\mu_1|\|\tilde{g}\|_{\Ll_v}+|\mu_1||\mu_2|.
\end{eqnarray}
And since $\|\cdot\|_{\Ll^1_u}$ and $\|\cdot\|_{\Ll_v}$ are norms and thus nonnegative, and $\k\geq 1,$ we immediately obtain $\|fg\|_{\B_{u,s}}\leq\|f\|_{\B_{u,s}}\|g\|_{\B_{u,s}}.$\\

The adjoint operation is an involution on $\B_{u,s}.$  The intersection of two Banach spaces is again a Banach space as we saw in Section \ref{BanHilSp}; thus $\B_{u,s}$ with norm \eqref{bnorm} is an involutive Banach algebra.\\

(b)  If $s>d,$ then $\t^{-1}_s(x)$ is integrable, so
\bea
\|f\|_{\Ll^1_u}&=&\|\tilde{f}+\mu e\|_{\Ll^1_u}\nn\\
&=&\underset{y\in\R^d}{\mbox{ess }\sup}\int_{\R^d}|\tilde{f}(t,y)|u(t-y)\t_s(t-y)\t^{-1}_s(t-y)\,dt+|\mu|\nn\\
&\leq&\|\tilde{f}\|_{\Ll_v}\,\,\int_{\R^d}\t^{-1}_s(t-y)\,dt\,+\,|\mu|\nn\\
&=&A\,\|\tilde{f}\|_{\Ll_v}+|\mu|.
\end{eqnarray}
Now let $A'=\max\{1,A\}.$  Then we have $\|f\|_{\Ll^1_u}\leq A'\,\|f\|_{\Ll_v}$
and thus the norms on $\B_{u,s}$ and $\Ll_v$ are equivalent:
\be
\|f\|_{\Ll_v}\leq\k\|f\|_{\Ll^1_u}+\|f\|_{\Ll_v}=\|f\|_{\B_{u,s}}\leq (\k A'+1)\|f\|_{\Ll_v},\nn
\end{equation}
and $\Ll_v\subset\Ll^1_u$ so $\B_{u,s}=\Ll_v.$
\end{myproof}
\end{lemma} 

\section{Symmetry, and decay of inverses}
\subsection{Symmetry of algebras}
Before we show the main result on the decay of $\star$-inverse kernels, we return briefly to the question of algebra symmetry, since the symmetry results shall be proven along the way.
Recall that $\s_\A(f)$ is the spectrum of $f$ {\em in the algebra} $\A,$ and $\s(f)$ is the spectrum of the operator $f$ acting on a Hilbert space $\Hi.$  In general, $\s_\A(f)\neq\s(f);$ if you can read French, see reference (9) in \cite{Hul72}.  The respective spectral radii are 
\bea
\r_\A(f)&=&\max\left\{\,|\la|\,:\,\la\in\s_\A(f)\right\}\nonumber\\
\r(f)&=&\max\left\{\,|\la|\,:\,\la\in\s(f)\right\}.\nonumber
\end{eqnarray}
If $f$ is self-adjoint, then $\r(f)=\|f\|_{op}$ (see, e.g., \cite{Nach01}).
\begin{lemma}\label{alginclspecsmaller}
If $ \A\subset{\mathcal B},$ then we have $\r_{\mathcal B}(f)\leq\r_ \A(f),\,\forall\,f\in \A.$\\
\begin{myproof}
For $\la\in\s_ \A(f)$ there may exist an element in ${\mathcal B}\setminus \A$ which inverts $\la I-f.$
\end{myproof}
\end{lemma}

A pair of nested algebras $\A\subset\B$ is called a {\it Wiener pair} if $f\in\A$ and $f^{-1}\in\B$ implies $f^{-1}\in\A.$  $\A$ is often termed {\it inverse closed} in $\B.$\\

A Banach algebra $\A$ is called {\it symmetric} if $\s_\A(f^*f)\subset[0,\infty)\,\forall\,f\in\A.$  Equivalently, $\A$ is symmetric $\Leftrightarrow\,f=f^*$ implies $\s_\A(f)\subset\R.$  Recall that if $f=f^*$ is a bounded linear operator on a Hilbert space, then $\s(f)\subset\R$ (see, e.g., \cite{Nach01}), so if we can establish for a particular $\A$ that $\forall\,\,f=f^*$ we have $\s_\A(f)=\s(f),$ we will have shown that $\A$ is symmetric.

\subsection{Crucial lemmas}\label{HBSlemmas}
\begin{lemma}\label{hulanicki}{\it (Hulanicki \cite{Hul72})}:  If $ \A$ is an involutive Banach algebra with identity contained in ${\mathcal B}({\mathcal H}),$ and if for all $f=f^*\,\in\,\A$ we have
\be
\r_ \A(f)=\r(f)=\|f\|_{op}
\end{equation}
then $\s_\A(f)=\s(f)$ for all $f=f^*\in\A.$
\end{lemma}

As a consequence of Hulanicki's lemma, we conclude that $\A$ is a symmetric Banach algebra.  In addition, $\A$ contains the inverses of each of its self-adjoint elements.  To see this, suppose $\forall\,f=f^*\in\A, \s(f)=\s_\A(f).$  Then if $f$ is not invertible on $\A,$ we have $0\in\s_{\A}(f)$ and thus $0\in\s(f)$ so $f$ is not invertible on $\B(\Hi).$   Thus by the contrapositive we conclude that if $f$ is invertible on $\B(\Hi),$ then $f$ is invertible on $\A.$  If we can extend this latter conclusion to all of $\A,$ then we will have shown that $\A$ and $\B(\Hi)$ form a Wiener pair.  

\begin{lemma}\label{barneslemma}{\it (Barnes \cite{Bar87})}: Suppose that $w(x)=(1+|x|)^\d$ for some $\d\in(0,1].$  Then for all $g\in\Ll^1_w$ we have
\be
\r_{\Ll^1_w}(g)=\r_{\Ll^1}(g)=\|g\|_{op}.\label{barnes}
\end{equation}
\end{lemma}
This demonstrates that $\Ll^1_w$ is symmetric, since we have 
\be
\r_{\Ll^1_w}(g)=\r(g)=\|g\|_{op}\label{barnes1}
\end{equation}
for all $g\in\Ll^1_w,$ and thus of course \eqref{barnes1} holds for self-adjoint elements as well, which implies satisfaction of Hulanicki's lemma for $\Ll^1_w.$  We can also conclude, using distributions closely related to pseudoinverses in a technique which will be described later in this section, that $\Ll^1_w$ is inverse closed in $\B(\Ltwo)$ for the above weights.  It is apparently an open question whether or not $\Ll^1$ is symmetric (see Section \ref{L1symm}).\\

We restate a lemma from \cite{GL03}:
\begin{lemma}(Gr\"ochenig, Leinert, Pytlik)\label{auxweights}
For any unbounded admissible weight function $v$ there exists a sequence of admissible weights $v_n,\,n\!\!\in\!\N,$ with the following properties:\\
(a)$\,\,v_{n+1}\leq v_n\leq v,\,\,\forall\,n$\\
(b)$\,\,\,\exists\,\,\,c_n>0$ such that $v\leq c_nv_n,$\\
(c)$\,\,\lim_{n\rightarrow\infty}v_n=1$ uniformly on compact sets of $\R^d,$ and\\
(d)$\,\,\,\exists\,\,\,\b_n$ such that $v_n(z)=e^{-n}v(z)$ for $\|z\|>\b_n.$
\end{lemma}

\subsection{$\Ll^1_v$ is inverse-closed in $\B(\Ltwo)$}
\begin{lemma}\label{specradlemma}
Suppose $v(x)\geq C(1+|x|)^\d$ with $\d>0.$  Under the hypotheses of Lemma \ref{auxweights} the following identities hold $\,\forall\,f=f^*\in\Ll^1_v:$
\bea
&&(a)\qquad\underset{n\rightarrow\infty}{\lim}\|f\|_{\Ll^1_{v_n}}=\|f\|_{\Ll^1},\mbox{ and }\label{limnorm}\\
&&(b)\qquad\r_{\Ll^1_v}(f)=\r_{\Ll^1}(f)=\|f\|_{op}\label{specrad}
\end{eqnarray}
\begin{myproof}
(a) Let $\e>0.$ Then if $f=f^*=\tilde{f}+\mu e\in\Ll^1_v,$ we have 
\be
\|\tilde{f}\|_{\Ll^1_{v_n}}=\underset{y\in\R^d}{\mbox{ess }\sup}\,\int_{\R^d}\,|\tilde{f}(x,y)|\,v_n(x\!-\!y)\,dx.\nn
\end{equation}
Part (d) of Lemma \ref{auxweights} states $v_n(z)=e^{-n}v(z)$ for $\|z\|>\b_n,$ so $\exists\,\, n_0=n_0(\e)\in\N$ such that 
\be
\underset{y\in\R^d}{\mbox{ess }\sup}\underset{x\in\R^d:\|x-y\|\geq\,\b_{n_0}}{\int}|\tilde{f}(x,y)|\,v_{n_0}(x\!-\!y)\,dx\leq e^{-n_0}\|\tilde{f}\|_{\Ll^1_v}<\e.\nn
\end{equation}
Invoking this and monotonicity $v_{n+1}\leq v_n\leq v,$ we get for all $n\geq n_0,$
\be
\underset{y\in\R^d}{\mbox{ess }\sup}\!\!\!\!\underset{x\in\R^d:\|x-y\|\geq\,\b_{n_0}}{\int}\!\!\!\!\!\!\!\!\!\!\!\!|\tilde{f}(x,y)|\,v_n(x\!-\!y)\,dx<\e.\nn
\end{equation}
By property (c) of the $v_n,$ we have that $v_n\rightarrow 1$ uniformly if $\|z\|\leq \b_{n_0},$ so for $n\geq\n_1=n_1(\e),$
\be
\underset{y\in\R^d}{\mbox{ess }\sup}\!\!\!\!\underset{x\in\R^d:\|x-y\|\leq\,\b_{n_0}}{\int}\!\!\!\!\!\!\!\!\!\!\!\!|\tilde{f}(x,y)|\,v_n(x\!-\!y)\,dx\,\leq\,(1+\e)\,\underset{y\in\R^d}{\mbox{ess }\sup}\int_{x\in\R^d}|\tilde{f}(x,y)|.\nn
\end{equation}
These last two inequalities yield, for $n>\max\,\{n_0,n_1\},$
\be
\|\tilde{f}\|_{\Ll^1_{v_n}}\leq\e+(1+\e)\|\tilde{f}\|_{\Ll^1}.\label{prelimit}
\end{equation}
We take the limit of \eqref{prelimit} to get $\underset{n\rightarrow\infty}{\lim}\|\tilde{f}\|_{\Ll^1_{v_n}}\leq\|\tilde{f}\|_{\Ll^1}.$  Since $\|f\|_{\Ll^1_{v_n}}=\|\tilde{f}\|_{\Ll^1_{v_n}}+|\mu|$ and $\|f\|_{\Ll^1_v}=\|\tilde{f}\|_{\Ll^1_v}+|\mu|,$ we obtain $\underset{n\rightarrow\infty}{\lim}\|f\|_{\Ll^1_{v_n}}\leq\|f\|_{\Ll^1}.$  It is clear that since $\,\forall\,n,\,v_n\geq 1,$ we have $\underset{n\rightarrow\infty}{\lim}\|f\|_{\Ll^1_{v_n}}\geq\|f\|_{\Ll^1},$
so we conclude $\underset{n\rightarrow\infty}{\lim}\|f\|_{\Ll^1_{v_n}}=\|f\|_{\Ll^1}.$\\
\indent (b) Since $v$ and $v_n$ are equivalent and $\r(g)\leq\|g\|$ for any Banach space norm, we have
\be
\r_{\Ll^1_v}(f)^k=\r_{\Ll^1_v}(f^k)=\r_{\Ll^1_{v_n}}(f^k)\leq\|f^k\|_{\Ll^1_{v_n}}\qquad\forall\,\,k,n,\in\N.\nn
\end{equation}
This and part (a) yield 
\be
\r_{\Ll^1_v}(f)^k\leq\underset{n\rightarrow\infty}{\lim}\|f^k\|_{\Ll^1_{v_n}}=\|f\|_{\Ll^1}\qquad\forall\,\,k,\in\N.\nn
\end{equation}
Taking $k^{th}$ roots gives, by continuity of $x^k,$
\be
\r_{\Ll^1_v}(f)\leq\underset{n\rightarrow\infty}{\lim}\|f^k\|_{\Ll^1_{v_n}}^{1/k}=\r_{\Ll^1}(f).\label{speclessop}
\end{equation}
Now $v(x)\geq C(1+|x|)^\d=\t_\d(x)$ and $\d>0\Rightarrow\Ll^1_v\subset\Ll^1.$  By Lemma \ref{alginclspecsmaller}, we have the reverse inequality of \eqref{speclessop}, namely $\r_{\Ll^1}(f)\leq\r_{\Ll^1_v}(f).$  Thus $\r_{\Ll^1}(f)=\r_{\Ll^1_v}(f).$  Barnes' Lemma \ref{barneslemma} implies $\r_{\Ll^1}(f)=\|f\|_{op}.$  We conclude $\r_{\Ll^1_v}(f)=\|f\|_{op}.$
\end{myproof}
\end{lemma}
This is an important lemma, since it affords us the following theorem.

\begin{theorem}\label{big1}
Assume that $v$ is an admissible weight which satisfies the weak growth condition
\be
v(x)\geq C(1+\|x\|)^\d\qquad\mbox{ for some }\d>0.\label{weakgrowth}
\end{equation}
Then $\,\,\forall\,\,f\in\Ll^1_v,\,\,\s_{\Ll^1_v}(f)=\s(f),$ thus $\Ll^1_v$ is a symmetric Banach algebra and inverse closed in $\B(\Ltwo).$\\

\begin{myproof}
By \eqref{specrad}, $\Ll^1_v$ satisfies the hypotheses of Hulanicki's Lemma \ref{hulanicki}, and we obtain $\s_{\Ll^1_v}(f)=\s(f),\,\,\forall\,\,f=f^*\in\Ll^1_v.$  Thus is $\Ll^1_v$ a symmetric Banach algebra.  Also, if $f=f^*\in\Ll^1_v$ is invertible on $\Ltwo,$ then $f$ is invertible as an element of $\Ll^1_v.$  To apply this conclusion to elements of $\Ll^1_v$ which are not self-adjoint (i.e., to conclude that $\Ll^1_v$ is inverse closed in $\B(\Ltwo)$), let $f\in\Ll^1_v$ be invertible on $\B(\Ltwo)$ but not necessarily self-adjoint.  Since $f$ is invertible on $\B(\Ltwo),$ we know $f^*f, ff^*$ are invertible on $\B(\Ltwo)$ as well.  Indeed, if $f$ is invertible $\exists\,\,g\,\ni\,fg=gf=e=\d(x-y).$  Now $\d(x-y)$ is self-adjoint, and by \eqref{fgaster} we have
\be
g^*f^*=(fg)^*=\d(x-y)=(gf)^*=f^*g^*\nn
\end{equation}
so $f^*$ has inverse $g^*.$  Then by \eqref{assoc}
\bea
(ff^*)(g^*g)&=&(f(f^*g^*)g)\nn\\
&=&f(\d(x-y)g)\nn\\
&=&fg=\d(x-y);\nn
\end{eqnarray}
the inverse of $f^*f$ is computed similarly.  We conclude that $f^*(ff^*)^{-1}$ and $(f^*f)^{-1}f^*$ are, respectively, right- and left-inverses of $f$ in $\Ll^1_v.$  Thus even if $f$ is not self-adjoint, if $f$ is invertible on $\B(\Ltwo),$ then $f$ is invertible on $\Ll^1_v,$ and thus $\Ll^1_v$ is inverse closed in $\B(\Ltwo).$  We conclude that if $f\in\Ll^1_v,$ then the right- and left-inverses of $f$ share the same decay properties as $f.$
\end{myproof}
\end{theorem}
\subsection{$\Ll_v$ is inverse-closed in $\B(\Ltwo)$, for $s>d$}\label{invclsdforsgd}
We next consider $\B_{u,s}$ and $\Ll_v.$  The main result:
\begin{theorem}\label{big2} Assume that $u$ is an admissible weight, $s>0,$ and $\,v=u\t_s.$  Then
\be
\r_{\B_{u,s}}(f)=\|f\|_{op}\qquad \forall\,\,f=f^*\in\B_{u,s}.
\end{equation}
Thus $\B_{u,s}$ is a symmetric Banach algebra, and inverse closed in $\B(\Ltwo).$\\
\begin{myproof}
Writing $v=(u\t_{\frac{s}{2}})\t_{\frac{s}{2}},$ we can assume without loss of generality that $u$ satisfies condition \eqref{weakgrowth}, i.e., $\Ll^1_u\subset\Ll^1_{\t_s/2}.$  Let $f=\tilde{f}+\mu e\in\B_{u,s}\subset\Ll_v.$ Then \eqref{cruxoflastpart} yields
\be
\|f^{2n}\|_{\Ll_v}\leq 2\k\|f^n\|_{\Ll^1_u}\|f^n\|_{\Ll_v},\qquad \forall\,n\geq 1.\nn
\end{equation}
which then implies 
\bea
\|f^{2n}\|_{\B_{u,s}}&=&\|\widetilde{f^{2n}}\|_{\Ll_v}+\|\widetilde{f^{2n}}\|_{\Ll^1_u}+|\mu|^{2n}\nn\\
&\leq& 2\k\|\widetilde{f^n}\|_{\Ll^1_u}\|\widetilde{f^n}\|_{\Ll_v}+\|\widetilde{f^n}\|^2_{\Ll^1_u}+|\mu|^{2n}\nn\\
&\leq& 2\k\|\widetilde{f^n}\|_{\Ll^1_u}\|\widetilde{f^n}\|_{\Ll_v}+2\k\|\widetilde{f^n}\|^2_{\Ll^1_u}+2\k|\mu|^{2n}\nn\\
&\leq&2\k(\|\widetilde{f^n}\|_{\Ll^1_u}+|\mu|^n)(\|\widetilde{f^n}\|_{\Ll^1_u}+\|\widetilde{f^n}\|_{\Ll_v}+|\mu|^n)\nn\\
&=& 2\k\|f^n\|_{\Ll^1_u}\|f^n\|_{\B_{u,s}}\label{ineqbusu}
\end{eqnarray}
The spectral radius formula is $\r_\A(g)=\underset{n\rightarrow\infty}{\lim}\|g^n\|^{1/n}_{\A}.$  So we have
\bea
\r_{\B_{u,s}}(f)&=&\underset{n\rightarrow\infty}{\lim}\|f^{2n}\|^{1/{2n}}_{\B_{u,s}}\nn\\
&\leq&\underset{n\rightarrow\infty}{\lim}(2\k)^{1/{2n}}\|f^n\|^{1/{2n}}_{\Ll^1_u}\|f^n\|^{1/{2n}}_{\B_{u,s}}\nn\\
&=&\sqrt{\r_{\Ll^1_u}(f)}\sqrt{\r_{\B_{u,s}}(f)},\nn
\end{eqnarray}
which implies
\be
{\r_\B}_{u,s}(f)\leq\r_{\Ll^1_u}(f).\nn
\end{equation}
Since $\Ll^1_u\subset\Ll^1_{\t_s/2},$ we can apply \eqref{specrad} of Lemma \ref{specradlemma} and get
\be
\r_{\B_{u,s}}(f)\leq\r_{\Ll^1_u}(f)=\|f\|_{op}.\nn
\end{equation}
It is clear that 
\be
\r_{\B_{u,s}}(f)\geq\r_{\Ll^1_u}(f)\nn
\end{equation}
since $\B_{u,s}\subset\Ll^1_u.$  Thus we have shown
\be
\r_{\B_{u,s}}(f)=\r_{\Ll^1_u}(f).\nn
\end{equation}
We conclude that $\s_{\B_{u,s}}(f)=\s(f)$ and $\B_{u,s}$ is symmetric and inverse closed in $\B(\Ltwo).$\end{myproof}
\end{theorem}
Finally if $s>d,$ the algebras $\B_{u,s}$ and $\Ll_v$ coincide, and we have the following:
\begin{corollary}
Suppose that $u$ is an admissible weight, $s>d,$ and $\,v=u\t_s.$  Suppose also that $f$ is invertible on $\Ltwo$ and satisfies the off-diagonal decay condition
\be
\underset{x,y\in\R^d}{\mbox{\textup{ess} }\sup}\,|\tilde{f}(x,y)|\,v(x\!-\!y)\,<\infty.\label{defLv}
\end{equation}
Then $g=f^{-1}$ satisfies the same condition, i.e.,
\be
\underset{x,y\in\R^d}{\mbox{\textup{ess} }\sup}\,|\tilde{g}(x,y)|\,v(x\!-\!y)\,<\infty.
\end{equation}
\begin{myproof}
Suppose first that $f^*=f.$  Since $s\!>\!d,\,\B_{u,s}\!=\!\Ll_v$ and Theorem \ref{big2} says $\s_{\B_{u,s}}(f)\!=\!\s(f).$  Thus since $0\!\not\in\!\s(f),$ we know that $0\!\not\in\!\s_{\B_{u,s}}(f).$  Now since $f$ is invertible in $\Ll_v,$ the inverse of $f$ must share the defining property of $\Ll_v,$ namely \eqref{defLv}.  To extend the result to all $f\in\Ll_v,$ we apply the same argument as in the proof of Theorem \ref{big1}.
\end{myproof}
\end{corollary}

We also note that if $f\in\Ll^1_v$ or $f\in\Ll_v,$ then in the case that $f$ is a positive operator, all powers $f^\a,\,\a\in\R$ are in $\Ll^1_v$ or $\Ll_v,$ respectively, by the Riesz functional calculus (see e.g., \cite{Rud91}).\\

\subsubsection{Template of an example}
As an illustration of how a specific example of an operator in $\Ll^1_v$ or $\Ll_v$ might be discovered, we investigate the convolution operator generated from $f(y)=\d(y)+2^{\frac{1}{4}}\,e^{-\pi y^2}.$  Note that we have $\wh{f}=1+\widetilde{f}.$
Let $W=W(\R^d)$ denote the Wiener algebra, for which $g\in W\Rightarrow g,\widehat{g}\in\Lone,$  and fix $h\in W.$  Then
\bea
\F\left(\frac{\wh{h}}{\wh{f}}\right)&=&\int_{\R^d}(1+\wt{f}(\w))^{-1}\left[\int_{\R^d}h(y)e^{-2\pi i\w\cdot y}dy\right]e^{2\pi i w\cdot x}d\w\allowdisplaybreaks\\
&=&\int_{\R^d}h(y)\int_{\R^d}(1+\wt{f}(\w))^{-1}e^{2\pi i w\cdot (x-y)}d\w\,dy\allowdisplaybreaks\\
&=&\int_{\R^d}h(y)\F^{-1}\left(\frac{1}{1+\wt{f}(\w)}\right)(x-y)\,dy\allowdisplaybreaks\\
&=&\int_{\R^d}h(y)\F^{-1}\left(1-\frac{\wt{f}(\w)}{1+\wt{f}(\w)}\right)(x-y)\,dy\allowdisplaybreaks\\
&=&\int_{\R^d}h(y)\F^{-1}\left(1-\frac{1}{e^{\pi\w^2}+1}\right)(x-y)\,dy\allowdisplaybreaks\\
&=&\int_{\R^d}h(y)\left[\d(x-y)-\F^{-1}\left(\frac{1}{e^{\pi\w^2}+1}\right)(x-y)\right]\,dy
\end{eqnarray}
$\Rightarrow$ the non-unit part of the kernel is $\wt{f^{-1}}(x,y)=-\F^{-1}\left(\frac{1}{e^{\pi\w^2}+1}\right)(x-y).$  Since for all $\w\in\R^d,$ we have
\be
\frac{e^{-\pi\w^2}}{2}\leq\frac{1}{e^{\pi\w^2}+1}<e^{-\pi\w^2}
\end{equation}
we conclude that $\F^{-1}\left(\frac{1}{e^{\pi\w^2}+1}\right)$ shares the same decay as $2^{\frac{1}{4}}\,e^{-\pi\w^2}=\F^{-1}\left(2^{\frac{1}{4}}\,e^{-\pi\w^2}\right).$  Of course $e^{-\pi\w^2}$ is a member of neither $\Ll^1_v$ nor $\Ll_v,$ but finding an explicit example in $\Ll^1_v$ or $\Ll_v$ has thus far proven elusive.

\np
\pagestyle{myheadings}
\chapter{Pseudodifferential Operators and Linear Time-Varying Systems}
\thispagestyle{myheadings} 
\markright{  \rm \normalsize CHAPTER 4. \hspace{0.5cm}
 PDOs AND LTVs}

The results in Chapter \ref{decaydistributions} should find application in the area of wireless communications channels.  Let $s(t)$ be a wireless signal transmitted through a channel ${\mathbf H},$ and let $s(t-\t)$ denote its $\t$-delayed version.  The received signal is
\begin{equation}
r(t) = ({\mathbf H}s)(t)+\e(t)
= \int_\R  h(t,\t) s(t-\t) d\t +\e(t),
\label{channeldecay}
\end{equation}
where $h(t,\t)$ is the impulse response of the time-varying channel
${\mathbf H}$ \cite{Bel63} and $\e(t)$ is zero-mean stationary white noise (\ref{channeldecay} will be explained fully in Chapter \ref{OBFDMinterpol}).  A typical channel may have a line-of sight component, and cause many reflections from objects in the vicinity of the transmitter and/or receiver (as well as Doppler effects from the relative motion between the transmitters and receivers).  Thus ${\mathbf H}$ will yield an $h$ which consists of the sum of $\d(t-\t)$ and some other function of $t$ and $\t$ possessing off-diagonal decay, quite possibly granting $h$ membership in $\Ll^1_v$ or $\Ll_v.$  The wireless channel applications of the work in Chapter \ref{decaydistributions} will be the subject of future investigation.\\

There is an interesting connection between pseudodifferential operators and linear time-varying systems, of which a wireless channel is one example (see \cite{Str04}), that will be instructive to describe as a point of departure for our excursion into signal design for wireless communications.  

\section{Pseudodifferential operators}
We will assume that all functions described in this section are in $\S(\R^d)$ so that the integrals described herein behave; i.e., we interpret all expressions in the distributional sense.  In the study of partial differential equations one considers equations of the form
\be
Af(x)=\sum_{|\a|\leq N}\s_\a(x)D^\a f(x)=g(x)\label{pde1}
\end{equation}
where $N$ is the order of $A,\,\a$ is a multi-index, and $\{\s_\a\}$ is the set of non-constant coefficients.  We have
\be
(D^\a f)^{\hat{ }}(\w)=(2\pi i\w)^\a\hat{f}(\w)\label{derivfourier}
\end{equation}
(see section \ref{tempereddistributions}).  Using \eqref{fourierinv} and \eqref{derivfourier}, one can write \eqref{pde1} as
\bea
Af(x)&=&\int_{\R^d}\biggl(\sum_{|\a|\leq N}\s_\a(x)(2\pi i\w)^\a \biggr)\hat{f}(\w)e^{2\pi i x\cdot\w}\,d\w\label{pde2}\\
&=&\int_{\R^d}\s(x,\w)\hat{f}(\w)e^{2\pi i x\cdot\w}\,d\w\label{pde3}
\end{eqnarray}
where $\s$ is the {\em symbol} of $A.$  One can admit more general symbols, yielding the definition of a pseudodifferential operator.
\begin{definition}\label{pdodef}
Let $\s$ be a measurable function or a tempered distribution on $\R^{2d}.$  Then define
\be
K_\s f(x):=\int_{\R^d}\s(x,\w)\hat{f}(\w)e^{2\pi ix\cdot\w}\,d\w.\label{pdo1}
\end{equation}
The operator $K_\s$ is called the {\em pseudodifferential operator} with symbol $\s.$
\end{definition}
With suitable restrictions on $\s,\,K^\s$ becomes a multiplication operator, or a Fourier multiplier, or a differential operator with constant coefficients.  Additionally, we can write $K^\s$ in several ways.  First, use \eqref{fourierinv} in \eqref{pdo1} to yield
\bea
K_\s f(x)&=&\int_{\R^d}\left(\int_{\R^d}\s(x,\w)e^{-2\pi i(y-x)\cdot\w}\,d\w\right)f(y)\,dy\label{pdoasconv}\\ 
&=&\int_{\R^d}k_\s(x,y)f(y)\,dy
\end{eqnarray}
where we deonte by $k_\s$ the kernel of the integral operator $K_\s.$  Note then that the results of Chapter \ref{decaydistributions} apply to pseudodifferential operators.  Denote by $L_a$ the coordinate transformation $L_aF(x,y)=F(x,y-x)$ and let $\F_2$ be the partial Fourier transform in the second variable.  Then $k_\s$ can be written as 
\be
k(x,y)=\F_2\s(x,y-x)=L_a\F_2\s(x,y).
\end{equation}
We can also interpret \eqref{pdoasconv} as the convolution of $f(y)$ with $\F_2^{-1}\s(x,y)$ from
\be
K_\s f(x)=\int_{\R^d}\F_2^{-1}\s(x,x-y)f(y)\,dy.\label{pdoasconvwinvfour}
\end{equation}

\section{Linear time-varying systems}
A {\em linear system} transforms an input signal $f$ into an output signal $Af$ for which $A(c_1f_1+c_2f_2)=c_1Af_1+c_2Af_2$ for all inputs $f_j$ and constants $c_j.$  For a time-invariant system, we have $A(T_xf)=T_xAf\,\,\forall\,\,x\in\R^d.$  Any linear operator satisfying this commutation relation can be written as a convolution operator
\be
Af=f*h\quad\mbox{where}\quad h=A\d.
\end{equation}
The frequency response of the system $A$ is
\be
A(e^{2\pi i\w\cdot t})(x)=\int_{\R^d}h(t)e^{2\pi i\w\cdot(x-t)}\,dt=\hat{h}(\w)e^{2\pi i\w\cdot x}
\end{equation}
and we conclude that the complex exponentials $e^{2\pi i\w\cdot x}$ are the eigenfunctions of $A$ with eigenvalues $\hat{h}(\w).$\\

A linear system $A$ may change with time, in which case $A$ is called a time-varying system.  In this case the convolver $h$ depends on time, and we have
\be
Af(x)=\int_{\R^d}h(x,x-y)f(y)\,dy.\label{tv1}
\end{equation}
By \eqref{pdoasconvwinvfour} we can interpret the linear time-varying system $A$ as a pseudodifferential operator $K_\s$ with the symbol $\s=\F_2h.$\\

In general one has $K^{-1}_\s\neq K_{1/\s},$ but for many pseudodifferential operators, in particular underspread operators, where the the channel does not spread the signal ``too much'' in the time-frequency plane, (see Chapter \ref{OBFDMinterpol}), one has $K^{-1}_\s\approx K_{1/\s}$ (see \cite{KM98}), and this is the case for a typical wireless communications channel.  Underspread operators have been studied, but these have almost always been Hilbert-Schmidt operators, which are not invertible.  This is too restrictive for the case of wireless communications, since it excludes the case where the channel operator is just a small perturbation of the identity.\\

A wireless communications channel is a linear time-varying system, and can thus be interpreted as a pseudodifferential operator.  Denote by $\Hchan_\s$ the channel whose characteristics are described by the symbol $\s.$  Strohmer proved in \cite{Str04} that if $\s$ is a wireless channel operator symbol, then $\s\in M^{\infty,1}_m,$ the space of symbols $\s$ for which
\be
\int_{\R^{2d}}\underset{X}{\sup}V_\Psi\s(X,\O)m(\O)\,d\O<\infty\label{sigweight}
\end{equation}
where $X,\O\in\R^{2d},$ $m$ is admissible in the sense of Section \ref{weights}, and $V$ denotes the $2d$-dimensional STFT.\\

In the case that the channel is underspread and the pulses used for encoding the data have good time-frequency localization or TFL (explained more rigorously in Chapter \ref{OBFDMinterpol}; for now we will simply say that a pulse has good TFL if it possesses good decay in both the time domain and frequency domain), we have $\Hchan^{-1}_\s\approx \Hchan_{1/\s}$ in the weak sense, meaning that we can trade invertibility for approximate invertibility and significantly less effort.  \\

In the sequel we shift the focus from the operator to the functions which make up the components of our signal.  To ensure that the spreading induced by a given channel operator is minimized, we consider, and make some demands of, the TFL of the pulses on which the data is encoded and sent through the channel.\\

In this regard, we have developed several improvements to a communications scheme called orthogonal frequency-division multiplexing, or OFDM.  This system not only shows great promise as a candidate for future wireless communications systems, but possesses a beautiful mathematical foundation as well.  The remainder of this thesis will be concerned with our improvements to the algorithm currently used in practical applications of OFDM.

\np
\pagestyle{myheadings}
\chapter{OFDM}\label{OBFDMinterpol}
\thispagestyle{myheadings} 
\markright{  \rm \normalsize CHAPTER 5. \hspace{0.5cm}
 OFDM}

\section{Introduction}
OFDM is an efficient technology for wireless data transmission \cite{Zou95}, and is currently used in the European digital audio broadcasting standard \cite{CFH89}, in digital terrestrial TV broadcasting, broadband indoor wireless systems \cite{Macedo97}, and an increasing number of wireless internet service providers (WISPs) in the U.S. \cite{Sel04}.  The advantages of OFDM include robustness against multipath fading and high spectral efficiency.\\

Transmission over wireless channels can be subject to multipath propagation (i.e., that the received signal is the sum of perhaps a line-of-sight component and several reflected delayed components), leading to time dispersion which causes intersymbol interference (ISI), and the Doppler effect (frequency compression or rarefaction due to the constancy of the speed of light and the relative motion between the transmitter and receiver) which causes dispersion in the frequency domain leading to interchannel interference (ICI)\cite{Rap96}.  The distortion resulting from time-frequency dispersive channels depends crucially on the TFL of the pulse.  For instance OFDM systems based on rectangular-type pulses with guard-time interval or cyclic prefix can prevent ISI, but can not effectively combat ICI, since they decay too slowly in the frequency domain (see Figure \ref{fig:OFDMtones}).\\

The design of time-frequency well-localized pulse-shaping OFDM filters is therefore an active area of research \cite{FAB95,VH96,HB97,KM98,Bol02}.  Some of the proposed pulse shape design methods involve the solution of highly non-linear optimization problems \cite{VH96,HB97}, other approaches lead to pulse shapes that can be easily computed, but little or no theoretical justification of the optimality of the so-obtained pulse shapes has been given \cite{FAB95,Bol02,KM98}. 

\section{OFDM basics}\label{OFDMbasics}
The principal idea of OFDM is to divide the available spectrum into several subchannels (subcarriers) so that over each subcarrier the data rate is slower, the symbol period is longer and thus the system is more tolerant of time dispersion.\\

\begin{figure}[!h]
\begin{center}
\epsfig{file=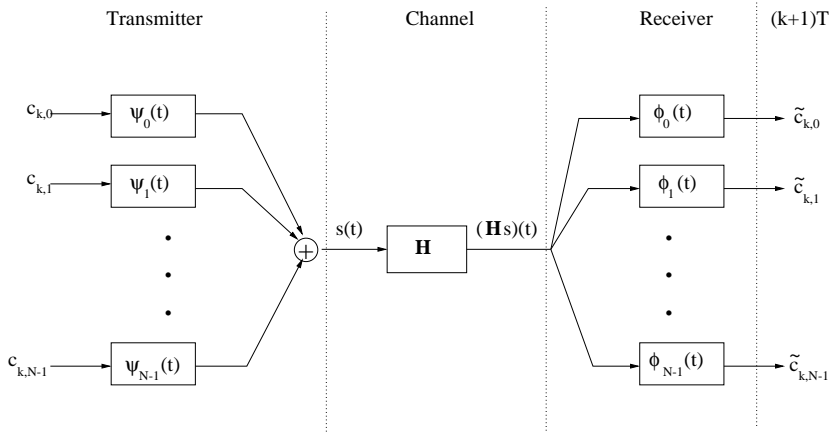,width=110mm}
\caption{A general baseband OFDM system configuration}
\label{fig:ofdm}
\end{center}
\end{figure}

The transmission functions of an OFDM system consist of translations and modulations of a single pulse shape $\psi$, i.e.,
\begin{equation}
\label{ofdmsystem}
\psi_{kl}(t)=\psi(t-kT) e^{2\pi i \l Ft}, 
\end{equation}
for some $T,F>0, k, \in \Z$ and $\l = 0,\dots,N-1$ or $\l \in \Z$.
We denote such a function system (whether orthogonal or not) generated via \eqref{ofdmsystem} by the triple $(\psi,T,F).$  We call $T$ the {\em symbol period} and $F$ the {\em carrier separation}.  Let $W$ be the given bandwidth, then $N_c=W/F$ is the number of subcarriers. Throughout our analysis it is convenient to assume $N=\infty,$ so the frequency index $\l$ in \eqref{ofdmsystem} runs through the integers instead of from $0$ to $N-1$.  But the results presented in this paper also hold for finitely many subcarriers, either verbatim or with a negligible approximation error.  A baseband OFDM system is schematically represented in Figure \ref{fig:ofdm}. 

A necessary condition for perfect reconstruction of the data is that the functions $\psi_{k,\l}$ are linearly independent.  This implies that $TF\ge 1$ which follows easily from the connection between OFDM systems and Gabor frames, e.g., see Section \ref{relevthmgabor} of this paper or \cite{KM98}.  We note that orthogonality of the functions $\psi_{k,\l}$ is not a requirement for perfect reconstruction, but minimizes the error caused by AWGN.\\

Let $c_{k,\l},\,k,\l\in \Z,$ be the data symbols to be transmitted. The OFDM signal is then given by
\begin{equation}
\label{OFDMsignal}
s(t)= \sum_{k,\l} c_{k,\l} \psi_{k,\l}.
\end{equation}
The received signal is
\begin{equation}
r(t) = ({\mathbf H}s)(t)+\e(t)
= \int_{\R^d}h(t,\t) s(t-\t) d\t +\e(t),
\label{eq5}
\end{equation}
where $h(t,\t)$ is the impulse response of the time-varying channel ${\mathbf H}$ \cite{Bel63} and $\e(t)$ is zero-mean stationary\footnote{A stationary stochastic process is defined by the statistical invariance of its properties under a shift of the time-origin.} white noise.  The data symbol $c_{k,\l}$ is transmitted at the $(k,\l)^{th}$ point of the rectangular time-frequency lattice $(kT,\l F), k,\l \in \Z$. Assuming that the $\psi_{k,\l}$ are orthonormal the transmitted data are recovered by computing the inner product of the received signal $r$ with the functions $\psi_{k,\l}$, i.e.,
$\tilde{c}_{k,\l}=\langle r ,\psi_{k,\l} \rangle$.  By \eqref{OFDMsignal} and \eqref{eq5} we have
\begin{equation}
\tilde{c}_{k,\l}=\langle \Hchan s ,\psi_{k,\l} \rangle + 
\langle \e,\psi_{k,\l} \rangle= \sum_{k',\l'} c_{k',\l'}\,\langle \Hchan \psi_{k',\l'},
\psi_{k,\l} \rangle + \langle \e,\psi_{k,\l} \rangle,\!\!\!\!\!\!
\label{innprod}
\end{equation}
where the interchange of summation and integration is justified under
mild smoothness and decay conditions on $\psi$ and $\Hchan$.
In the ideal (but unrealistic) case the wireless channel does not
introduce any distortion. Thus the orthogonality between the functions
$\psi_{k,\l}$ is preserved, and if $\e=0$ the data $c_{k,\l}$ can be recovered 
exactly. However in practice wireless channels are time-dispersive and 
frequency-dispersive. Hence when $\psi_{k,\l}$ passes through the time-varying 
channel $\Hchan$ it is spread out in time and frequency, causing 
ISI and/or ICI in the transmitted data. It is customary 
(cf.~e.g. \cite{FAB95,Rap96,HB97,KM98}) to describe the channel by the 
{\em maximal excess delay} $\t_0$ and the {\em maximal Doppler spread} 
$\nu_0$ of the channel scattering function.  It is clear from \eqref{innprod} that the interference of the data symbol 
$c_{k,\l}$ with neighbor symbols will be smaller with better 
TFL of $\psi_{k,\l}$ and with larger distance 
between $c_{k,\l}$ and adjacent data symbols in the time-frequency domain.\\

As shown in \cite{HB97} for a large number of subcarriers the {\em spectral efficiency} $\r$ of the OFDM system in terms of data symbols per second per Hertz\footnote{In the remainder of this paper we will use $\r$ to denote spectral efficiency, since it is standard notation, and there will be no confusion with the spectral radius $\r(f)$ of Chapter \ref{decaydistributions}.} is approximately given by $\rho=(TF)^{-1}$.  Since $TF \ge 1$ we have of course that the maximal spectral efficiency is $\r_{\max}=1$.  Thus any increased distance between adjacent symbols in the time-frequency domain costs the designer in terms of spectral efficiency.\\

\section{Dispersion in the wireless channel}\label{multnDoppl}

\subsection{The Doppler effect and ICI}
The formula for the Doppler shift is given by
\be
f_d=\frac{v}{\la}\text{cos}(\th)
\end{equation}
where $f_d$ is the change in frequency of the signal due to $v,$ the relative velocity of the transmitter and receiver and $\la,$ the wavelength of the signal (for simplicity, we assume a pure frequency), and $\th$ is the instantaneous angle of approach.  For a speed of 1100 mph ($\approx$492 m/s, the approximate relative velocity of two approaching airliners), on direct approach ($\th=0$), for a 5 GHz signal ($\la=6$ cm) the Doppler shift is then 8.196 kHz.  But for IEEE 802.11a, the subcarrier spacing is 312.5 kHz.  So how can such a small relative shift in frequency lead to significant ICI at the receiver?\\

\begin{figure}[!h]
\begin{center}
\epsfig{file=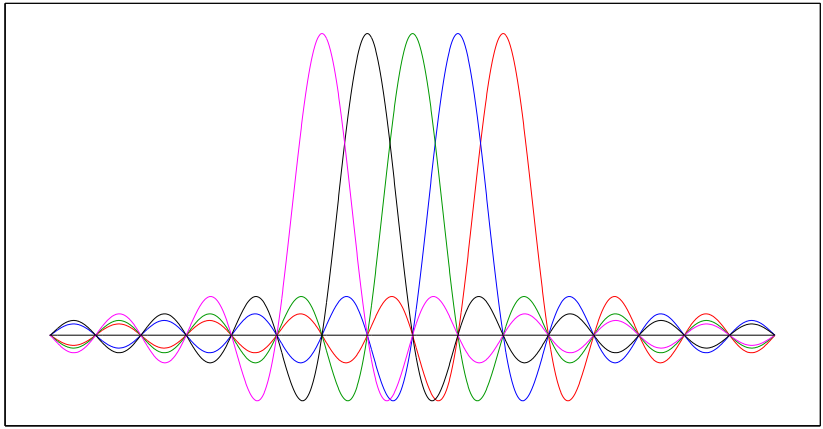,width=110mm}
\caption{Frequency domain representation of five standard (square-wave-carried) OFDM tones.  Note that the peak of each tone is at a zero of every other tone.}
\label{fig:OFDMtones}
\end{center}
\end{figure}

Figure \ref{fig:OFDMtones} illustrates the frequency domain representation of the modulated square-wave pulses used in standard OFDM.  Note that the peak amplitude of each frequency is located at the zeros of all other frequencies.  This is what distinguishes OFDM from FDM, a system wherein the designer precludes interference only by maintaining a large guard interval between frequencies, which can reduce spectral efficiency dramatically.  OFDM requires no frequency guard interval, since ideally the subcarriers do not interfere with each other.  But as one might suspect, this system is quite sensitive to frequency perturbations.  When the receiver samples at the center frequency of each subcarrier, it gathers only energy from the desired signal, provided that there is no frequency dispersion.  This is referred to as zero-delay orthogonality.  It is clear from Figure \ref{fig:OFDMtones} that it is not necessary for the Doppler shift to be anywhere near half of the subcarrier separation to cause significant interference.  One cannot speculate as to the relative velocities which communications systems might encounter in the future.

\subsection{Multipath delay and ISI}
\begin{figure}[!h]
\begin{center}
\epsfig{file=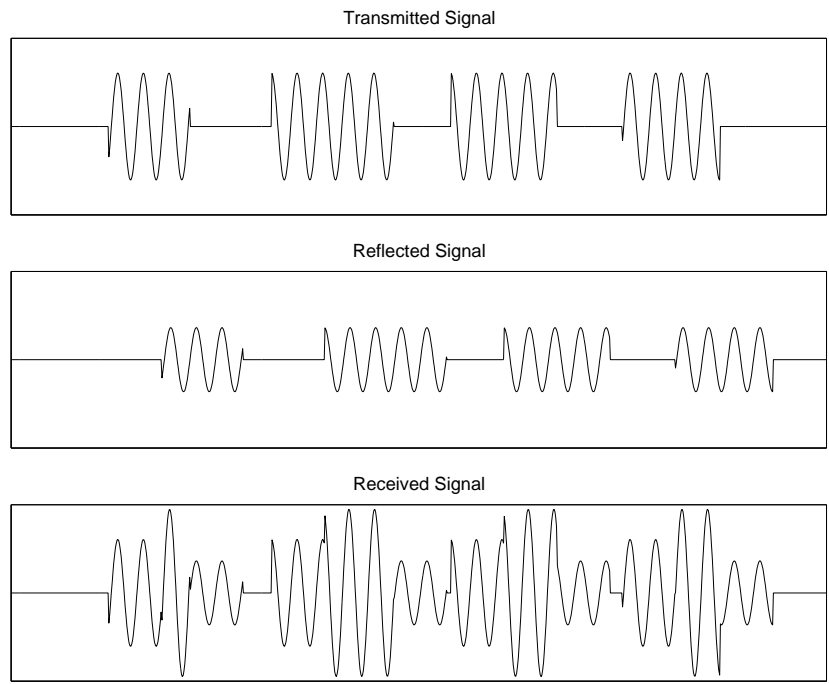,width=110mm}
\caption{Multipath delay spread affects a signal in a reflective environment.}
\label{fig:delay}
\end{center}
\end{figure}

A typical received signal consists of (perhaps) a direct line-of-sight signal from transmitter to receiver, plus reflected signals arriving at later times due to the necessarily greater path lengths travelled, spreading the signal energy in time.  Figure \ref{fig:delay} shows this effect.

\section{Time-frequency localization}
Ideally we would like to construct an OFDM system $(\psi,T,F)$ that 
satisfies the following three requirements simultaneously:
\begin{itemize}
\item[\quad (i)] the functions $\psi_{k,\l}$ are mutually orthogonal,
\item[\quad (ii)] $\psi$ is well-localized in time and frequency domain,
\item[\quad (iii)] $TF=1$, i.e., maximal spectral efficiency.
\end{itemize}

As is well-known from Gabor theory these three conditions cannot be
satisfied simultaneously due to the Balian-Low Theorem \cite{Dau92} the typical statement of which precludes, in essence, simultaneous satisfaction of (ii) and (iii).  It should be noted, however, that this holds true, even if we relax condition (i) and consider biorthogonality instead (which however would increase the sensitivity of the system
to AWGN).\\

Hence we relax condition (iii) and allow $TF>1$. This loss in spectral efficiency is usually an acceptable price to pay to mitigate ISI/ICI (provided $TF$ is not too large).  As pointed out above in a dispersive channel the transmitted signal is spread in time and in frequency, causing ISI and ICI. This interference can be reduced if the pulse shape possesses good TFL. The classical way to measure the TFL involves the Heisenberg Uncertainty Principle. More precisely, if $f\in\LtR$ and $\tau,\omega\in\R$ are arbitrary,
\begin{equation}
\left(\int_{\R} (t-\tau)^2 |f(t)|^2 dt\right)^{\frac{1}{2}}\left(\int_{\R} (\nu-\omega)^2 |\hat{f}(\nu)|^2 d\nu\right)^{\frac{1}{2}}\geq\frac{\|f\|^2}{4\pi}. 
\label{up0}
\end{equation}
Equality holds if and only if $f$ is of the form $f(t)=ce^{2\pi i \omega(t-\tau)}
e^{-\pi \alpha (t-\tau)^2}$ for $\tau, \omega \in \R, c\neq 0,$ and $\alpha>0$, i.e.,
if $f$ is a (possibly translated and/or modulated and/or dilated) Gaussian. \\

In \eqref{up0} time and frequency are treated separately. It seems more appropriate in our context, in particular with respect to time-frequency dispersive channels, to measure the TFL by considering time and frequency simultaneously.  Moreover in light of \eqref{innprod} we are naturally interested in measuring the cross-interference between different transmitter functions $\psi_{k,\l}$ and $\psi_{k',l'}$.  Thus we replace the pair $(f, \hat{f})$ by a joint time-frequency representation of two functions $f$ and $h$, such as the cross-ambiguity function $A(f,h)$.\\

In light of the Heisenberg Uncertainty Principle it is not surprising
that $A(f,h)$ cannot be arbitrarily well concentrated.
This can be expressed for instance by the fact that $A(f,h)$ cannot
be compactly supported, see \cite{Jan98}. An alternative formulation
is the following Heisenberg-type inequality \cite{Coh01,BDJ03}.
\begin{theorem}
\label{th:amb}
For every $f, h \in \LtR$ there holds
\begin{gather}
\int_{\R}\int_{\R}
(t-\tau)^2 |A(f,h)(t,\nu)|^2 dt \,d \nu \cdot
\int_{\R}\int_{\R}
(\nu-\omega)^2 |A(f,h)(t,\nu)|^2 dt \,d \nu \ge 
\frac{\|f\|^4 \|h\|^4}{4\pi^2}. \label{up2a} 
\end{gather}
Equality holds if and only if 
\begin{gather}
f(t)=c_1e^{2\pi i\beta t} e^{-\pi \alpha (t-s)^2},
\notag \\
h(t)=c_2e^{2\pi i(\omega+\beta)t}e^{-\pi \alpha (t-\tau-s)^2},
\label{up2min}
\end{gather}
with $\alpha >0,c_1,c_2\neq 0,$ and $\beta, s,\tau \in \R$.
\end{theorem}

Thus again, translated, modulated, and dilated versions of the Gaussian lead to optimal TFL.  But Gaussians cannot be orthogonal under translation, so we appeal to a technique first discovered in the finite-dimensional discrete case by the chemist L\"owdin and here termed {\em L\"owdin orthogonalization}.  By applying an infinite-dimensional continuous version of L\"owdin's technique to form orthogonal pulses from the Gaussian, we retain as much of the Gaussian's optimal TFL as is theoretically possible.

\np
\pagestyle{myheadings}
\chapter{Frames, Gabor frames, and L\"owdin orthogonalization}\label{FrmnGabFrm}
\thispagestyle{myheadings} 
\markright{  \rm \normalsize CHAPTER 6. \hspace{0.5cm}
 FRAMES, GABOR FRAMES, AND L\"OWDIN ORTHOGONALIZATION}

To construct an OFDM communications system which is robust against both ISI and ICI, one must begin with a function $h$ possessing good decay in both time and frequency, and generate a set $\{h_{k,\l}\}$ of functions by integral time- and frequency-shifts of $g.$  This setup is almost identical to a general Gabor system (defined in Section \ref{gabfr}), with the exception of the phase factor present in a Gabor system due to the application of modulation first, and then translation to the Gabor system generator, but with the additional constraint that every function in the sequence $\{h_{k,\l}\}$ be orthogonal to all other members of $\{h_{k,\l}\}.$  Since any system $\{g_{k,\l}\}$ generated from the Gaussian $g$ obviously does not yield an orthogonal set of transmission functions defined in this manner, and nor do any dilations of the Gaussian\footnote{Since for all $\s>0\,,\,g^\s(t)=(2\s)^{1/4}e^{-\pi\s t^2}>0$ for all $t$, $\ggam(t)$ can never by orthogonal to any of its shifted copies, a similar argument holds of course for modulations of $\ggam.$}, we use a standard method to ``orthogonalize'' the set $\{g_{k,\l}\}.$  In connection with pulse shape design it has been used implicitly 
in \cite{FAB95} (see also \cite{BJ00}), and explicitly in \cite{BDH00,Bol02,Str01}.  We describe this algorithm first in finite dimensions as an illustration of its elegance.

\section{Discrete L\"owdin orthogonalization in finite dimensions}
Suppose $\Phi\,\in\C^{m\times n},$ with $m\geq n$ and rank$(\Ph)=n.$  Consider a nonsingular linear transformation $A\in\C^{n\times n}$ which, when applied to $\Phi,$ yields an orthonormal basis $\Ps= \Ph A$ of $\C^n.$  We have
\be
I_{n\times n}=\Ps^*\Ps=(\Ph A)^*(\Ph A)=A^*\Ph^*\Ph A,
\end{equation}
and we define $\Delta\in\C^{n\times n}$ as the Gram matrix $\Ph^*\Ph.$  By the properties of $\langle\cdot,\cdot\rangle,\,\Delta$ is Hermitian, and by Proposition 5.1.5 in \cite{Gro01}, $\Delta$ is a positive invertible operator.  Consequently, $\exists \,\,U$ unitary which diagonalizes $\Delta := U^*\La U,$ and the elements of $\La$ are positive.  We can define an arbitrary function $f$ on $\Delta$ via 
\be
f(\Delta)=U^* f(\La)U,\label{arbfcnonDel}
\end{equation}
provided that the Taylor series for $f(\La)$ converges, since \eqref{arbfcnonDel} is always true for polynomial $f.$  In particular, for $f(x)=x^{-1/2},$ we note that $\Delta$ is nonsingular, so the Taylor series for $\Delta^{-1/2}$ converges just as the series for $f(x)=x^{-1/2}$ converges for $x>0.$  I.e., if $x\in(0,\infty],$ then since the mapping $f:x\rightarrow f(x)$ is continuous and infinitely differentiable on $(0,\infty],$ the Taylor series for $f$ converges at $x.$  Thus we have $\Delta^{-1/2} = U^*\La^{-1/2}U,$ and $\Delta^{-1/2}$ is Hermitian as well.\\

Define $B=\Delta^{1/2}A,$ so that $A=\Delta^{-1/2}B$ and $A^*=B^*(\Delta^{-1/2})^*=B^*\Delta^{-1/2}.$  Thus
\begin{equation}
A^*\Delta A=B^*\Delta^{-1/2}\Delta\Delta^{-1/2}B=B^*B
\end{equation}
and we see that the equation $A^*\Delta A=I$ has a solution if $B$ is an arbitrary unitary matrix.  The choice of $B=I$ leads to {\it L\"owdin} or {\it symmetric} orthogonalization (or orthonormalization): $\Psi=\Phi\Delta^{-1/2}.$\\

An important property of L\"owdin orthogonalization is that of all possible orthogonalizations it disturbs the original basis the least under the Frobenius norm, defined by
\be
\|\Phi\|_F=\sum_{i=1}^m\sum_{j=1}^n |a_{ij}|^2.\nn
\end{equation}

\begin{theorem} Let $\Ph$ be a linearly independent set of $n$ vectors in $\C^m.$  Considering $\Ph$ as an $m\times n$ matrix of complex numbers, we have the singular value decomposition $\Ph=UDV^*,$ where $U\in \C^{m\times n},\, V\in\C^{n\times n},$ $U$ and $V$ are unitary, and $D\in\C^{n\times n}$ contains the singular values of $\Ph.$  Then if $Q$ is unitary, $\|\Ph-Q\|_F$ is minimized by $Q=P:=\Ph(\Ph^*\Ph)^{-1/2}=UV^*.$\\

\begin{myproof}  The strategy is to show first that $P=UV^*$ and thus that $P$ is unitary, and then to show that $P$ minimizes the above expression.
\bea
\Ph = UDV^*\,\Rightarrow \Ph(\Ph^*\Ph)^{-1/2} &=& UDV^*(VDU^*UDV^*)^{-1/2}\\
&=& UDV^*(VD^2V^*)^{-1/2}\nonumber\\
&=& UDV^*(VDV^*VDV^*)^{-1/2}\nonumber\\
&=& UDV^*(VDV^*)^{-1}\nonumber\\
&=& UDV^*VD^{-1}V^*\nonumber\\
&=& UV^*\nonumber
\end{eqnarray}
Thus $P=UV^*$ and $P$ is unitary.  To show that $Q=P$ minimizes $\|\Ph-Q\|_F$ among unitary $Q,$ note that\\
$$\|\Ph-Q\|_F=\|UDV^*-Q\|_F=\|D-U^*QV\|_F,$$
which is minimized if $U^*QV=I,$ or $Q=UV^*=P.$
\end{myproof}
\end{theorem}
Having treated the discrete finite case, we turn to the continuous infinite case, beginning with an introuction to frames.

\section{Frames}
There are many ways to reconstruct $\Ltwo$-functions using sets of functions in $\Ltwo.$  We can use wavelets, or even more exotic sets of functions.  Orthonormal bases are usually ideal in some sense, but often enough we have no stake in uniqueness of coefficients.  This particular indifference leads to the use of what are called {\it frames}.  A sequence $\{e_j\}_{j\in J},$ over some index set $J,$ in a Hilbert space ${\mathcal H},$ is a frame for ${\mathcal H}$ if and only if
\be
\label{frm}
\exists\, A,B\in\R_+\ni\forall f\in {\mathcal H} \qquad A\|f\|^2\leq\underset{j\in J}{\sum}|\langle f,e_j\rangle |^2\leq B\|f\|^2.
\end{equation}
For any subset $\{e_j\}_{j \in J} \subset {\mathcal H}$, we have the {\it analysis} and {\it synthesis operators}
\begin{eqnarray*}
Cf &=& \{\langle f,e_j \rangle |\,j \in J \} \\
Dc &=& \begin{Huge}\sum\end{Huge}_{j \in J}c_j e_j.
\end{eqnarray*}
It can be easily shown that $D=C^*.$  The action of the {\it frame operator} $S=C^*C=DD^*$ on $f\in {\mathcal H}$ is
\be
Sf = \sum_{j\in J}\langle f,e_j\rangle e_j.
\end{equation}
The condition $A\!>\!0$ ensures that $S$ is injective; the condition $B\!<\!\infty$ ensures that $S$ is bounded and thus continuous (see any elementary text discussing operators on Banach spaces).  Since $\langle Sf,f\rangle=\langle C^*Cf,f\rangle=\langle Cf,Cf\rangle\geq 0,\,S$ is a positive operator.  If $A=B$ then $\{e_j\}$ is a {\it tight} frame.  Tight frames are important to us because they can lead to orthonormal systems, which we will see later.  

\begin{lemma}\label{tgtfrmS12}
If $\{e_j\,|\,j\in J\}$ generates a frame, then the set $\{S^{-\frac{1}{2}}e_j\,|\,j\in J\}$ generates a tight frame.\\
\begin{myproof}$\,$ (From \cite{Gro01})
\bea
&&f=S^{-\frac{1}{2}}SS^{-\frac{1}{2}}f=\sum_{j\in J}\langle f,S^{-\frac{1}{2}}e_j\rangle S^{-\frac{1}{2}}e_j\\
&\Rightarrow&\qquad\langle f,f\rangle=\sum_{j\in J}|\langle f,S^{-\frac{1}{2}}e_j\rangle|^2.
\end{eqnarray}
\end{myproof}
\end{lemma}

We will refer to frames $\{d_j\,|\,j\in J\}$ and $\{e_j\,|\,j\in J\}$ for which 
\bea
f&=&\underset{j\in J}{\sum}\langle f\,,\,e_j\rangle d_j\label{dual1}\\
&=&\underset{j\in J}{\sum}\langle f\,,\,d_j\rangle e_j\label{dual2}
\end{eqnarray}
as {\em dual frames with respect to the identity}, or simply {\em dual frames}.  If
\bea
S^\a f&=&\underset{j\in J}{\sum}\langle f\,,\,e_j\rangle d_j\label{dual3}\\
&=&\underset{j\in J}{\sum}\langle f\,,\,d_j\rangle e_j\label{dual4}
\end{eqnarray}
then we refer to $\{d_j\}$ and $\{e_j\}$ as {\em dual frames with respect to $S^\a$}.

\subsection{Powers of the frame operator}
Assume that $\{e_j\}_{j\in J}$ generates a frame.  Consider the nonorthogonal expansions 
\bea
f&=&\underset{j\in J}{\sum}\langle f\,,\,S^{-\frac{1}{2}}e_j\rangle S^{-\frac{1}{2}}e_j\label{tgtidentity}\\
&=&\underset{j\in J}{\sum}\langle f\,,\,e_j\rangle S^{-1}e_j\label{dualidentity1}\\
&=&\underset{j\in J}{\sum}\langle f\,,\,S^{-1}e_j\rangle e_j.\label{dualidentity2}
\end{eqnarray}\nonumber\\
We can use these to derive several results about the frame operator $S^{ ^{\scriptstyle p}}.$
\begin{itemize}
\item  For which $p$ does $S^pe_j$ generate a frame?  It follows from Appendix A.7 in \cite{Gro01} that $\forall\,p\!\in\!\R,\,S^{ ^{\scriptstyle p}}e_j$ generates a frame.
\item  What are the frame elements which generate the frame operator $S^{ ^{\scriptstyle q}}$?  If we operate on $f$ in \eqref{tgtidentity}, we have
\bea
S^{ ^{\scriptstyle q}}f&=&S^{ ^{\scriptstyle r}}(S^{ ^{\scriptstyle {q-r}}}f)\nn\\
&=&\underset{j\in J}{\sum}\langle S^{ ^{\scriptstyle r}}f\,,\,S^{-\frac{1}{2}}e_j\rangle S^{ ^{\scriptstyle {q-r}}}(S^{-\frac{1}{2}}e_j)\nn\\
&=&\underset{j\in J}{\sum}\langle f\,,\,S^{ ^{\scriptstyle r-\frac{1}{2}}}e_j\rangle S^{ ^{\scriptstyle q-r-\frac{1}{2}}}e_j\label{alsodual}
\end{eqnarray}
where we used the fact that powers of $S$ are self-adjoint, and the expression \eqref{tgtidentity} (as well as \eqref{dualidentity1} and \eqref{dualidentity2}) converge unconditionally.  We conclude that the $S^{\frac{q-1}{2}}e_j$ generate the frame for $S^{ ^{\scriptstyle q}}\,(r=q/2),$ and that the $S^{ ^{\scriptstyle r-\frac{1}{2}}}e_j$ and $S^{ ^{\scriptstyle {q-r}}}e_j$ generate frames which are dual with respect to $S^{ ^{\scriptstyle q}}.$  Equivalently, if the frame elements are $S^te_j,$ then the frame operator is $S^{(2t+1)/2}.$ 
\item  What is the dual frame for $\{S^{ ^{\scriptstyle p}}e_j\}_{j\in J}$?  We have
\bea
f&=&S^{ ^{\scriptstyle -p}}(S^{ ^{\scriptstyle p}}f)\\
&=&\underset{j\in J}{\sum}\langle S^{ ^{\scriptstyle p}}f\,,\,e_j\rangle S^{ ^{\scriptstyle -p}}(S^{-1}e_j)\nn\\
&=&\underset{j\in J}{\sum}\langle f\,,\,S^{ ^{\scriptstyle p}}e_j\rangle S^{ ^{\scriptstyle -p-1}}e_j
\end{eqnarray}
which shows that the dual of $S^{ ^{\scriptstyle p}}$ is $S^{ ^{\scriptstyle -(p+1)}}.$
\end{itemize}

\subsection{Riesz bases}
Let $J$ be some index set.
\begin{definition}
A set $\{e_j\}_{j\in J}$ is a {\it Riesz basis} for ${\mathcal H} \Leftrightarrow$
\be
\label{rb}
\exists\,A',B'\,\in\,\R_+\quad\mbox{ such that }\quad A'\|c\|^2_2\leq\bigl\|\underset{j\in J}{\sum}c_je_j\bigr\|^2\leq B'\|c\|^2_2
\end{equation}
for all finite sequences $c.$
\end{definition}
What is the essential difference between a Riesz basis and a frame?  A Riesz basis yields unique coefficients while a frame does not; alternatively a Riesz basis is a spanning set of minimal cardinality while a frame is redundant.  It is instructive to present an example distinguishing frames and Riesz bases.  Consider the finite Hilbert space $\R^2$ and the frame (but not Riesz basis) $G_1=\{(1,0),(0,1),(1,1)\}.$  $\,G_1$ is not a Riesz basis since we cannot find an $A'>0$ for which \eqref{rb} holds (example: $c_1=c_2=1,\,c_3=-1$ yields an inner sum of 0).  But $G_1$ is a frame since the inner sum in the inequality \eqref{frm} is positive $\forall f\in \R^2,$ and obviously finite.  Clearly any Riesz basis is a frame since we have uniqueness of coefficients and can take $c_j=\langle f,e_j\rangle.$

\section{Lattices and Gabor frames}
The main distinction between a frame and a Gabor frame is that a Gabor frame is necessarily defined on a symplectic lattice in $\R^{2d}.$  
\subsection{Lattices, symplectic lattices, and adjoint lattices}
We briefly review some basic facts about symplectic lattices and lattices in general, see \cite{CS93b} for details.  The {\em generator matrix} $L$ of a lattice $\La$ contains the coordinates of the column vectors that describe the fundamental domain of a lattice; i.e., $\La=L(\Z^d).$  The density $\xi(\La)$ of a lattice $\La$ is given by $\xi(\La)=1/\det(L).$  Write $\la=(x,\w)$ for $x,\w\in\R^d.$  The {\em symplectic form} is $[\la_1,\la_2]=x_2\cdot\w_1-x_1\cdot\w_2.$  A lattice $\La$ with generator matrix $L$ which preserves the symplectic form, i.e., $[L\la_1,L \la_2]=[\la_1,\la_2],$ is called symplectic, and is an element of the {\em symplectic group} {\em Sp}$(d).$  If $\La\in$ {\em Sp}$(d),$ then det$(L)=\pm 1,$ so {\em Sp}$(1)=${\em SL}$(2,\R).$
\begin{definition}
Any lattice $\La$ whose generator matrix $L=cA\,,\,A\in$ {\em Sp}$(d)\,,\,c\in\R,$ is called a {\em symplectic lattice}.
\end{definition}
\begin{definition}
The {\em adjoint lattice} of a symplectic lattice is $\La^\circ=\xi(\La)^{-\frac{1}{d}}\La,$ where $L$ is the generator matrix for the original lattice $\La.$
\end{definition}
A {\em separable} or {\em rectangular lattice} has generator of the form $L=\begin{bmatrix}\a&0\\0&\b\end{bmatrix}\,\a,\b\in\Z^d.$ 
In the sequel, all lattices under consideration will be symplectic, unless otherwise indicated.

\subsection{Gabor frames}\label{gabfr}
Let $T_x$ denote the translation operator whose action on $f$ is given by $T_xf(t)=f(t-x),$ and denote by $M_\w$ the modulation operator, the action of which on $f$ is given by $M_\w f(t)=f(t)e^{2\pi i\w\cdot t}.$  We then have the {\em commutation relations}
\be
T_xM_\w=e^{-2\pi i\w\cdot t}M_\w T_x,\label{commrels}
\end{equation}
and it is quite easy to demonstrate that
\be
(T_xf)^{\hat{ }}=M_{-x}\hat{f}\label{transfourier}
\end{equation}
and
\be
(M_\w f)^{\hat{ }}=T_\w\hat{f}\label{modfourier}
\end{equation}
\begin{definition}
Fix both nonzero $g,\in\Ltwo$ and a lattice $\La,$ with generator matrix $L$ and det$(L)>0.$ Write $ z(\la)=T_xM_\w\,,\,\la=(x,\w)\in\La).$  The set of time-frequency shifts
\be
(g,\La) = \{ z(\la)g\,|\,\la\in\La\}\label{Gaborsystem}
\end{equation}
is called a {\em Gabor system}, or {\em Gabor frame} if a frame.  
\end{definition}
If $\g=g$ we say that $(g,\La)$ is a Gabor frame.  If $\La$ is separable with $L=\begin{bmatrix}\a&0\\0&\b\end{bmatrix},\,\a,\b\in\Z^d,$ then we write
\be
(g,\a,\b) = \{g_{k\l}:=T_{\a k}M_{\b\l}g\,|\,k,\l\in\Z^d\}.\label{Gaborrectsystem}
\end{equation}
Note the similarity of an OFDM system to a Gabor frame; the structure is very similar.  A Gabor system includes a phase factor which occurs since modulation is performed first, then translation, while the order of operations is reversed for OFDM, so for OFDM there is no phase factor.  If the product $x\cdot\w\in\Z,$ the phase factor disappears.  Aside from the phase factor mentioned at the beginning of this Chapter, there is additional difference between Gabor frames and OFDM systems, namely that they are dual in a certain sense.  Analysis and synthesis by a Gabor frame is
\be
f\in \LtR \overset{C}{\rightarrow}\{d_{kl}\}\in \ell^2({\Bbb Z}^2)\overset{C^*}{\rightarrow}\tilde{f} \in \LtR
\end{equation}
where the frame must be complete (i.e., the closed linear span of the frame is the entire space) and may be linearly independent for perfect reconstruction of $f.$  OFDM eoncoding, transmission, and data recapture is
\be
\{c_{k\ell}\}\in\ell^2({\Bbb Z}^2)\overset{C^*}{\rightarrow}s\in\LtR\overset{C}{\rightarrow}\{\tilde{c}_{kl}\}\in\ell^2({\Bbb Z}^2)
\end{equation}
where the $\psi_{kl},$ of which $s$ is a linear combination, must be linearly independent and may be complete for perfect reconstruction of $\{c_{kl}\}.$\\

Proposition 9.4.4 in \cite{Gro01} states that a Gabor system on a symplectic lattice is the image of a Gabor system on a separable lattice under a unitary operator.\\

\begin{definition}
The {\em Gram matrix} $R_{g,\La}$ associated with $(g,\La)$ is
defined via
\begin{equation}
\{R_{g,\La}\}_{\la,\la'} := \langle g_{\la'},g_{\la}\rangle.
\label{gram1}
\end{equation}
\end{definition}
If $\La$ is rectangular, we write $R_{g,T,F}.$  $R_{g,\La}$ may be invertible, as is the case when $\{g_\la\}$ is a Riesz basis, or singular, as when $\{g_\la\}$ is a redundant (possibly tight) frame.\\

The action of the Gabor frame operator $S$ is:
\be
Sf =\sum_{\la\in\La}\left\langle f, z(\la)g\right\rangle z(\la)g\label{gabfropdef}
\end{equation}
or
\be
Sf =\sum\limits_{k,\l \in \Z^d} \langle f,g_{k\l} \rangle g_{k\l}\label{gabrectfropdef}
\end{equation}
if $\La$ is separable.  If we wish to emphasize the dependence of $S$ on the function $g$ and the lattice $\La,$ we write $S_{g,g,\La}.$

\begin{definition}
Let $g,\g\in\Ltwo,$ and suppose that every $f\in\Ltwo$ can be expanded as
\bea
f&=&\sum_{\la\in\La}\langle f, z(\la)g\rangle z(\la)\g\\
&=&\sum_{\la\in\La}\langle f, z(\la)\g\rangle z(\la)g.
\end{eqnarray}
Then $(g,\g,\La)$ is a Gabor frame, and in this case we call $\g$ (resp. $g$) the {\em dual Gabor frame} or simply the {\em dual frame} to $g$ (resp. $\g$).
\end{definition}
It is shown in \cite{Gro01} that if $S$ is a Gabor frame operator, both $S$ and the inverse frame operator $S^{-1}$ commute with time-frequency shifts which generate $S,$ i.e.,  
\be
ST_xM_\w=T_xM_\w S\qquad\mbox{and}\qquad S^{-1}T_xM_\w=T_xM_\w S^{-1} \qquad\forall\,\,(x,\w)\in\La.
\end{equation}
Thus we have
\bea
f=SS^{-1}f&=&\sum_{\la\in\La}\langle S^{-1}f, z(\la)g\rangle z(\la)g\nn\\
&=&\sum_{\la\in\La}\langle f, z(\la)S^{-1}g\rangle z(\la)g\label{Scommute}
\end{eqnarray}
where \eqref{Scommute} hold since $S,$ and hence $S^{-1},$ is self-adjoint.  Thus $S^{-1}g$ is a dual window for $g,$ called the {\em canonical} dual window.  We conclude that the dual frame $\G$ of a Gabor frame $G$ has the same Gabor structure as $G.$  The same property holds for $S^\nu,\,\nu\neq 0,$ and in particular for $S^{-\frac{1}{2}}.$  If we denote by $S^\nu(g,\La)$ the action of the $\nu^{th}$-power of the Gabor frame operator $S$ on the entire Gabor frame $(g,\La),$ this property is denoted
\be
S^\nu(g,\La)=(S^\nu g,\La).
\end{equation}
We demonstrate this in Appendix \ref{spowercommute}.  

\section{L\"owdin orthogonalization in the continuous setting}\label{contLowd}
We describe explicitly the role that L\"owdin orthogonalization plays in our OFDM scheme.  Fix a function $f$ with good decay, let $L=\begin{bmatrix}T&0\\0&F\end{bmatrix}$ be the generator matrix for $(f,T,F),$ and assume that $(f,T,F)$ is linearly independent.  Since $(f,T,F)$ is a countable set, we can index its elements by $\Z;$ it is however more transparent to index the elements of $(f,T,F)$ by $\Z\times\Z.$  Define $\{f_{j,k}\}_{j,k\in\Z}:=(f,T,F).$  Since $\{f_{j,k}\}_{j,k\in\Z}$ is linearly independent, it forms a Riesz basis for its closed linear span, and we conclude that there exist constants $A', B'\in\R_+$ such that
\be
A'\|c\|_2\leq\bigl\|\sum_{j,k\in\Z}c_{j,k}f_{j,k}\bigr\|\leq B'\|c\|_2\label{rieszbasiscondition}
\end{equation}
holds for all {\em finite} sequences $\{c_{j,k}\}.$  Let $J$ be some finite index set.  Then for any finite sequence $\{c_{j,k}\},$ we have 
\be
\langle Rc,c\rangle=\sum_{j,k\in J}\langle f_{j',k'},f_{j,k}\rangle c_{j',k'}\overline{c}_{j,k}=\bigl\|\sum_{j',k'\in J}c_{j',k'}f_{j',k'}\bigr\|^2.\label{gramriesz}
\end{equation}
Since any $c_0\in\ltZ$ can be approximated arbitrarily closely by some finite sequence $c,$ the action of $R$ is bounded above and below on the unit ball in $\ltZ$ by \eqref{gramriesz}.  Also $\langle Rc\,,\,c\rangle=c^*Rc>0,$ and we conclude that \eqref{rieszbasiscondition} implies that the Gram matrix $R$ is a positive invertible operator.  From this we can infer the existence of $R^{-\frac{1}{2}}.$  We can construct an orthonormal
system $(\phi,T,F)$ from $(f,T,F)$ via
\begin{equation}
\label{o1}
\phi = \sum_{k,\l} R^{-\frac{1}{2}}_{k,\l,0,0} f_{k,\l},
\end{equation}
where
\begin{equation}
\label{o2}
R_{k,\l,k',l'}=\langle f_{k',l'},f_{k,\l} \rangle.
\end{equation}
We will write
\begin{equation}
\label{o3}
\phi = L\ddot{o}\,(f,T,F)
\end{equation}
for the function computed by this procedure.  The proof that $\phi$ generates an orthonormal system is the subject of the next section.

\section{L\"owdin orthogonalization via Gabor analysis}\label{invsqrtgram}
We derive here the general symplectic case of L\"owdin orthogonalization.

\subsection{Relevant theorems and lemmas of Gabor analysis}\label{relevthmgabor}
In this section $\La$ will always be a symplectic lattice in $\R^d.$
\begin{definition}
Let $f,\g\in\Ltwo,$ and suppose that $(f,\g,\La)$ is a Gabor frame.  The action of frame operator on $h\in\Ltwo$ is  
\be
S_{f,\g,\La}(h)=\xi(\La)\underset{\la\in\La}{\sum}\langle h, z(\la)f\rangle z(\la)\g.
\end{equation}
\end{definition}
\begin{theorem}\label{jansrep}
{\em (Janssen's representation) \cite{Jan95}}
\be
S_{f,\g,\La}(h)=\xi(\La)\underset{\la_\circ\in\La^\circ}{\sum}\left\langle\g, z(\la_\circ)f\right\rangle z(\la_\circ)h.
\end{equation}
\end{theorem}
We will also cite the following important theorems:
\begin{theorem}\label{WRbr}
{\em (Wexler-Raz biorthogonality relations) \cite{WR90}} Suppose that $D_g$ and $D_\g$ are bounded on $\ll.$  Then the following conditions are equivalent:\\
\indent (a) $S_{f,\g,\La}=S_{\g,f,\La}=I\quad\mbox{on}\,\,\Ltwo$\\
\indent (b) $\xi(\La)^d\langle \g, z(\la_\circ)f\rangle =\d_{\la_\circ 0}\quad\mbox{for}\,\la_\circ\in\La^\circ.$
\end{theorem}
\begin{corollary}\label{WRtgtfrm}
A Gabor system $(f,\La)$ is a tight frame if and only if $(f,\La^\circ)$ is an orthogonal system.  In this case the frame bound $A$ satisifes $A=\xi(\La)^d.$

\begin{myproof}(From \cite{Gro01}.)  We have
\be
\frac{1}{A}S_{f,f}=I
\end{equation}
where $I$ is the identity on $\Ltwo.$   Using $\g=\frac{1}{A}f$ in Theorem \ref{WRbr} yields
\be
\frac{\xi(\La)^d}{A}\left\langle z(\la_\circ)f, z(\la'_\circ)f\right\rangle=\d_{\la_\circ,\la'_\circ}.
\end{equation}
For $\la_\circ=\la'_\circ$ we get $A=\xi(\La)^d\|f\|^2_2$ for the frame bound.  Conversely, if $(f,\La^\circ)$ is an orthogonal system, Theorem \ref{jansrep} yields $S_{f,f}=\xi(\La)^d\|f\|_2^2\,I.$  Thus $(f,\La)$ is a tight frame.
\end{myproof}
\end{corollary}

\begin{theorem}\label{RSdp}
{\em (Ron-Shen duality principle) \cite{RS97a}} Let $f\in\Ltwo$ and $\a,\b>0.$  Then the Gabor system $(f,\La)$ is a frame for $\Ltwo$ if and only if $(f,\La^\circ)$ is a Riesz basis for its closed linear span.
\end{theorem}
It is interesting to note that this duality relation can already be extracted from a paper of Rieffel on non-commutative tori \cite{Rie88}.

\begin{lemma}\label{frmtgtfrmorth}
Let $\Hi_0$ denote the closed linear span of $(f,\La^\circ).$ If $(f,\La^\circ)$ is a Riesz basis for $\Hi_0,$ then $(S_{f,f,\La}^{-\frac{1}{2}}f,\La^\circ)$ is an orthonormal basis for $\Hi_0.$

\begin{myproof}\bea
&&(f,\La^\circ) \text{ is a Riesz basis for }\Hi_0\nn\\
&\Rightarrow& (f,\La) \text{ is a frame for } \Ltwo\label{tgf2}\\
&\Rightarrow& (S_{f,f,\La}^{-\frac{1}{2}}f,\La)\text{ is a tight frame for }\Ltwo\label{tgf3}\\
&\Rightarrow& (S_{f,f,\La}^{-\frac{1}{2}}f,\La^\circ)\text{ is an orthonormal basis for }\Hi_0\label{tgf4}
\end{eqnarray}
where we used Theorem \ref{RSdp} for \eqref{tgf2}, \eqref{tgf3} derives from Lemma \ref{tgtfrmS12}, and \eqref{tgf4} is due to Corollary \ref{WRtgtfrm}.
\end{myproof}
\end{lemma}

\begin{lemma}\label{density}
If $(f,\La)$ is a frame for $\Ltwo,$ then $\xi(\La)\geq 1.$
\end{lemma}
This lemma is proved in Section 9.4 of \cite{Gro01}.  One of the consequences of Lemma \ref{density} and Lemma \ref{frmtgtfrmorth} is 

\begin{corollary}
If $(f,T,F)$ is an orthonormal system in $\Ltwo,$ then $TF\geq 1.$ \label{density2}
\end{corollary}

\subsection{L\"owdin orthogonalization and the tight frame operator}\label{LoRb}

\begin{lemma}\label{grail}
Let $\La$ be a lattice in $\R^{2d}$ with $\xi(\La)>1.$  Fix $f\in\Ltwo,$ and denote by $f_\la$ members of the Gabor system ($f,\La),$ and by $f_{\la_\circ}$ members of the Gabor system ($f,\La^\circ).$  Suppose $(f,\La^\circ)$ is a Riesz basis for its closed linear span, denoted $\Hi_0.$  Denote by $({R_{f,\La^\circ}})_{\la_\circ \la_\circ'}=\langle f_{\la_\circ'},f_{\la_\circ}\rangle$ the Gram matrix of $(f,\La^\circ).$  Then
\be
\phi:=\sum_{\la_\circ\in\La^\circ} ({R^{-\frac{1}{2}}_{f,\La^\circ}})_{_{\scriptstyle{\la_\circ, 0}}}\,\, f_{\la_\circ}\,=S_{f,f,\La}^{-\frac{1}{2}}f.\label{aileq}
\end{equation}
Furthermore, either side of (6.19) generates an orthonormal system on $\La^\circ.$ 
\end{lemma}

\begin{myproof}
The right-hand side is in fact an orthonormal system on $\La^\circ,$ by Lemma \ref{frmtgtfrmorth}.  It remains to prove that the sides of (6.19) are equal.  For convenience we will denote $S_{\La,f,f}$ by $S_1.$  We have, (cf. \eqref{alsodual}),
\bea
S_1^{-\frac{1}{2}} z(\la)f&=&\underset{\la'\in\La}{\sum}\left\langle S_1^{-\frac{1}{2}} z(\la)f\,,\,S_1^{-\frac{1}{2}}z(\la')f\right\rangle S_1^{-\frac{1}{2}}z(\la')f\\
&=&\underset{\la'\in\La}{\sum}\left\langle z(\la)f\,,\,S_1^{-\frac{3}{4}}z(\la')f\right\rangle S_1^{-\frac{3}{4}}z(\la')f\label{tgtfrmop}
\end{eqnarray}
since $\forall h\in{\mathcal H},$ and for any frame operator $S,\,\, S^{-\frac{1}{2}}h=S^{\a} S^{-\frac{1}{2}}S^\a h$ and powers of $S$ are self-adjoint.  $S_1^{-\frac{3}{4}}f$ generates a frame by \eqref{alsodual}.  Janssen's representation of the action of the frame operator $S^{-\frac{1}{2}}_1$ on arbitrary $h$ is
\be
\label{jansreptgt}
S_1^{-\frac{1}{2}}(h)=\xi(\La)^d\underset{\la_\circ\in\La^\circ}{\sum}\left\langle S_1^{-\frac{3}{4}}f, z(\la_\circ)S_1^{-\frac{3}{4}}f\right\rangle z(\la_\circ)h.
\end{equation}
Write $R_1:=R_{f,\La^\circ},\,\upsilon:=\xi(\La)^{ ^{\scriptstyle{-\frac{d}{4}}}}S_1^{-\frac{3}{4}}f.$  If we can show that
\be
R_1^{-\frac{1}{2}}=R_2:=\xi(\La)^d R_{\upsilon,\La^\circ}=\xi(\La)^d \{\left\langle z(\la_\circ)\upsilon, z(\la'_\circ)\upsilon\right\rangle\}
\end{equation}
our task will be complete.  The $(\la_\circ,\la_\circ''')$ entry of $\xi(\La)^{2d}R_1R_2^2$ is
\bea
&&\xi(\La)^{2d}\underset{\la_\circ'\in\La^\circ}{\sum}\bigl\langle z(\la_\circ)f, z(\la'_\circ)f\bigr\rangle\underset{\la_\circ''\in\La}{\sum}\bigl\langle z(\la'_\circ)\upsilon, z(\la''_\circ)\upsilon\bigr\rangle\bigl\langle  z(\la''_\circ)\upsilon, z(\la'''_\circ)\upsilon\bigr\rangle\label{r1r22}\\
&=&\xi(\La)^{\frac{3d}{2}}\underset{\la_\circ'\in\La^\circ}{\sum}\bigl\langle  z(\la_\circ)f, z(\la'_\circ)f\bigr\rangle\bigl\langle z(\la'_\circ)S_1^{-\frac{1}{2}}\upsilon, z(\la'''_\circ)\upsilon\bigr\rangle\label{altinnprodtgt}\\
&=&\xi(\La)^{\frac{5d}{4}}\bigl\langle z(\la_\circ)S_1(S_1^{-\frac{5}{4}}f), z(\la'''_\circ)\upsilon\bigr\rangle\\
&=&\xi(\La)^d\bigl\langle z(\la_\circ)S_1^{-\frac{1}{4}}f, z(\la'''_\circ)S_1^{-\frac{3}{4}}f\bigr\rangle\label{altinnprodfrm}\\
&=&\delta_{\la_\circ\la_\circ'''}\label{oddduals}
\end{eqnarray}
where \eqref{altinnprodtgt} derives from \eqref{tgtfrmop} and Proposition 2.6 in \cite{Jan95}, \eqref{altinnprodfrm} is due again to Proposition 2.6 in \cite{Jan95}, and \eqref{oddduals} holds by \ref{WRbr} on the duals in \eqref{alsodual}.  Thus we have proven that
\be
I = \xi(\La)^{2d}R_1R_2^2,\quad\Rightarrow\quad R_1^{-\frac{1}{2}}=\xi(\La)^d R_2,\nn
\end{equation}
i.e., Lemma \ref{grail} holds by Theorem \ref{jansreptgt} with $h=f.$ 
\end{myproof}

\vspace{.15in}

In the sequel we shall call OFDM using the above algorithm {\em pulse-shaping OFDM}, or PS-OFDM.

\section{Optimality of $\phi=L\ddot{o}\,(f,\La)$}\label{s:pulse}
The following theorem shows that the orthogonalization procedure described
above yields pulse shapes that are optimal in a certain sense.
\begin{theorem}
\label{th:min}
Let $\{f_{\la}\,|\,\la\in\La\}$ with $f_{\la}(t)=[z(\la)f](t)$ and $\|f\|=1$ be a
Riesz basis for the separable Hilbert space ${\cal H}$.\\
(i) The function $\phi=L\ddot{o}\,(f,\La)$ is optimally
close to $f$ in the sense that $\phi$ solves the optimization problem
\begin{gather}
\arg \underset{\phi}{\min} \,\|Af - Ah\|_2 \qquad
\text{when $(h,\La)$ is an orthonormal basis for ${\cal H}$.}
\label{opt}
\end{gather}
(ii) Furthermore
\begin{equation}
\|Af - A\phi \| \le \max\,\{|1 - \theta_{\min}^{-\frac{1}{2}}|^2,
|1 - \theta_{\max}^{-\frac{1}{2}}|^2\},
\label{spect}
\end{equation}
where $\theta_{\min}$ and $\theta_{\max}$ are the minimal and
maximal eigenvalues respectively of $R$.
\end{theorem}

\begin{myproof}
By the orthogonality relations of the ambiguity function (e.g., see Chapter 4.2 in \cite{Gro01}) we have that $\|Af - Ah\| = \|f - h\|^2$. Now \eqref{opt} follows immediately from Theorem 2.1 in \cite{JS00} by using the duality relations between the systems $(f,\La)$ and $(f,\La_\circ)$.\\

To show part (ii) we again use the duality relations.  Recall from Section \ref{LoRb} that we can write $\phi = \xi(\La)^{-\frac{1}{2}}\SQI f$ where $S$ is the Gabor frame operator of the frame $(f,\La_\circ)$ (the normalization by $\xi(\La)^{-\frac{1}{2}}$ yields $\|\phi\|=1$). There holds $\theta_{\min}\le \xi(\La)S \le \theta_{\max}$, see Theorem 3.1 in \cite{Jan95}. The estimate \eqref{spect} follows from the calculation
\bea
\|f - \psi\|&=&\|(I-\xi(\La)^{-\frac{1}{2}}\SQI)f\|\le \|I-\xi(\La)^{-\frac{1}{2}}\SQI \|_{op}\nn \\
&\le& \max\big\{|1-1/\sqrt{\theta_{\min}}|, |1-1/\sqrt{\theta_{\max}}|\big\}.
\label{map}
\end{eqnarray}
\end{myproof}

\section{Invariance of $\ph=L\ddot{o}\,(g,T,F)$ under the Fourier transform}
In this section we shall treat rectangular lattices only.  Recall that $T$ is the symbol length and $F$ is the carrier separation for an OFDM system. \\

From here on we will always denote Gaussians by $\ggam$, i.e.,
\begin{equation}
\label{gauss}
\ggam(t):=(2\sigma)^{1/4} e^{-\pi \sigma t^2}=\D_{\sqrt{\s}}(2^{1/4}e^{-\pi t^2}).
\end{equation}
Since the pulse shape obtained from the Gaussian via \eqref{o3} plays
a major role throughout the paper we introduce a separate notation for
this pulse to distinguish it from other pulses one might use for OFDM. We denote
\begin{equation}
\label{o4}
\psigam := L\ddot{o}\,(\ggam,T,F),
\end{equation}
and if $\sigma=1$ we skip the subscripts and simply write $\ph$ or $g,$ respectively.

\begin{figure}[!ht]
\begin{center}
\epsfig{file=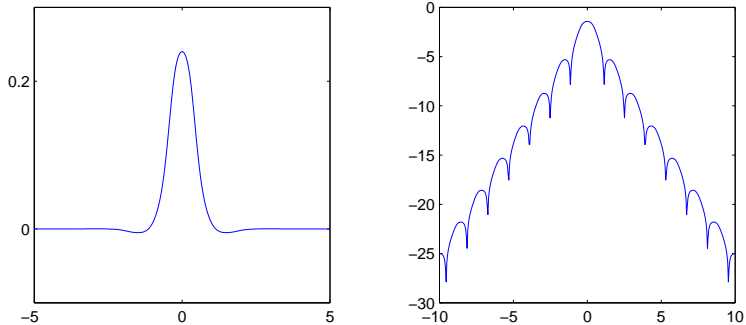,width=100mm}
\caption{OFDM window produced by L\"owdin orthogonalization (left) and its log-plot (right).  Note the similarity to the Gaussian, except for the zero-crossings, and the exponential envelope decay evident in the log-plot.}
\label{fig:fig8}
\end{center}
\end{figure}

In Section \ref{s:pulse} we made no mention of the TFL properties $f$ must have for $\phi$ to be a good candidate for an OFDM system.  Theorem \ref{th:min} says that the better the TFL of our Gabor atom $f,$ the better the TFL of the L\"owdin orthogonalized pulse $\phi.$  It is thus natural to begin with a normalized Gaussian pulse\footnote{We shall show in Chapter \ref{GenPulseDesign} that in certain situations, $g$ should be chirped and/or dilated.} $g.$  In addition to possessing the best possible TFL, $g$ is invariant under the Fourier transform $\F$ as well.   This property is inherited by $\ph:=L\ddot{o}\,(g,T,F)$ provided that $T=F$ and $TF\in\Z.$  Indeed, for all $f\in (L^1\cup L^2)(\R),$we have $(M_yT_xf)^{\hat{ }} =T_yM_{-x}\hat{f}=e^{2\pi ixy}M_{-x}T_y\hat{f}.$  Thus we see that under the Fourier transform, the roles of $T$ (symbol length) and $F$ (carrier separation) are rotated by $\frac{\pi}{2},$ and a phase factor depending on the product $xy$ is introduced.  So, if $T=F$ and $xy=kT\ell F\in\Z,$ the entire Gabor basis $(g,T,F)$ is invariant under $\F$ (since $g$ is invariant under $\F$).  Or, if we have only that $xy=kT\ell F\in\Z,$ then still $\ph$ is invariant under ${\mathcal F}.$  In this case,
\begin{eqnarray}
{\mathcal F}\ph&=&{\mathcal F}\sum_{k,\ell} R^{-\frac{1}{2}}_{k,\ell,0,0} g_{k,\ell}\\
&=&\sum_{k,\ell} R^{-\frac{1}{2}}_{k,\ell,0,0} {\mathcal F} g_{k,\ell} \qquad \text{(Fubini's Theorem)} \\
&=&\sum_{k,\ell} R^{-\frac{1}{2}}_{\ell,k,0,0} g_{\ell,k}\qquad\quad\text{(R is Hermitian)}\\
&=&\ph.\nn
\end{eqnarray}

\section{Condition numbers and TFL}\label{condnTFL}
Theorems \ref{th:amb} and \ref{th:min} imply that the OFDM pulse shapes
$\psigam=L\ddot{o}\,(\ggam,T,F)$ are optimally time-frequency localized
among all possible OFDM pulse shapes in the sense that they are closest 
(in the $L^2$-sense) to the Gaussian $\ggam$.
Unlike $\ggam$ we see from Example~2 in Section 6 of \cite{BJ00} that
$\psigam$ has ``only'' exponential decay in time and frequency. Here
the constant in the exponent depends on the condition number of 
$R_{\ggam,T,F}$. In Figure \ref{fig:condR} we display the behavior of 
$\condR$ in dependence of the spectral efficiency 
$(TF)^{-1}$ for $(TF)^{-1} = \frac{p}{p+1}$ with $p=1,2,\dots,25$. 
When $(TF)^{-1} \rightarrow 1$ we have $\condR \rightarrow \infty$
as a consequence of the Balian-Low Theorem and the functions $\psigam$
become increasingly ill-localized, cf. also Theorem 7.3 in \cite{FZ98}.
Fortunately for typical values of $(TF)^{-1}$ in OFDM such as
$TF \in [\frac{1}{2},\frac{4}{5}]$ (and $\sigma=F/T$) the condition number
of $R$ is small. For instance for $TF=2$ one can show 
analytically \cite{Jan98b} that $\cond(R_{\ggam,T/\sqrt{\sigma},2\sqrt{\sigma}/T})=\sqrt{2}$ and for 
$TF=\frac{3}{2}$ a numerical inspection yields 
$\cond(R_{\ggam,T/\sqrt{\sigma},\frac{3}{2}\sqrt{\sigma}/T})=\sqrt{3}$. 
The second graph in~Figure \ref{fig:condR} corresponds to the
condition number of the Gram matrix of a modified OFDM system,
which we will introduce in Chapter \ref{GenPulseDesign}.\\

For a general theoretical investigation of the decay properties of 
pulse shapes generated via \eqref{o3} refer to \cite{Str01}. 
Efficient numerical algorithms to compute $\psigam$ that are also
applicable if $TF \neq \N$ can be easily derived from the algorithms
in \cite{Str97a}.\\

It has been shown in \cite{BJ00} that for certain pulse shapes the so-called 
IOTA approach in \cite{FAB95} is equivalent to \eqref{o3}. In \cite{FAB95} 
the authors computed numerically the TFL of $\ph$ via \eqref{up0} as 
$\Delta \tau_{\ph} \cdot \Delta \nu_{\ph}=1.024\cdot\frac{1}{4\pi}$ 
(as compared to $\frac{1}{4\pi}$ for the Gaussian) and conjectured that this 
is the optimal value for all pulse shapes that create an OFDM system with 
spectral efficiency $\frac{1}{2}$. Theorem \ref{th:min} provides additional
evidence for the correctness of their conjecture (as long as we consider 
only rectangular time-frequency lattices).

\section{Biorthogonal frequency-division multiplexing and OFDM-BFDM interpolates} \label{ss:bfdm}

Recently biorthogonal Frequency-Division multiplexing (BFDM) has been introduced \cite{KM98}. The reason for considering BFDM is that by giving up the orthogonality condition one can improve the TFL of the transmitter pulse shape for $TF>1$. This comes however at the cost of increased sensitivity to AWGN, since we no longer have an orthogonal set of transmitter functions. Nevertheless the concept of using general time-frequency lattices in order to optimally adapt the transmission functions for time-frequency dispersive channels can be applied without any modification to BFDM. \\

BFDM and OFDM can be seen as special cases of the following general framework which allows one to fine-tune the balance between AWGN sensitivity and ISI/ICI robustness. Assume we are given a non-orthogonal linearly independent set of transmission functions $(h, \La)$ (here $h$ may be the Gaussian function) and let $R_{h,\La}$ denote the Gram matrix as defined in \eqref{gram1}.\\

Let $R=R_{g,\La^\circ}.$  A pair of perfect reconstruction transmitter and receiver function sets $(\phi,\La)$ and $(f,\La)$ can be constructed by defining (\cite{Jan95})
\begin{equation}
\phi = \sum_{k,\l} R^{-p}_{k,\l,0,0} h_{k,\l}, \quad
f = \sum_{k,\l} R^{-1+p}_{k,\l,0,0} h_{k,\l}, \quad p\in [0,1].
\label{GFDM}
\end{equation}
Clearly the choice $p = \frac{1}{2}$ leads to a standard OFDM system (with $\phi=f$) with the pulse shape construction as used throughout the paper. The choice $p=0$ or $p=1$ corresponds to BFDM.  For $0<p<1/2$ (or $1/2<p<1$ respectively) we obtain a transmission system that ``interpolates'' between OFDM and the standard setup of BFDM.  We term these such systems {\em extended-biorthogonal frequency-division multiplexing}, or EBFDM, systems.  Since $\phi$ is the pulse whose modulations and translations are subject to the vagaries of the wireless channel, if $\phi$ has good TFL the signal is likely to be more robust against multipath and Doppler dispersion.  The EBFDM systems $(\phi,\La)$ and $(f,\La)$ become ``more and more orthogonal'' as $p\rightarrow\frac{1}{2}$ (the condition number of $R_{\phi,\La}$ and $R_{f,\La}$ approaches 1), but the TFL of $\phi$ might degrade.  Hence by varying $\alpha$ we should be able to emphasize either the AWGN performance or the ISI/ICI performance.

\section{Numerical simulations}
For the computer simulations we assume the channel to be wide-sense stationary with uncorrelated scattering (WSSUS), cf. \cite{Bel63}.  We consider a Jakes-type Doppler behavior and exponential decay in the temporal direction for the scattering function \cite{Rap96}. In addition we assume that $\tau_0 \nu_0 < 1$, i.e., the channel is underspread \cite{Bel63}.  Following other work on pulse shaping OFDM \cite{VH96,HB97,KM98,BDH00,FAB95}, we consider transmission systems with spectral efficiency $\rho = 0.5$ (in practice a higher spectral efficiency is certainly preferable).
\\

Recall from \eqref{innprod} that the received data $\tilde{c}_{k,\l}$ can be
written as 
\be
\tilde{c}_{k,\l}=\sum_{k',l'} c_{k',l'}\langle \Hchan\phi_{k',l'},\phikl \rangle + \langle \eps, \phikl \rangle.
\end{equation}

Let $S_{\Hchan}$ denote the time-varying transfer function (also known as Weyl symbol) of $\Hchan$, see \cite{Bel63}. For underspread channels and pulse shapes that are well concentrated in time and frequency it has been shown in \cite{KM98} for rectangular time-frequency lattices that 
\be
\langle \Hchan \phi_{k,\l}, \phi_{k,\l} \rangle =
S_{\Hchan}(\la_k,\nu_l), \quad (\la_k,\nu_l)\in\La
\label{}
\end{equation}
and that $\langle \Hchan \phi_{k',l'}, \phi_{k,\l} \rangle \approx 0$
for $k \neq k', l \neq l'$. This result also holds for general lattices.  We thus can use a
one-tap equalizer (simple division by the value of $S_{\Hchan}$ at one point in the TF-plane) to approximately recover the data.  The received data can be equalized by computing
\begin{equation}
d_{kl} = \frac{\tilde{c}_{kl}}{S_{\Hchan}(\la_k,\nu_l)+\sigma^2},
\label{equal}
\end{equation}
where $\sigma$ is the variance of the AWGN.\\

Furthermore, for a fair evaluation of the robustness of CP-OFDM, PS-OFDM (and later LOFDM, defined in the next chapter) against ISI/ICI, we must exclude any additional effects resulting from coding or channel estimation. 
Thus we assume that the channel is known at the receiver and do not perform any coding/decoding of the data.\\

We measure performance by the received signal's {\em signal-to-noise-plus-interference ratio}, defined by 
\be
\text{SNIR}(\Hchan,s)=-10\log_{10}\left(\frac{\|c-d\|_2}{\|c\|_2}\right),\label{snirdef}
\end{equation}
where $c$ is the vector of transmitted symbols, and $d$ is the vector of received, equalized (by use of \eqref{equal}) symbols.

\subsection{Our numerical implementation of OFDM and PS-OFDM}\label{numimplofdm}
We discuss briefly the numerical implementation of an OFDM system based on an orthogonalized Gaussian.  First, one should have the discretized Gaussian invariant under the DFT, in order that the L\"owdin-orthogonalized pulse $\ph$ be invariant under the DFT (recall Section \ref{contLowd} for the continuous case).  Thus can we consistently reflect the continuous signal built from $\ph$ in both the time and frequency domains, and avoid any issues stemming from different sampling intervals in the time and frequency domains.  Note that if the sampling interval of the underlying continuous function is $\Delta t,$ then DFT$\{f(n\Delta t)\}=\text{DFT}\{f\}(k/(N\Delta t)),$ where $N$ is the DFT length and $n,k$ are time and frequency indices, respectively. Thus, to enforce DFT-invariance of discretized $\ph,$ we must have equal spacing in the time and fequency domains, and so we demand $\Delta t=(N\Delta t)^{-1},$ which yields $\Delta t=1/\sqrt{N}.$\\

We choose a pulse length of $N=512$ yielding $\Delta t=1/(16\sqrt{2}).$  Let $(g,T,F)$ with $T=F$ and $g=\sqrt[4]{2}e^{-\pi t^2}$ be a Gabor system with oversampling by a factor of two; i.e., $T=F=1/\sqrt{2},$ which yields a system with spectral efficiency $\r=\frac{1}{2}$ (of course in practice a higher spectral efficiency is desirable, but $\r=\frac{1}{2}$ was chosen for ease of programming).  Then by the argument in Section \ref{invsqrtgram}, our OFDM system is $(\ph,\sqrt{2},\sqrt{2}),$  with $\ph=\sum_{k,\l}R^{-1/2}_{k,\l,0,0}g_{k,\l}$ with $g_{k,\l}=g(t-kT ')e^{2\pi i\l F't}.$  We then have $T'=T_0\Delta t,$ where $T_0=32=F_0$ are the values we use for $T'$ and $F'$ while programming, since the underlying continuous time and frequency are irrelevant on the computer.  Let $N_c$ be the number of subcarriers, and $N_k$ the number of time-slots.  Let $c$ be the vector of data symbols, and upsample $c$ by $F$ to yield $c_F.$  Taking the inverse DFT of $c_F$ yields (we shall not complicate things by adopting or discussing the indexing or normalization conventions of MATLAB's FFT function here)
\be
\hat{c}_F[j]=\sum_{\l=1}^{N_c}c_{F_\l} e^{2\pi ij\l F_0/N}=\sum_{\l=1}^{N_c}c_\l e^{2\pi i\l Ft}
\end{equation}
which is the sum of the symbols times the modulators since $F=\sqrt{2}\,,\,t=j\,dt,$ and $dt=1/\sqrt{N}.$  Multiplication by the time-shifted window completes the encoding operation.  The signal $s$ is then subjected to a dispersion subroutine which imparts to $s$ the effects of multipath delay and Doppler dispersion.  Next we add noise to $s,$ demodulate via the DFT, equalize using \eqref{equal}, and compute the SNIR.  

\section{Comparison of performance of pulse-shaping OFDM with cyclic-prefix OFDM}\label{psofdmcpofdm}
It is well-known that inserting a cyclic prefix (CP) into the OFDM symbol at each time-slot, while reducing the spectral efficiency (e.g., by 20 percent in IEEE 802.11a), effectively precludes ISI, as long as the CP length is longer than the effective temporal support of the channel scattering function.  CP-OFDM is of no use in combatting ICI resulting from frequency dispersion due to the Doppler effect from moving transmitters and/or receivers.  In this section we present the results of numerical simulations to this effect.\\

We have the freedom of several parameters in our simulations.  A noise level of $10^{-n}$ shall mean that AWGN with a variance of $10^{-n}$ was added to the signal before it entered the receiver.  A noise vector ${\mathcal N}_{sig}$ of unit variance was generated in MATLAB, and multiples of that same ${\mathcal N}_{sig}$ were used in every experiment in this section and Section \ref{pbofdmcpofdm}.  The simulated channel scattering function is governed by four parameters, $\t_0,\nu_0,\a,$ and $\b.$  The value of $\t_0$ represents the maximum multipath delay (now in units of $\Delta t$), and $\nu_0$ represents the maximum Doppler spread (now in units of $\Delta f=\Delta t$).  For $\t_0=\nu_0<22$ the channel is underspread (the IEEE 802.11a standard cyclic prefix length is 20$\%,$ this is $\t_0=\nu_0<7$ in our setting).  The parameter $\a<1$ is the maximum (attenuated) amplitude, relative to the line-of-sight strength, of the reflected part of the signal arriving at a delay of $\t_0,$ while $\b$ determines the proximity of the maximal frequency shifts determined by $\nu_0$ to the asymptotes of the Jakes PSD (see  \cite{Rap96}); alternatively, $\b$ determines how tightly clustered is the incoming signal around zero Doppler.  The channel {\em spreading function} is the scattering function multiplied by a normally distributed noise matrix ${\mathcal N}_{scat}$  with variance 1.  Again, that same ${\mathcal N}_{scat}$ was used in every experiment in this section and Section \ref{pbofdmcpofdm}.

\begin{figure}[!ht]
\begin{center}
\subfigure[]{
\epsfig{file=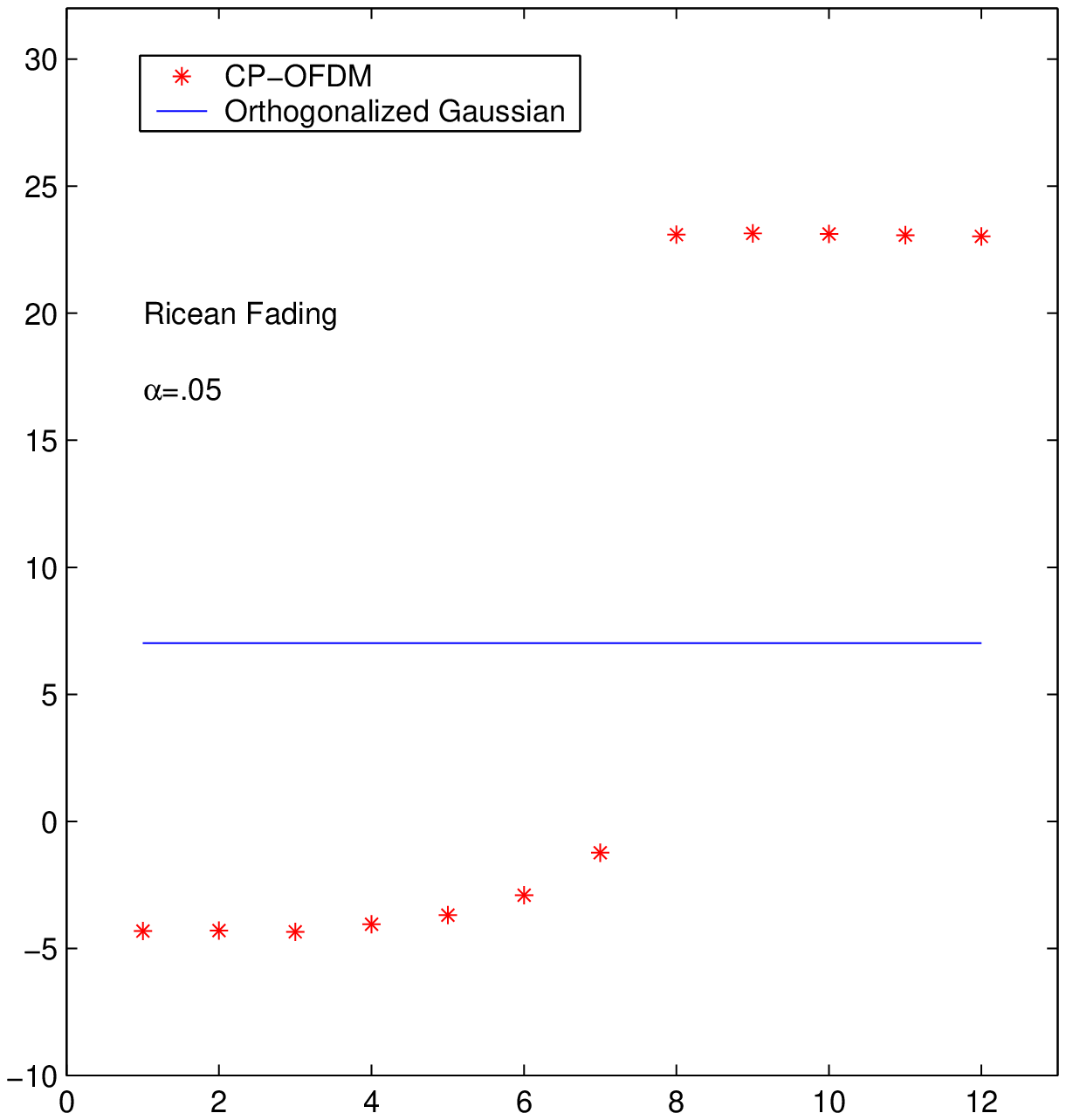,width=65mm,height=45mm}}
\subfigure[]{
\epsfig{file=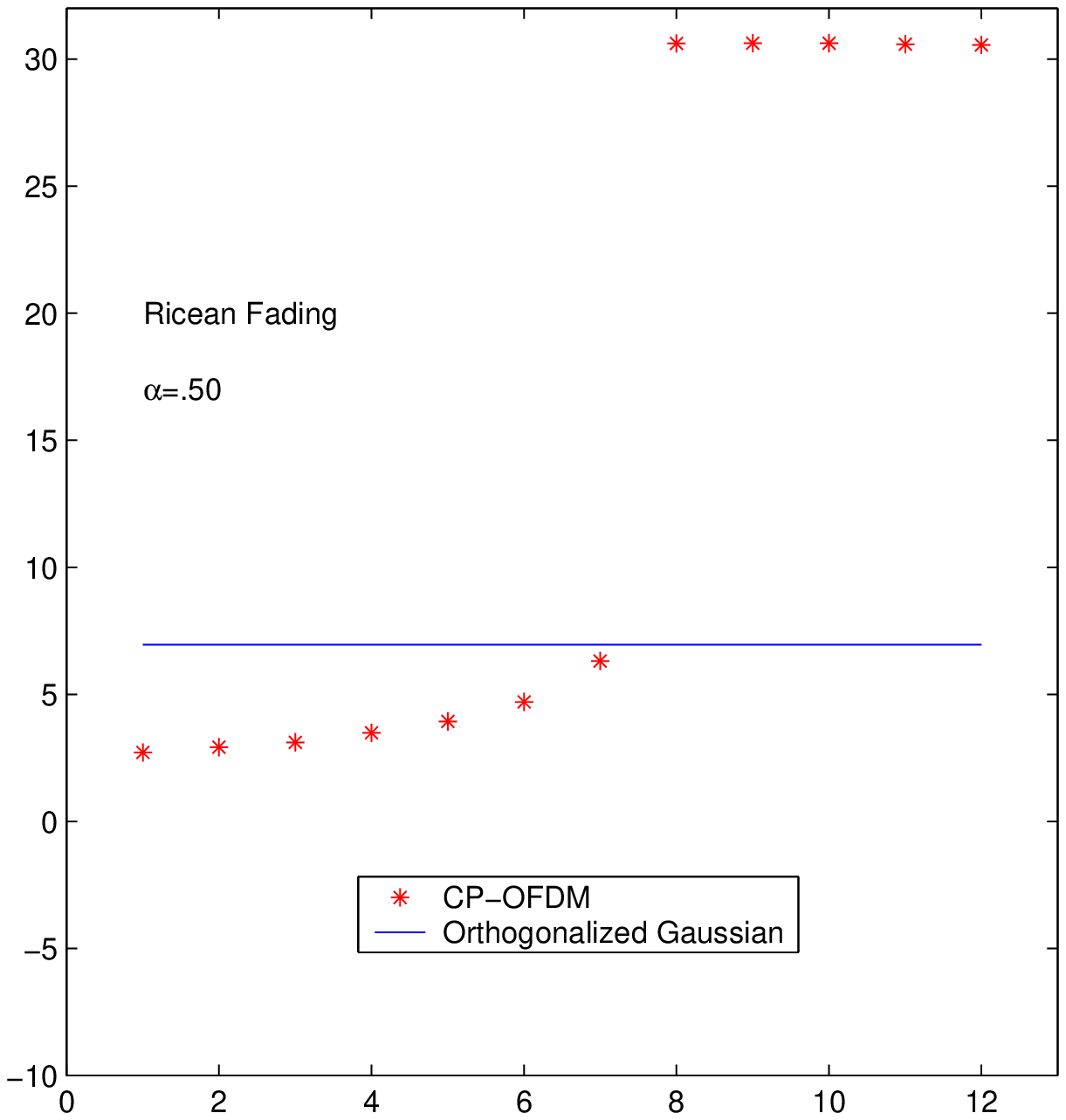,width=65mm,height=45mm}}
\subfigure[]{
\epsfig{file=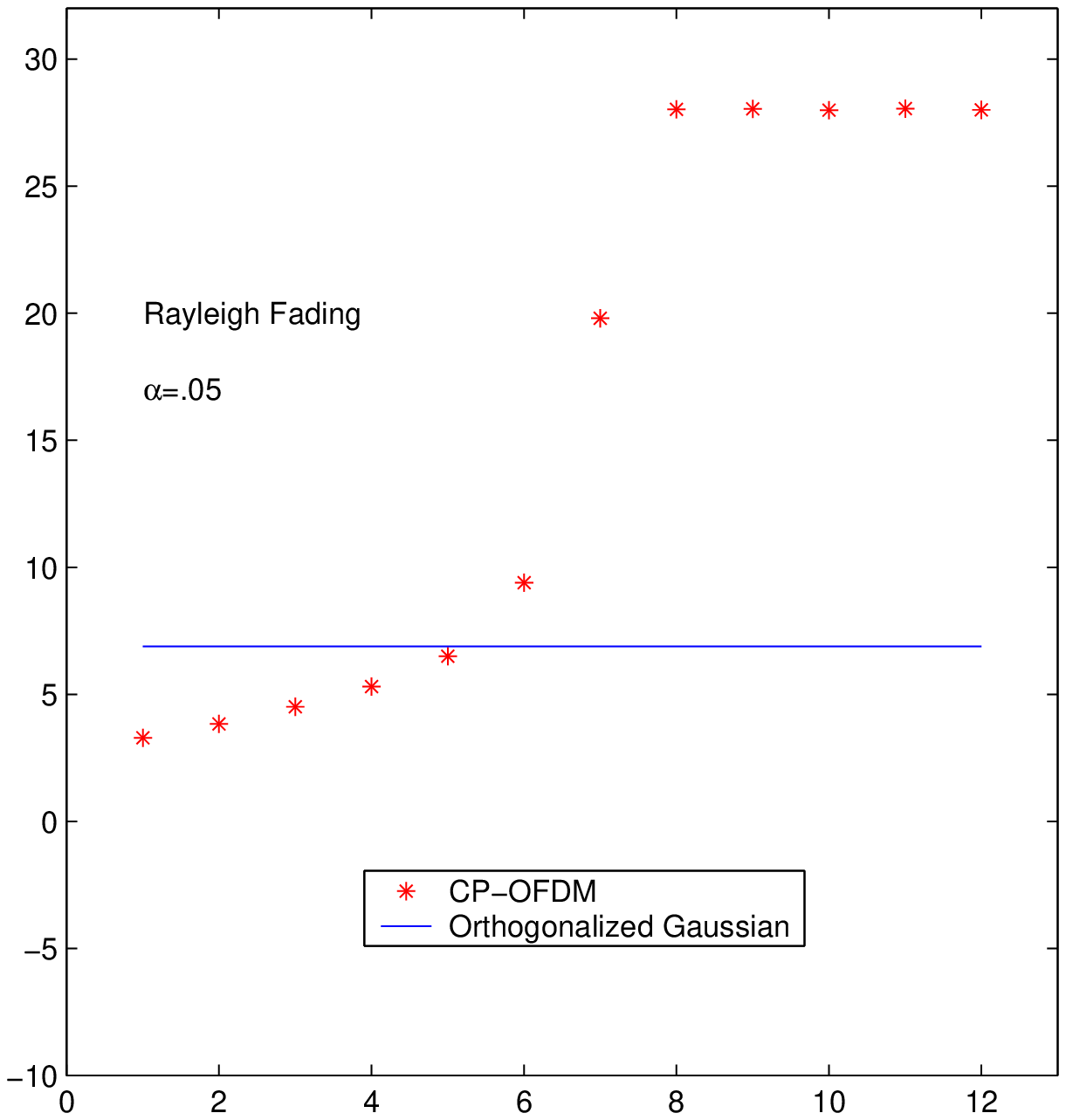,width=65mm,height=45mm}}
\subfigure[]{
\epsfig{file=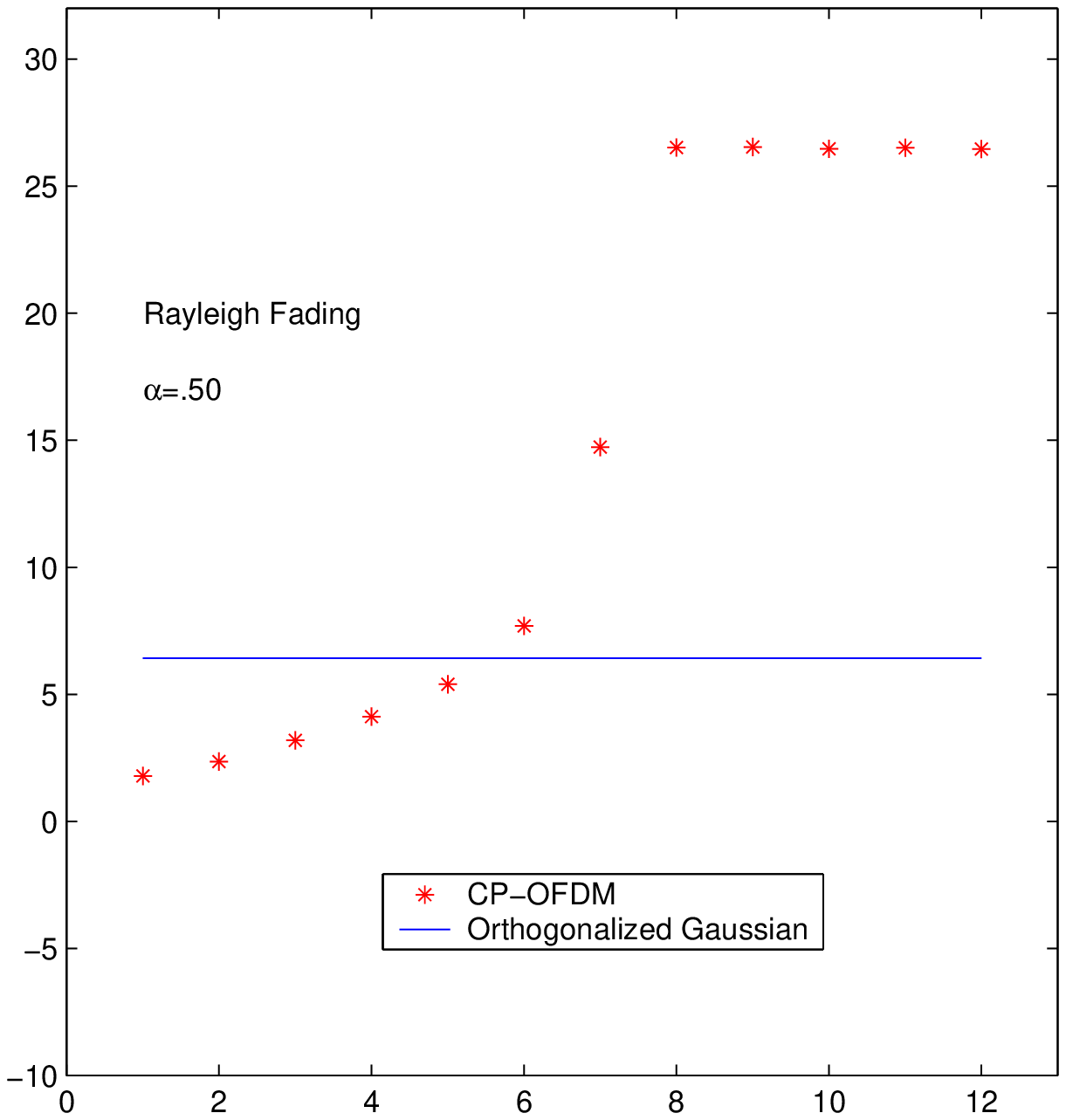,width=65mm,height=45mm}}
\caption{SNIR vs. cyclic prefix length for both rapid ($\a=.05$) and slow ($\a=.5$) temporal decay of channel scattering function, under both Ricean (LOS-dominant) fading and Rayleigh fading.  The SNIR of the L\"owdin-orthogonalized Gaussian pulse $\ph$ is shown by the solid line.  The noise level is $10^{-3},$ and $\t_0=8$ but $\nu_0=0.$}
\label{fig:rayfadeofdm}
\end{center}
\end{figure}

Note that even under the absence of the Doppler effect, $\ph,$ optimally adapted to channel conditions where $\t_0=\nu_0,$ usually outperforms the square pulse with cyclic prefix if the cyclic prefix length is smaller than $\t_0.$  If however, we introduce the condition that $\nu_0=\t_0,$ CP-OFDM is defenseless.\\

\begin{figure}[!ht]
\begin{center}
\subfigure[]{
\epsfig{file=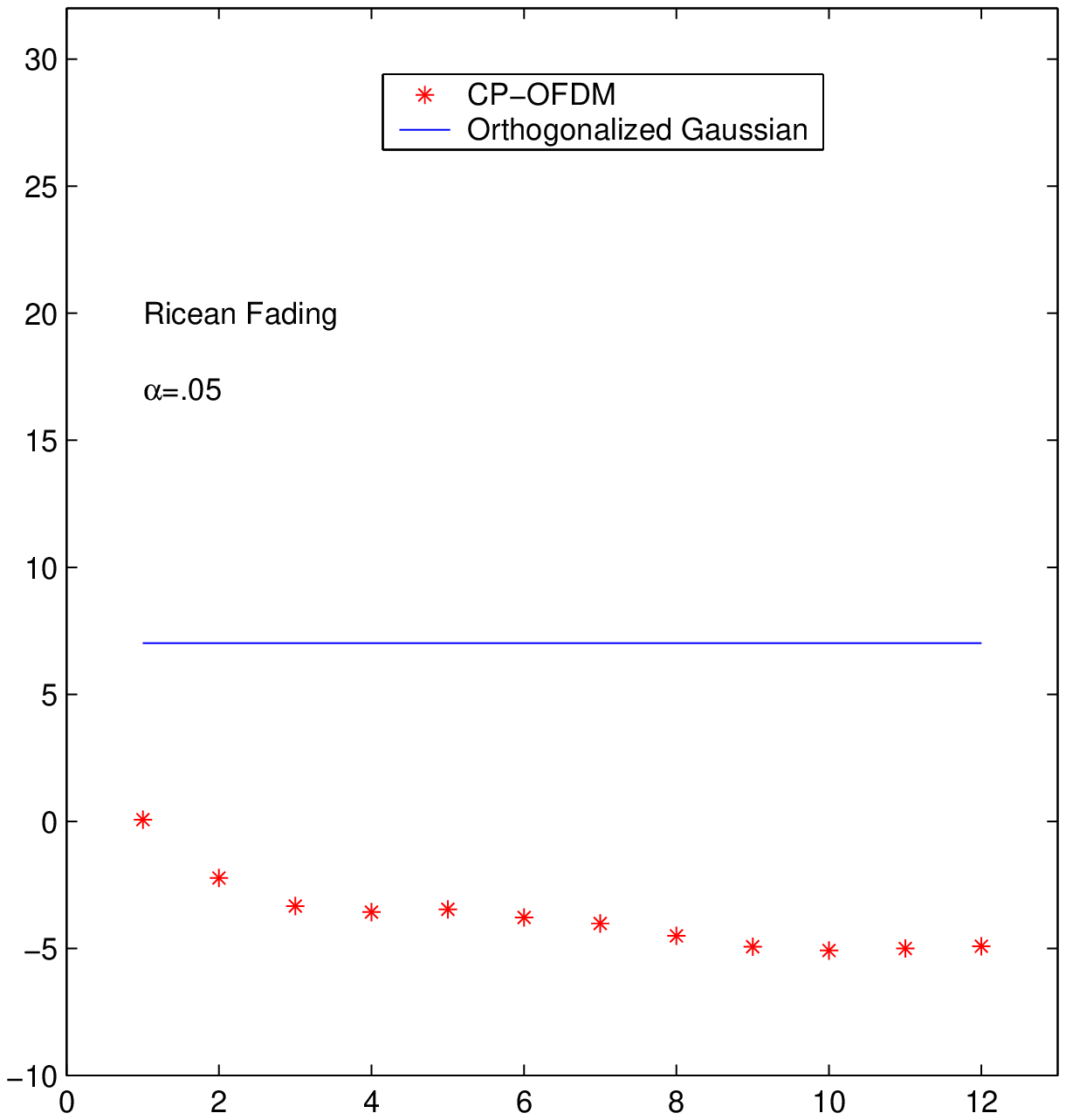,width=65mm,height=45mm}}
\subfigure[]{
\epsfig{file=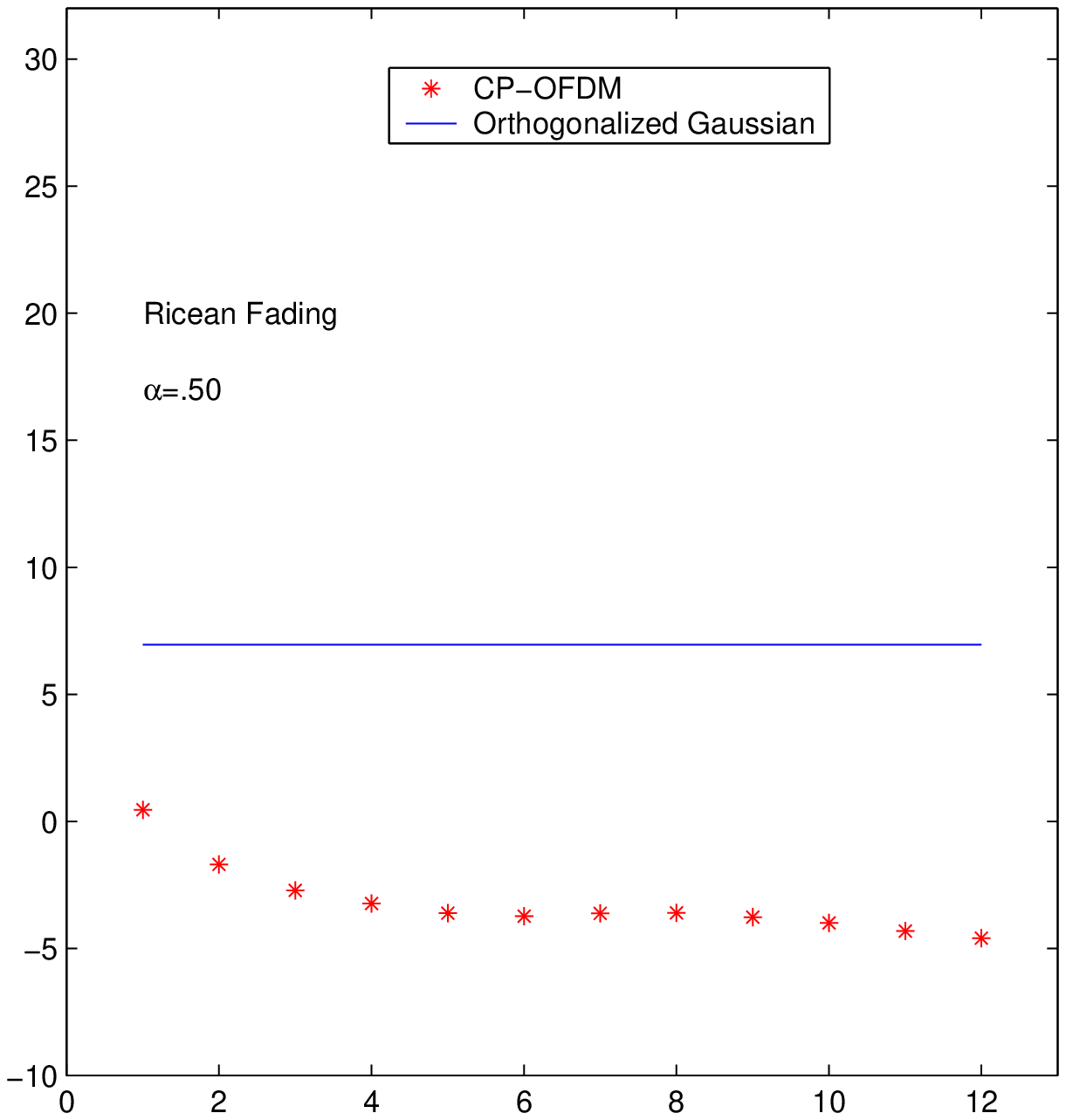,width=65mm,height=45mm}}
\subfigure[]{
\epsfig{file=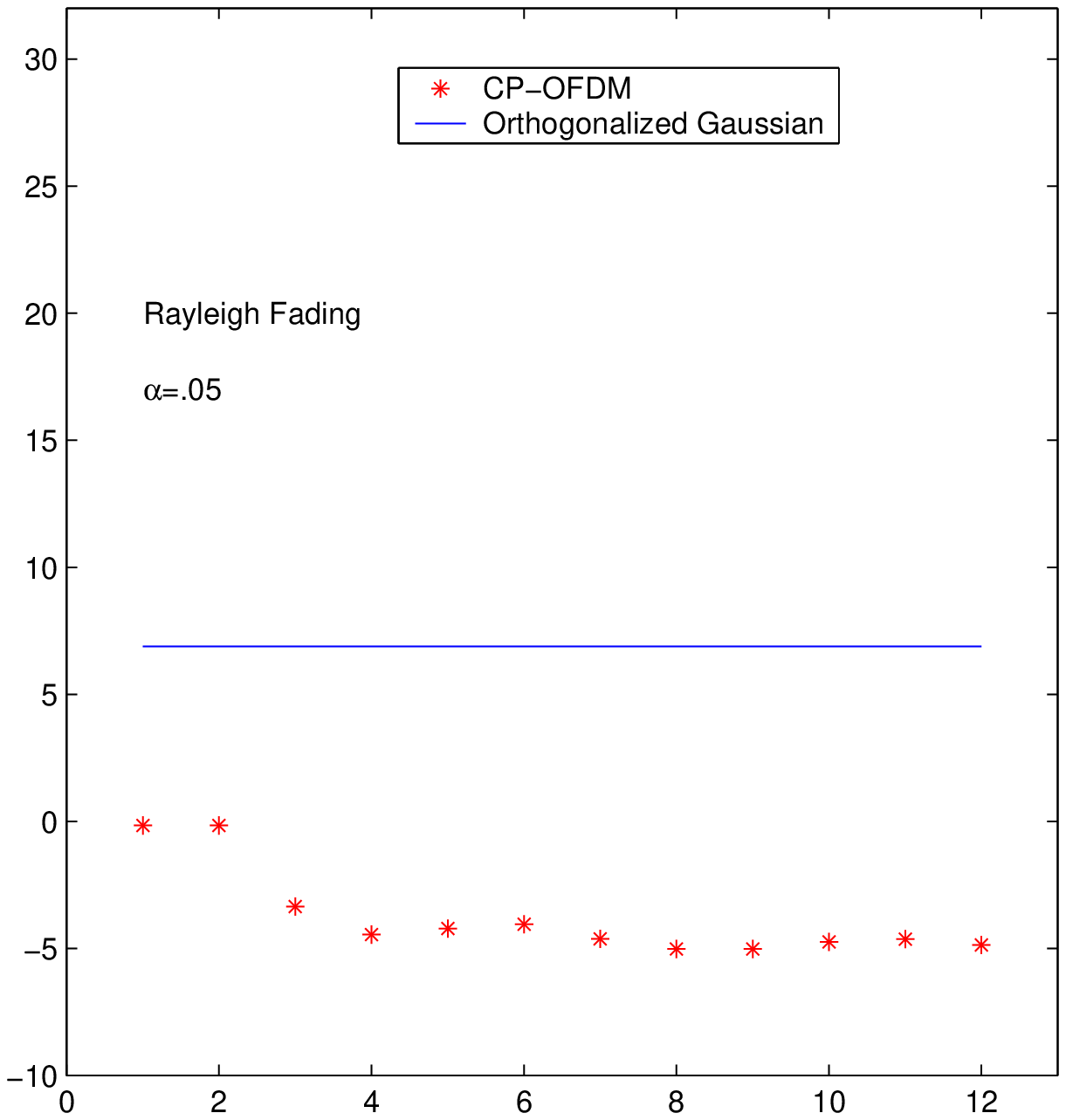,width=65mm,height=45mm}}
\subfigure[]{
\epsfig{file=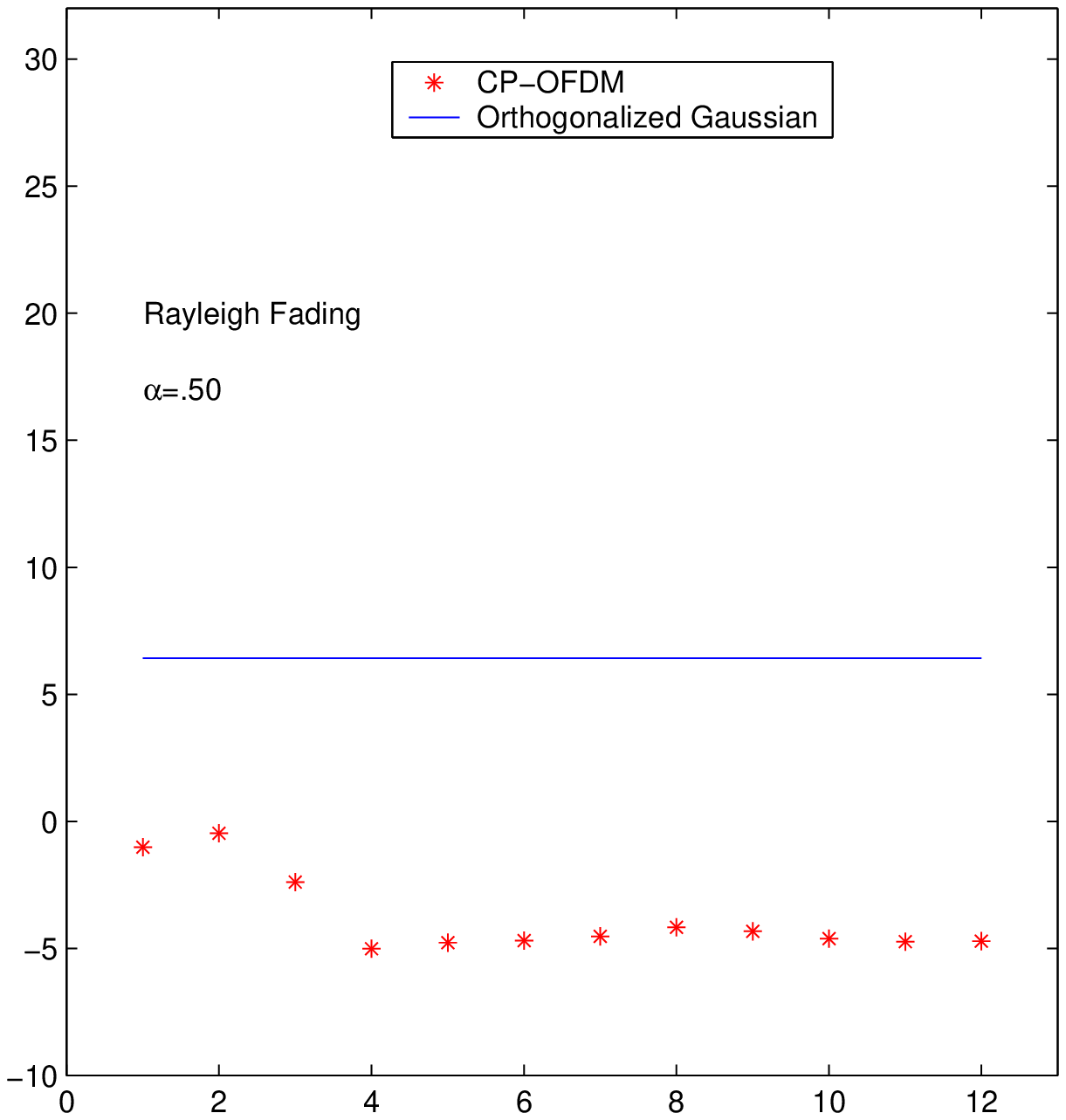,width=65mm,height=45mm}}
\caption{The same parameters as in Figure \ref{fig:rayfadeofdm}, except $\nu_0=\t_0=8.$}
\label{fig:rayfadeofdmnu}
\end{center}
\end{figure}
Where the Doppler effect is present, cyclic-prefix square-wave OFDM cannot mount any resistance to its frequency dispersion, and the resultant interference destroys the signal.  But the signal generated from $\ph$ wards off the dispersive effects of Doppler interference.

\section{Performance of EBFDM}\label{pbofdmcpofdm}
In this section we consider the performance of EBFDM systems as a function of the exponent $p.$  Is it well-known in communication theory that orthogonality ($p=\frac{1}{2}$) yields optimal robustness against AWGN.  

\begin{figure}[!ht]
\begin{center}
\epsfig{file=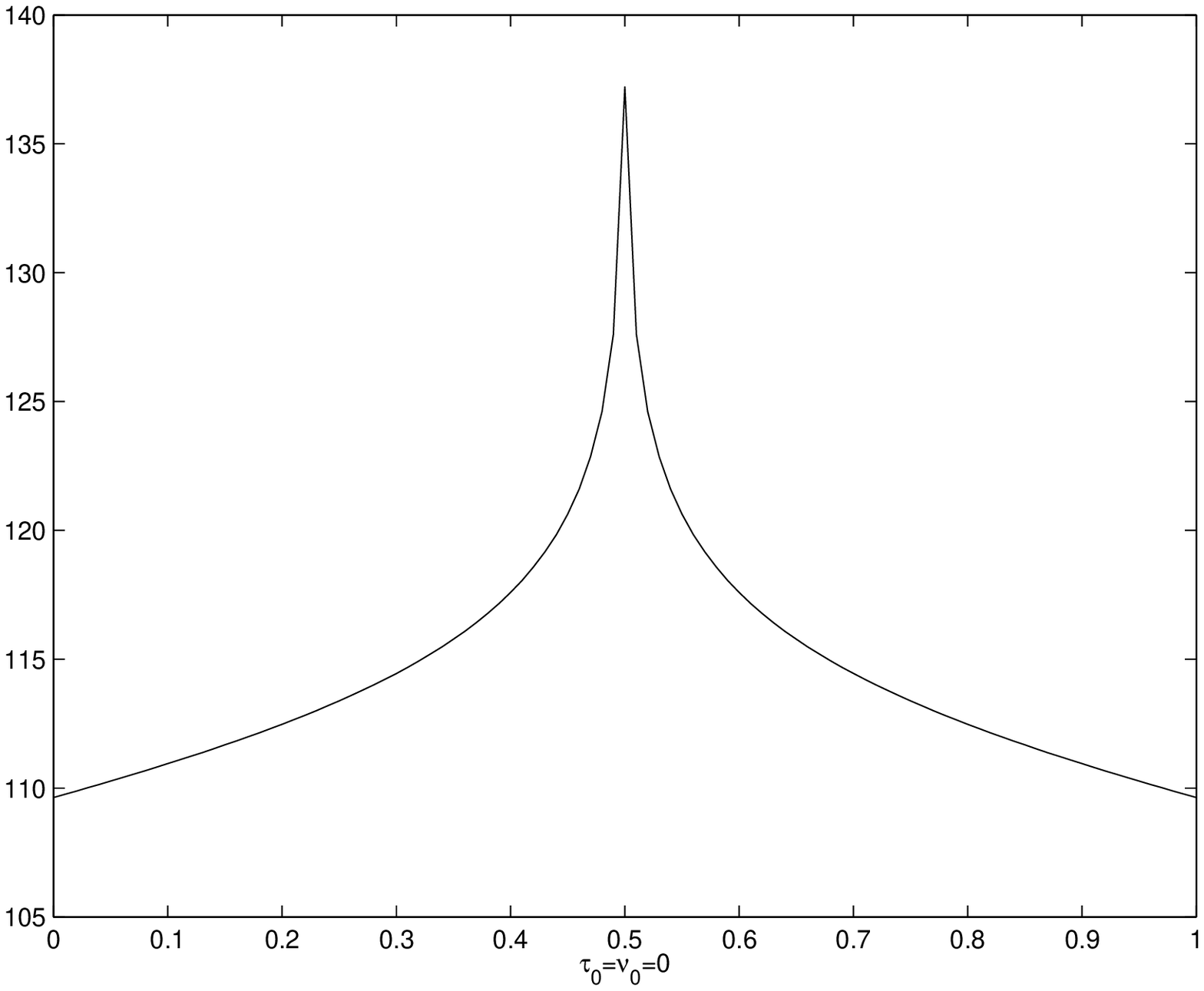,width=85mm}
\caption{SNIR vs. $p$ for noise level $10^{-16}$ and $\t_0=\nu_0=0.$}
\label{fig:nodisp1}
\end{center}
\end{figure}

\begin{figure}[!ht]
\begin{center}
\subfigure[]{
\epsfig{file=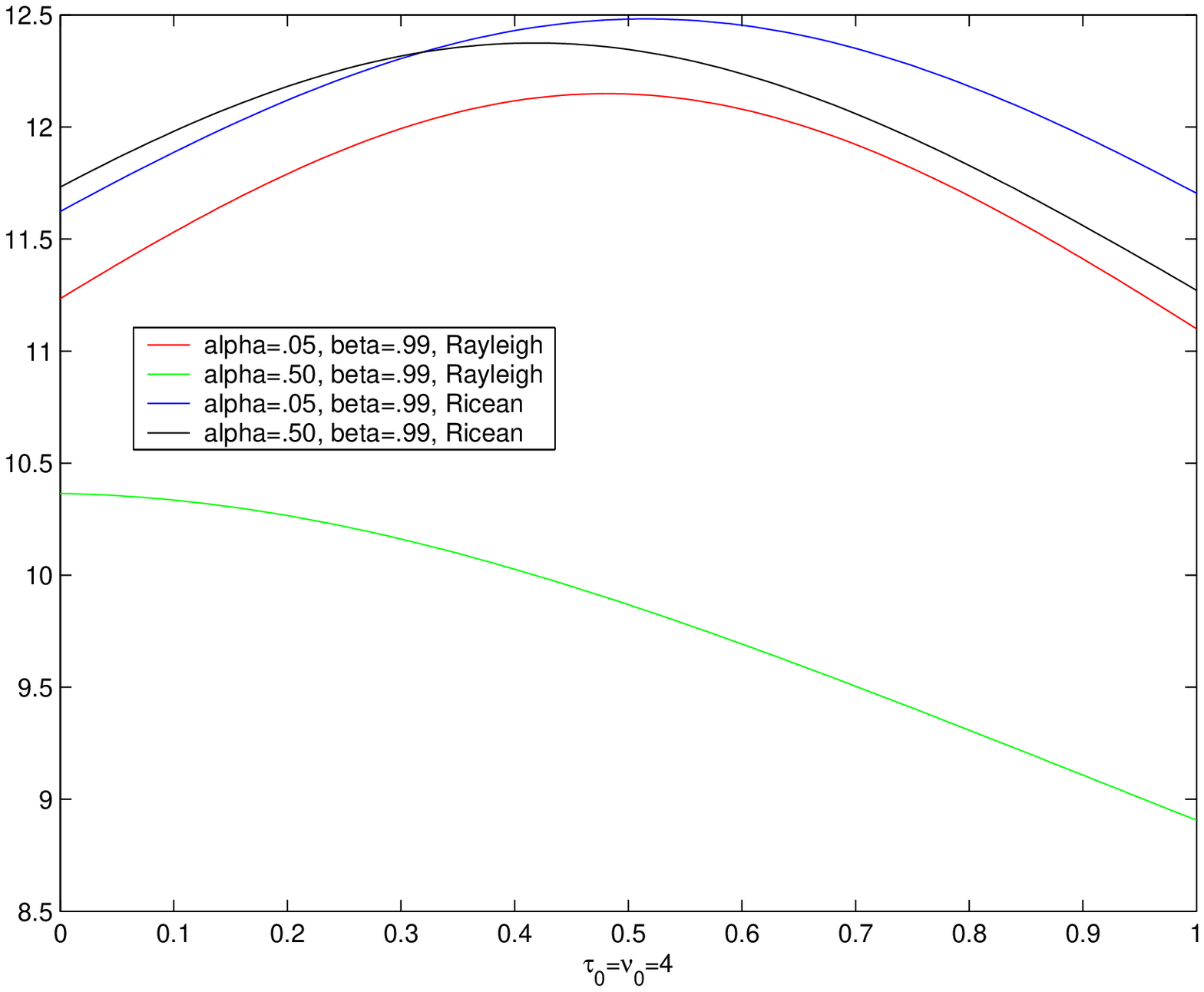,width=65mm,height=45mm}}
\subfigure[]{
\epsfig{file=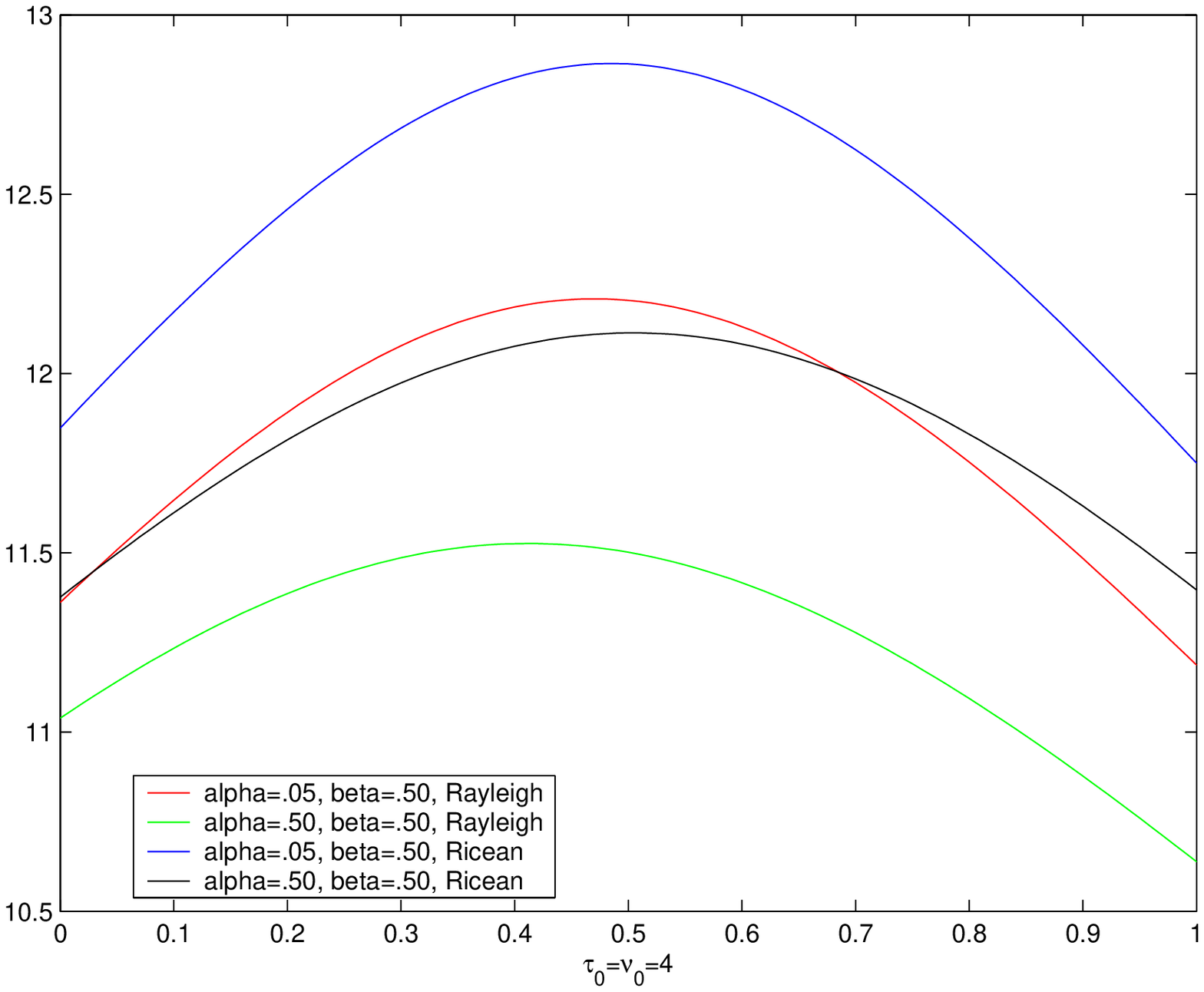,width=65mm,height=45mm}}
\subfigure[]{
\epsfig{file=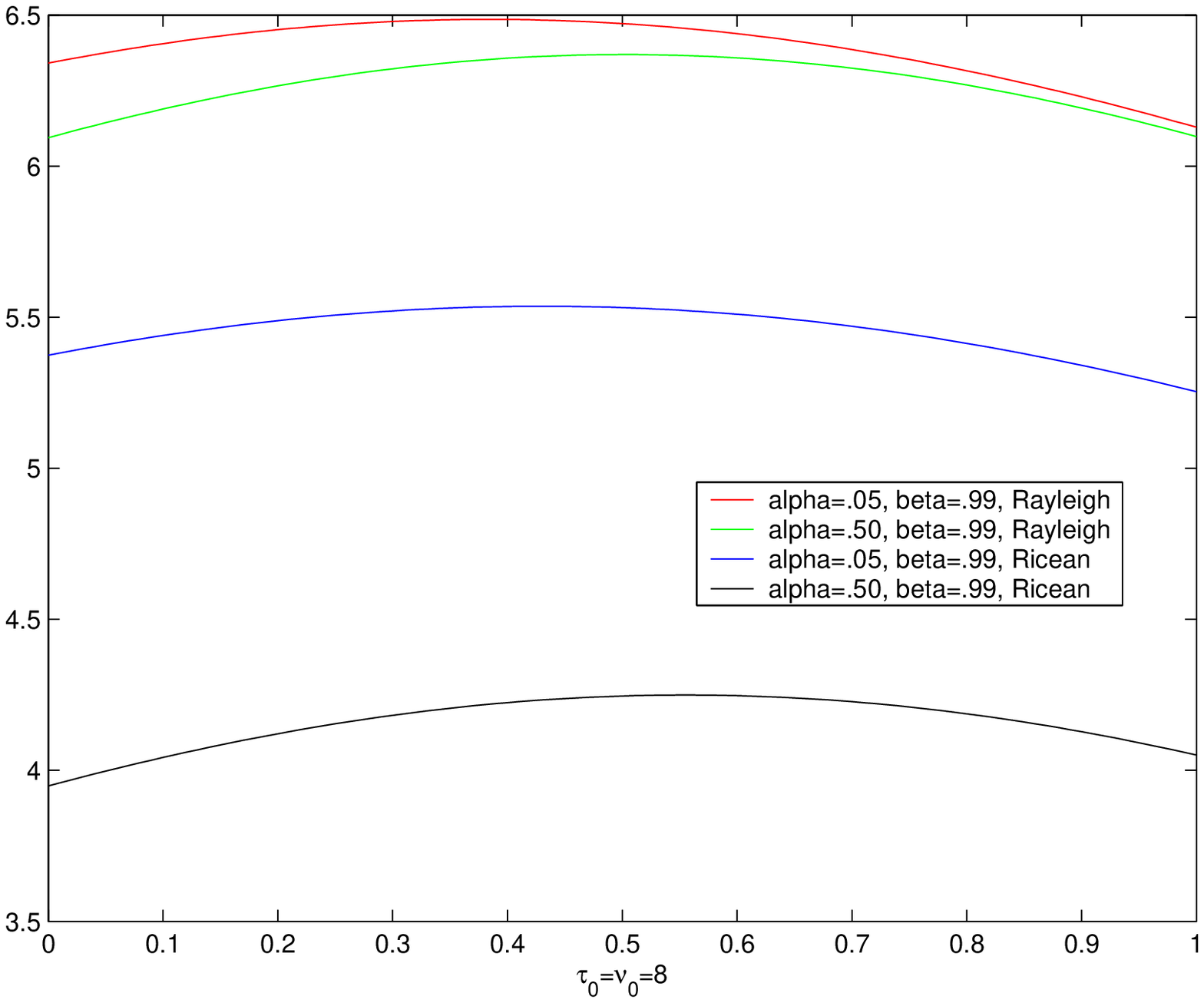,width=65mm,height=45mm}}
\subfigure[]{
\epsfig{file=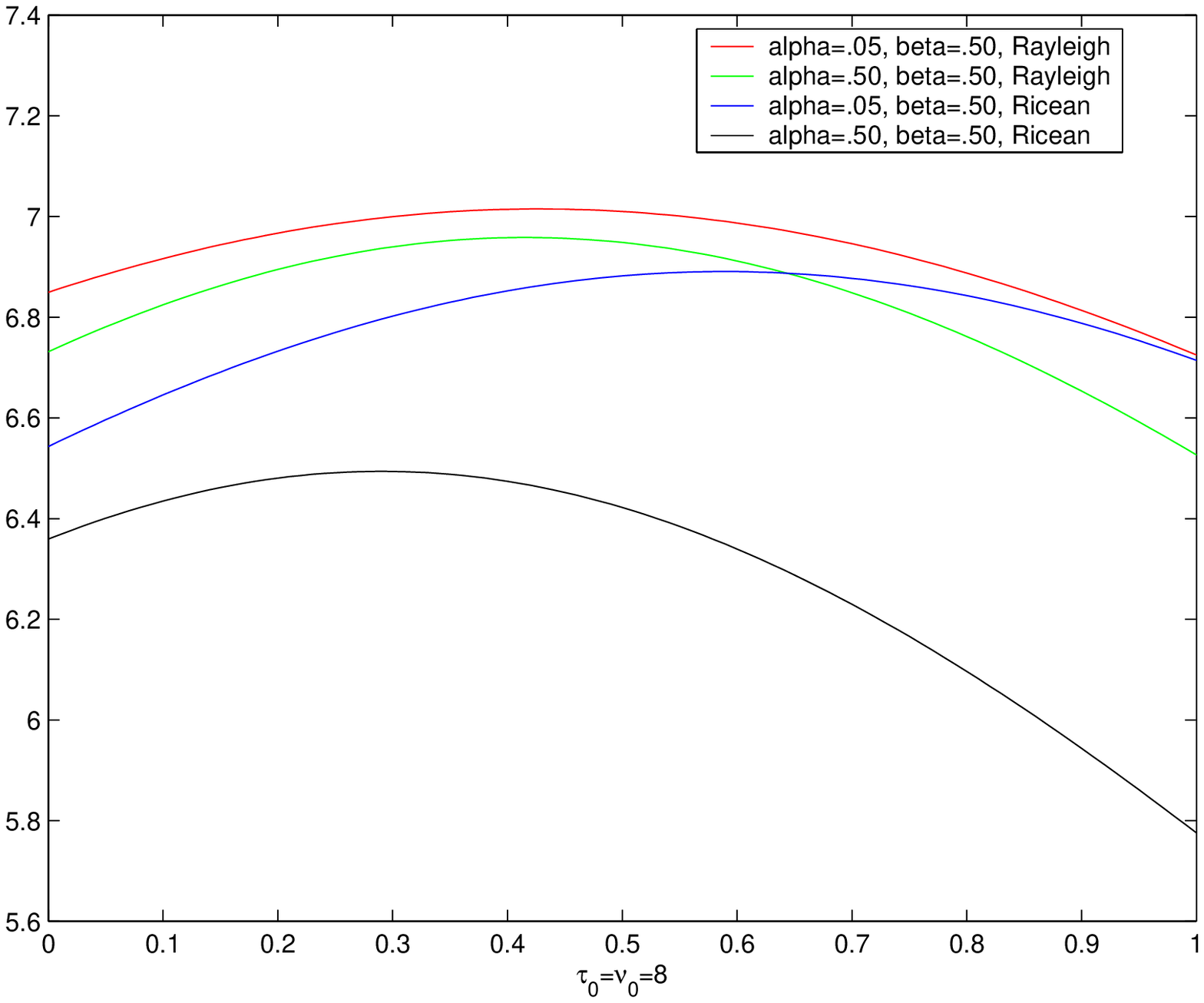,width=65mm,height=45mm}}
\caption{SNIR vs. $p$ for noise level $10^{-6}$ and $\t_0=\nu_0=4$ (Figures \ref{fig:psbfdm2} (a) and (b)) and $\t_0=\nu_0=8$ (Figures \ref{fig:psbfdm2} (c) and (d)) for both rapid ($\a=.05$) and slow ($\a=.5$) temporal decay of channel scattering function, and wide (Figures \ref{fig:psbfdm2} (a) and (c)) and narrow (Figures \ref{fig:psbfdm2} (b) and (d)) distributions of Doppler shifts.}
\label{fig:psbfdm2}
\end{center}
\end{figure}

\begin{figure}[!ht]
\begin{center}
\subfigure[]{
\epsfig{file=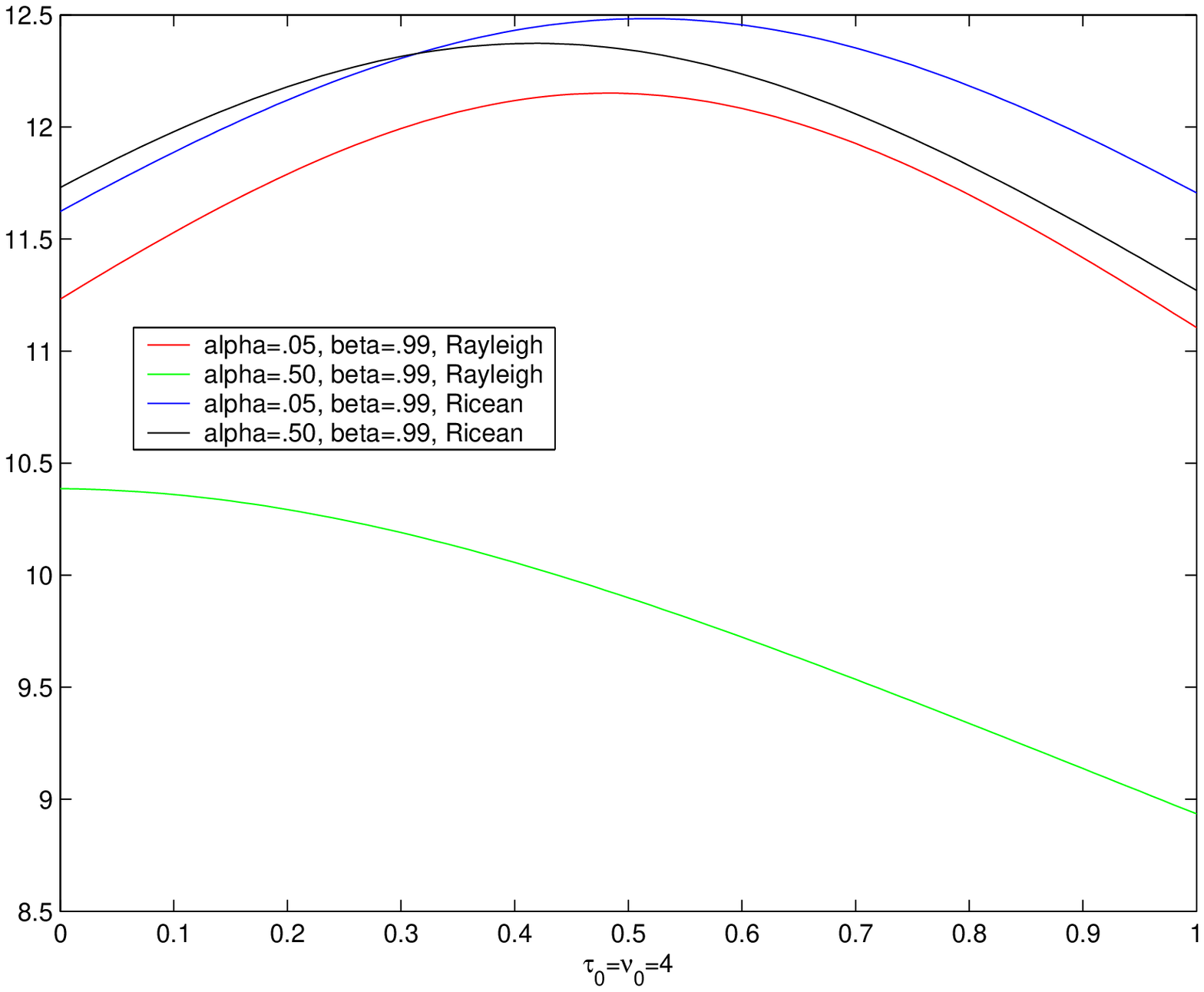,width=65mm,height=45mm}}
\subfigure[]{
\epsfig{file=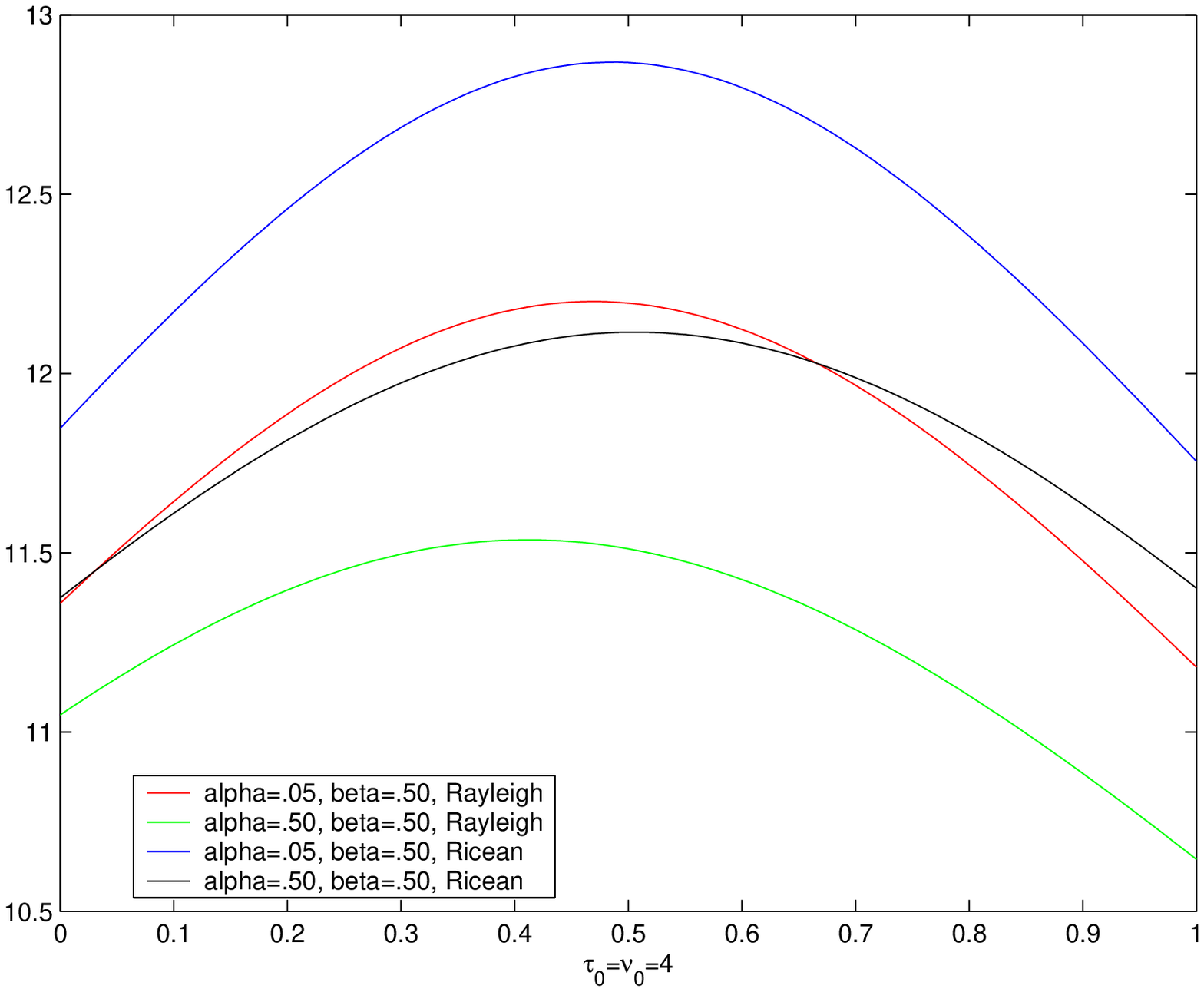,width=65mm,height=45mm}}
\subfigure[]{
\epsfig{file=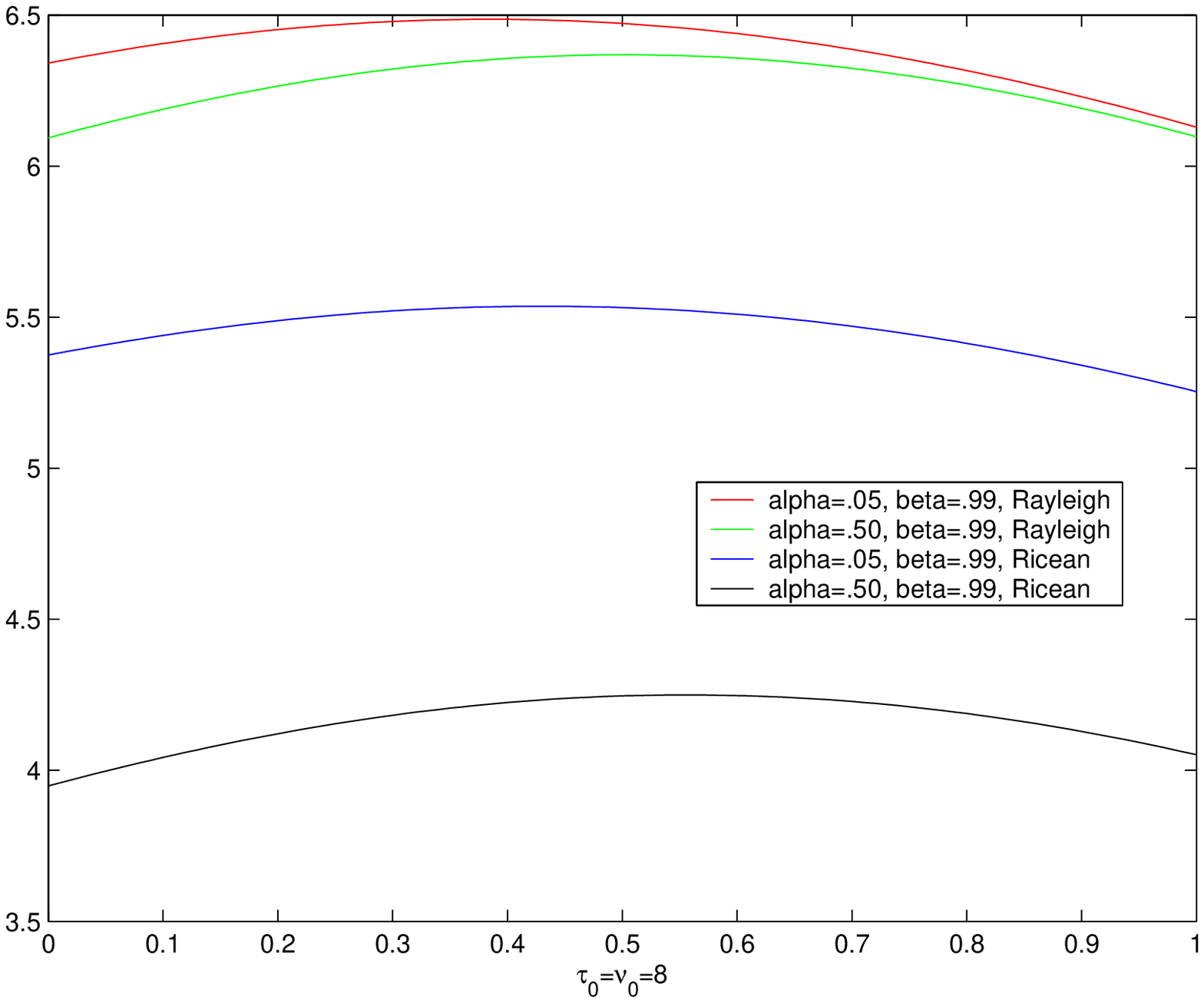,width=65mm,height=45mm}}
\subfigure[]{
\epsfig{file=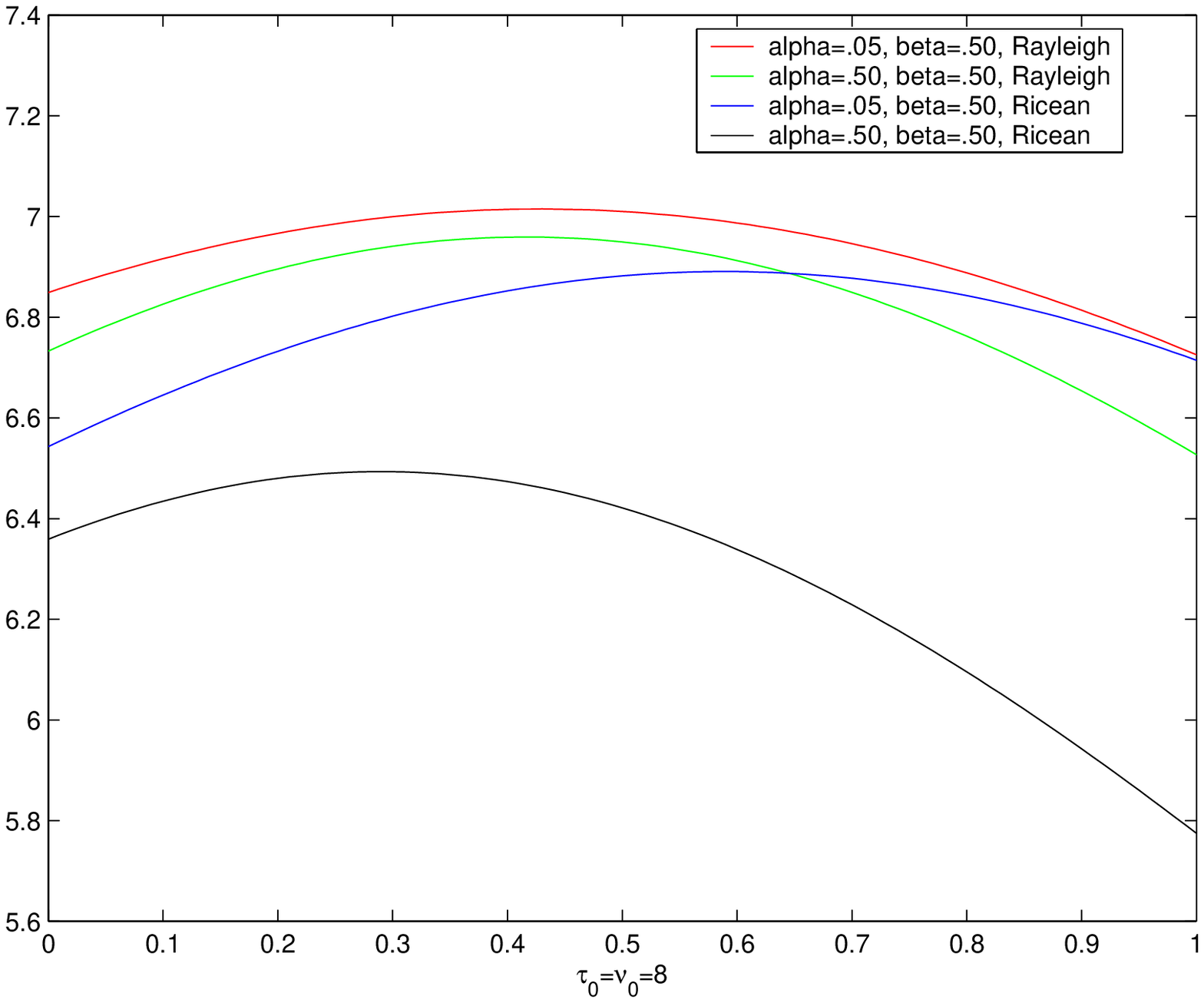,width=65mm,height=45mm}}
\caption{The same parameters as for Figure \ref{fig:psbfdm2} except that the  noise level is $10^{-3}.$}
\label{fig:psbfdm1}
\end{center}
\end{figure}

However, for nonzero dispersion the deviation from optimality of $p=\frac{1}{2}$ for nonzero dispersion is often small, but thus far unexplained.  Random effects are a possible culprit, although lengthening the signal has produced essentially the same results.  Additionally, it is not yet clear what, if any, dependence the SNIR has on $\a$ or $\b,$ nor is clear the difference in character between the cusp-like maximum in Figure \ref{fig:nodisp1} and the apparently smooth maxima in e.g., Figure \ref{fig:psbfdm2}; additionally the local concavity in Figure \ref{fig:nodisp1} fails to mimic the convexity in e.g., Figure \ref{fig:psbfdm2}.  We note that in Figures \ref{fig:psbfdm2} and \ref{fig:psbfdm1}, the difference in SNIR for the two noise levels above is negligible.  But for a noise level of $10^{-1}$, the noise clearly plays a significant role in lowering the values of the SNIR.

\begin{figure}[!ht]
\begin{center}
\subfigure[]{
\epsfig{file=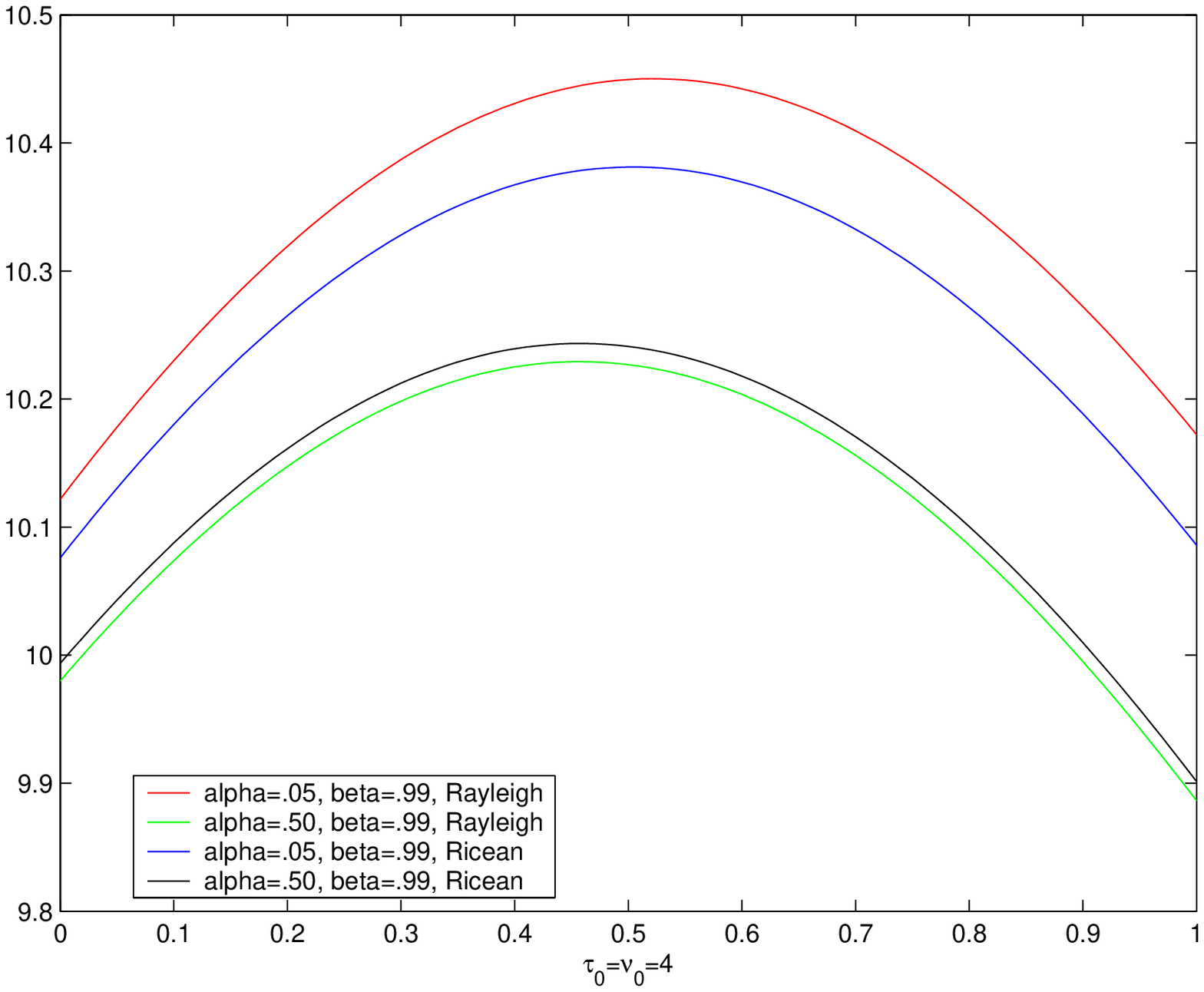,width=65mm,height=45mm}}
\subfigure[]{
\epsfig{file=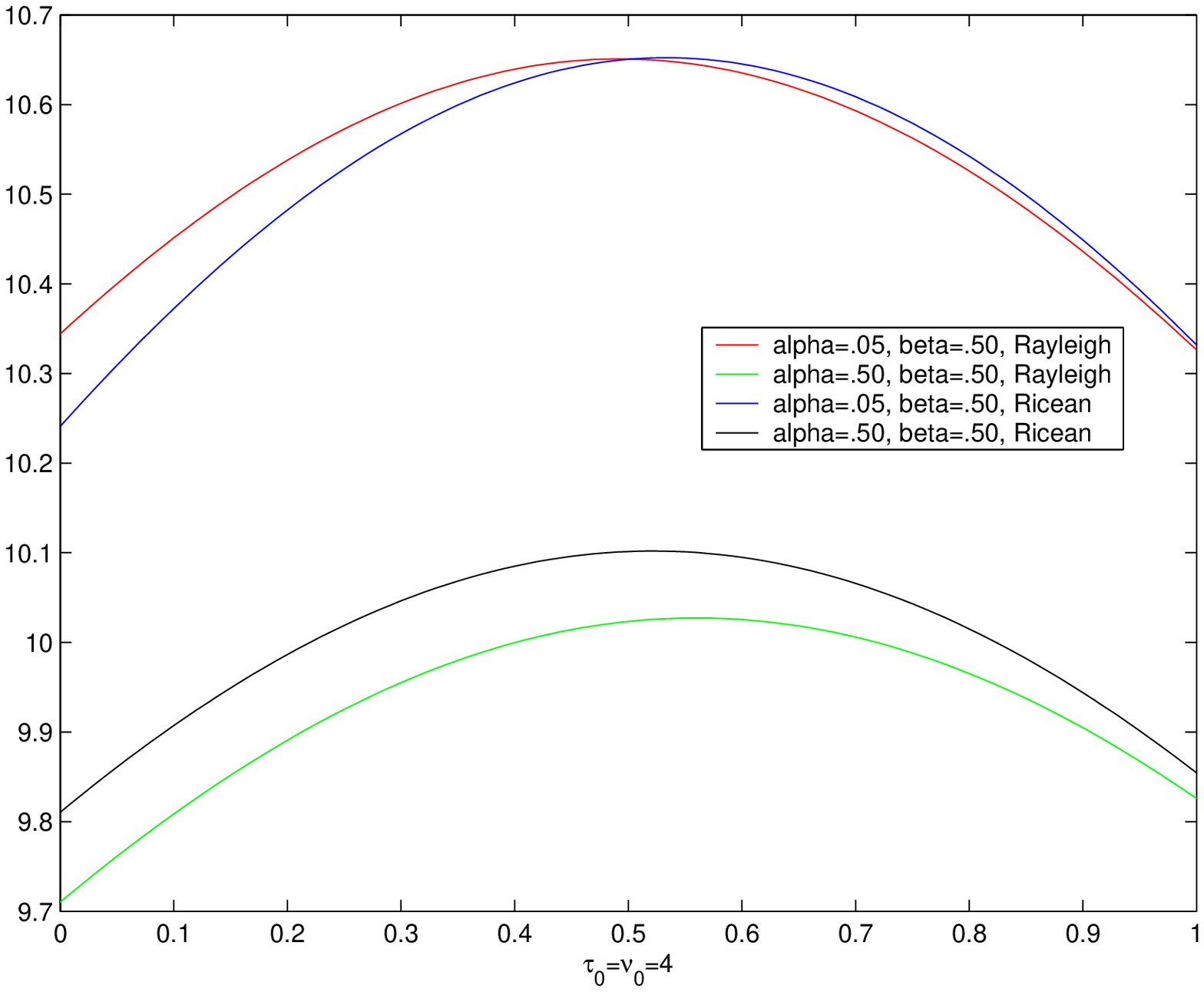,width=65mm,height=45mm}}
\subfigure[]{
\epsfig{file=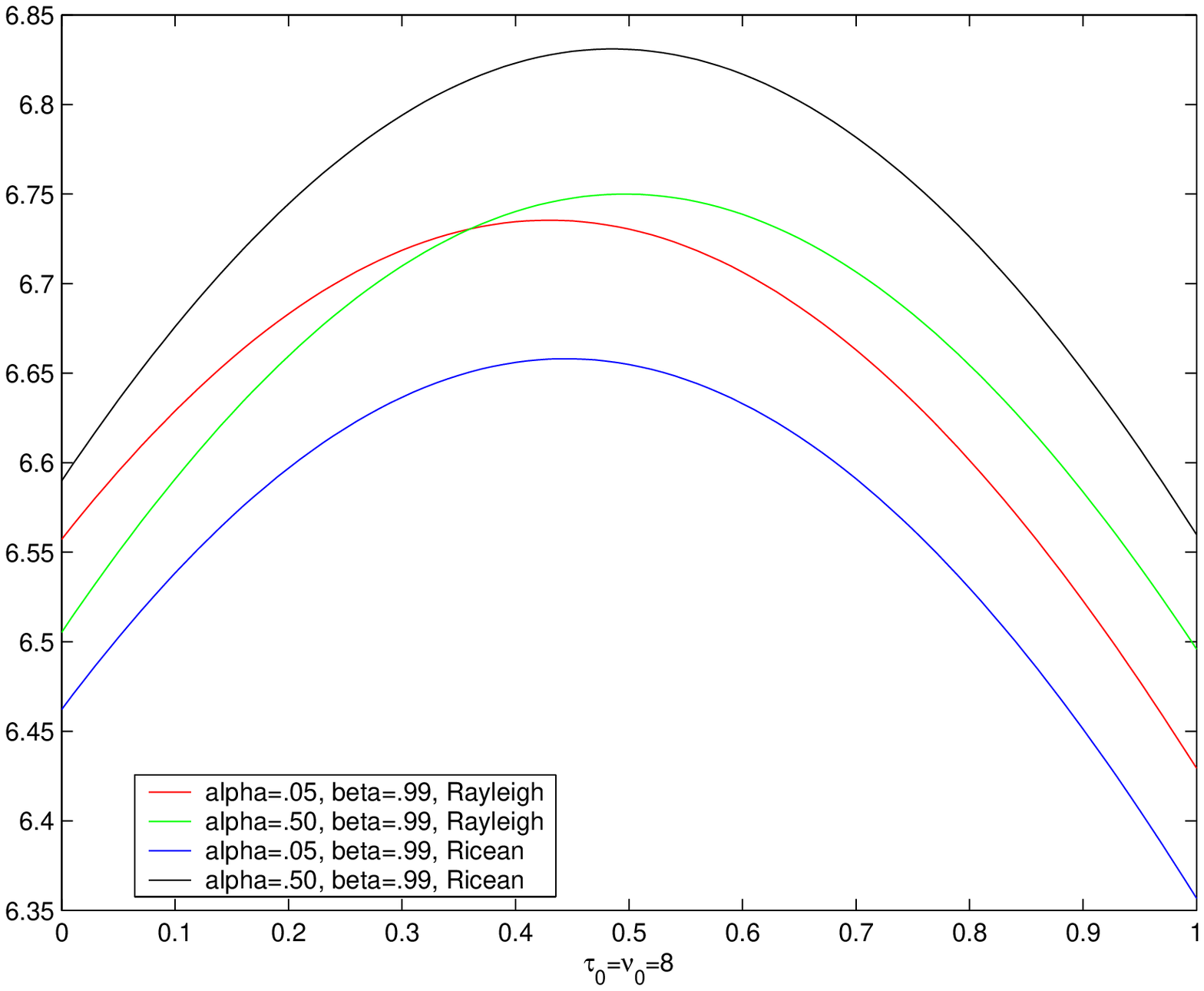,width=65mm,height=45mm}}
\subfigure[]{
\epsfig{file=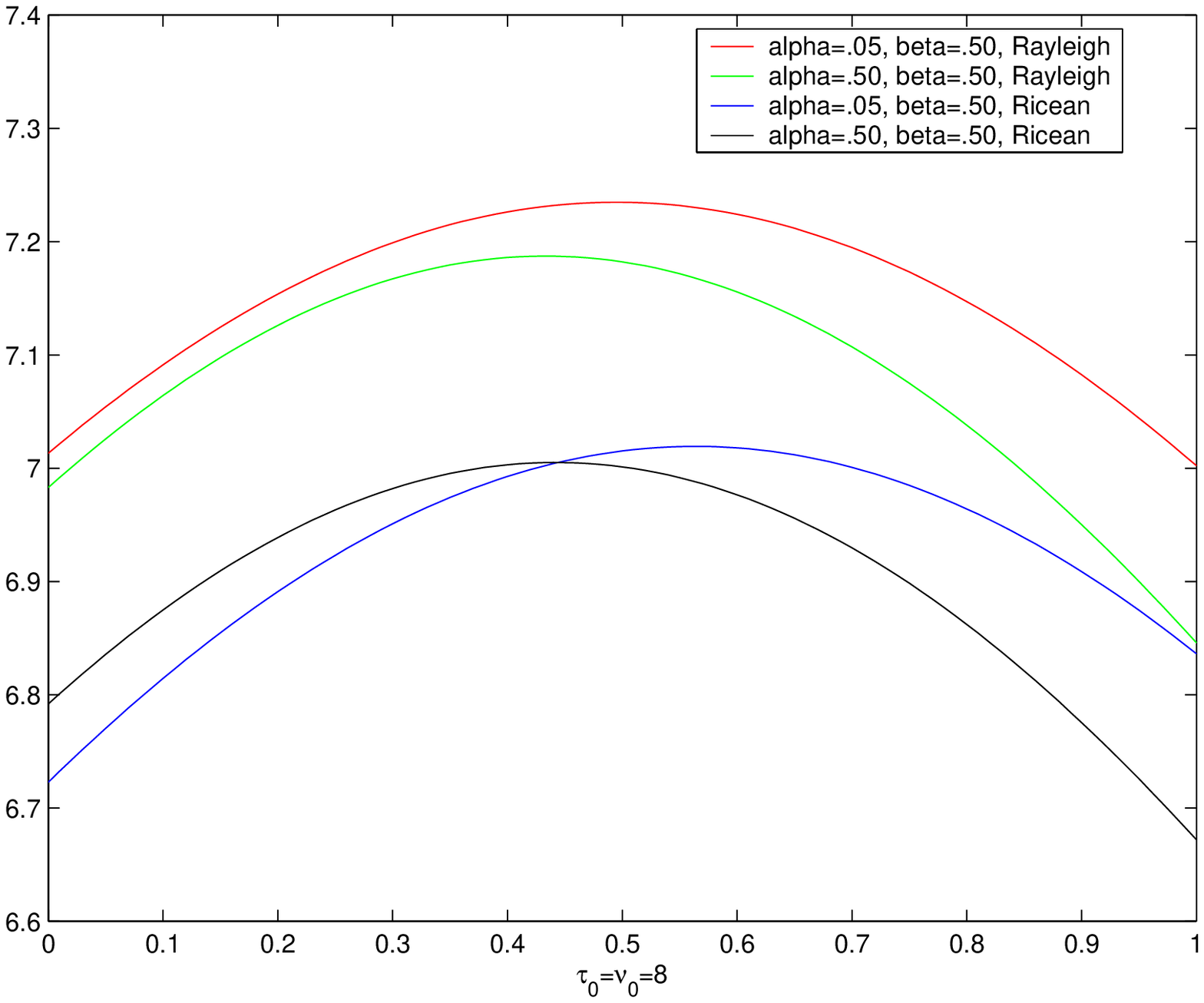,width=65mm,height=45mm}}
\caption{The same parameters as for Figures \ref{fig:psbfdm1} and \ref{fig:psbfdm2} except that the  noise level is $10^{-1}.$}
\label{fig:psbfdm3}
\end{center}
\end{figure}

\np
\pagestyle{myheadings}
\chapter{A General Framework for Pulse Shape Design}\label{GenPulseDesign}
\thispagestyle{myheadings} 
\markright{  \rm \normalsize CHAPTER 7. \hspace{0.5cm}
A GENERAL FRAMEWORK FOR PULSE SHAPE DESIGN}

\section{Adaptation of the pulse to channel conditions}\label{adapttochannel}
It has often been suggested that the time duration and frequency duration of the transmitter functions should be adapted to the channel properties, described by the channel scattering function $C_{\mbox{\boldmath $H$}}(x,\omega).$ The ambiguity function $A\psi_{k,\l}(x,\omega)$ is the inner product of a receiver function centered at $(kT,\l F)$ and a channel-scattered transmitter function centered at $(kT \pm x,\l F \pm \omega).$  Thus $A\psi_{k,\l}$ describes the ability of an OFDM system to detect data symbols in the presence of time- and frequency-dispersion.  For a function $f$ and $\eps >0$ define the {\em effective support} of the
ambiguity function $Af$ by
\begin{equation}
\supp_{\eps}(A f) = \{(t,\omega): |A f(t,\omega)| > \eps\}.
\label{supp}
\end{equation}

\begin{figure}[htb]
\begin{center}
\epsfig{file=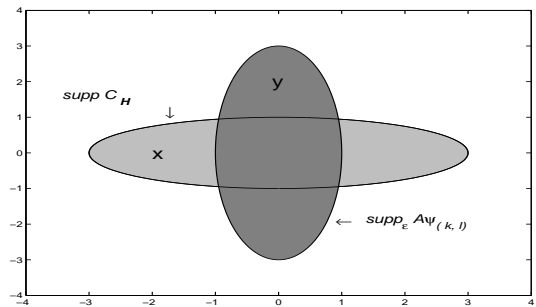,width=70mm,height=40mm}
\caption{Ambiguity function poorly matched to channel scattering function.}
\label{fig:scamb}
\end{center}
\end{figure}

Consider Figure \ref{fig:scamb}.  In region $Y,$ $|A\psi_{k,\l}|\geq \epsilon,$ but the channel scatters no input to $Y$ so designing a receiver that can detect a $Y$-scattered pulse is wasteful. Conversely, if the receiver cannot detect a pulse scattered to the region $X,$ data will be lost. Thus we should adapt the pulse so that $\supp_{\eps} A \psikl$ coincides, at least roughly, to the support of the channel scattering function.  \\

In order to reduce ISI and ICI it often has been suggested to choose the time duration $\Dtf$ and the frequency duration $\Dnf$ of the transmission pulse as well as the parameters $F$ and $T$ such that they are adapted to the channel properties, cf. e.g. \cite{FAB95,HB97,KM98}
and for a more theoretical approach see \cite{Koz98}.  This means that for given maximum multipath delay $\tau_0$ and maximum Doppler spread $\nu_0$ 
we choose $\Dtf$, $\Dnf$, $T$, and $F$ such that 
\be
\frac{\tau_0}{\nu_0} \approx \frac{\Dtf}{\Dnf} \approx \frac{T}{F}.
\label{adapt}
\end{equation}
For the Gaussian-based window $\psigam$ this means that for given spectral efficiency $\rho,$ we set $\sigma \approx \frac{\nu_0}{\tau_0}$ and $T=\sqrt{\frac{1}{\sigma \rho}}$,
$F=\sqrt{\frac{\sigma}{\rho}}.$  Using basic properties of the dilation operator $\Dil$ we obtain the relation 
\begin{equation}\label{orthdil}
\psigam=L\ddot{o}\,(\ggam,1/\sqrt{\sigma \rho},\sqrt{\sigma/\rho})=
{\cal D}_{\sqrt{\sigma}} \big(L\ddot{o}\,(g,1/\sqrt{\rho},1/\sqrt{\rho})\big)=
{\cal D}_{\sqrt{\sigma}} \ph.
\end{equation}
A detailed proof of \eqref{orthdil} is in Appendix \ref{dilorth}.\\

Thus the transmission pulse $\psigam$ that is matched to the channel parameters
$\tau_0$ and $\nu_0$ is obtained from the ``standard'' pulse $\ph$
by an appropriate dilation.  Hence at least theoretically it is very easy to adapt the OFDM system to time-varying channel conditions. We only have to know how the 
ratio of $\tau_0$ and $\nu_0$ changes and adapt the lattice and $\psigam$ 
accordingly by a simple dilation. Of course, in practice one must also deal with issues such as time and frequency synchronization when one changes the lattice parameters, which may be a non-trivial task.

\section{Lattice-OFDM, wireless channels, and interference}\label{s:lofdm}
The robustness of an OFDM system against ISI and ICI depends essentially
on two factors:
\begin{itemize}
\item[\quad (i)] the TFL of the functions $\psi_{kl}$;
\item[\quad (ii)] the distance between adjacent lattice points
$(kT,\l F)$, $((k\pm1) T,((l\pm1) F)$.
\end{itemize}
We have seen in Section \ref{adapttochannel} how to construct for given channel
conditions OFDM pulse shapes with optimal time-frequency concentration.
If we want to further improve the stability of the OFDM system against
interference we have to increase the distance between adjacent data
symbols by increasing the distance between adjacent grid points.
However increasing $T$ and/or $F$ results in an undesirable loss of spectral
efficiency.  Can we increase the distance between adjacent grid points
without reducing the spectral efficiency?\\

We will show in this section that this is indeed possible.  The solution lies in the introduction of a more general and flexible OFDM setup which we call {\em Lattice-OFDM} (LOFDM).

\subsection{Lattice-OFDM}
\label{ss:lofdm}

The generator matrix $L_R$ for the rectangular lattice 
$\La_R$ displayed in Figure \ref{fig:rect} is 
\begin{equation}
L_R=\begin{bmatrix} T & 0 \\0 & F 
\end{bmatrix},
\label{rectmat}
\end{equation}
the density of which is $\xi(\La_R) = 1/\det(L_R)=1/(TF)$. \\

The hexagonal lattice
$\La_H$ (shown in Figure \ref{fig:hex} with $T=F$) with generator matrix
\begin{equation}
L_H=
\begin{bmatrix}
\frac{\sqrt{2}}{\sqrt[4]{3}}T & \frac{\sqrt{2}}{2\sqrt[4]{3}}T \\
0 & \frac{\sqrt[4]{3}}{\sqrt{2}} F
\end{bmatrix}
\label{hexmat}
\end{equation}
has the same density, viz.\ $\xi(\La_H)=1/\det(L_H)=1/(TF)$.

\begin{definition}
Let $\La$ be a lattice with generator matrix
\begin{equation}
\label{genmat}
L=
\begin{bmatrix} 
x & y \\
0 & z 
\end{bmatrix},\qquad  x,z \neq 0, 
\end{equation}
and let $\phi \in \LtR$ be given. As in \eqref{Gaborsystem},we denote by $(\phi,\La)$ the 
function system
\begin{equation}
\label{lofdm}
\phi_{k\l}(t)=\phi(t-\lambda_k) e^{2\pi i \mu_\l t},
\qquad (\la_k,\mu_\l) \in \La, 
\quad k \in \Z,\, \l\!=\!0,\dots,N\!-\!1 \,\text{or}\,\, \l \in \Z.
\end{equation}
If the functions $\phi_{k\l}$ are mutually orthonormal we call 
$(\phi, \La)$ a Lattice-OFDM (LOFDM) system.\footnote{We have restricted in the definition of LOFDM to lattices with upper-triangular generator matrices, instead of general invertible
$2\times 2$ generator matrices. The reason is that unlike general lattices those with upper-triangular generator matrix - though not necessarily separable - still possess a distinctive orientation along the time-axis and the frequency-axis respectively. This orientation is very useful for an efficient numerical implementation (see also Section \ref{ss:eff}) and it keeps the time delay introduced by the non-separability of the lattice small. All results that follow could be derived (with some minor adaptations) for general lattices, however with little or no practical benefit.}
\end{definition}
It is clear that the spectral efficiency $\rho$ of an LOFDM system
$(\phi,\La)$ coincides with the density of the lattice;
i.e., $\rho = 1/\det(L)$.

\subsection{Sphere-packings and minimization of ISI/ICI} \label{ss:sphere}
In this section we show how to design pulse shapes and lattices
which provide optimal robustness for time-frequency dispersive channels.
We start with a specific but important example. The analysis of the general
case will be based on this example.\\

Assume that the parameters of the channel scattering 
function satisfy $\tau_0 = \nu_0$ (we will discuss more general
choices for $\tau_0, \nu_0$ later). By the results in Section \ref{s:pulse}
the optimal pulse shape $\ph$ in this case is given by 
$\ph = L\ddot{o}\,(g,T,F)$ with $g(t)=2^{\frac{1}{4}} e^{-\pi t^2}$
(since $\sigma=\frac{\nu_0}{\tau_0}=1$) and $T=F$. \\

\begin{figure}[!ht]
\begin{center}
\subfigure[]{
\epsfig{file=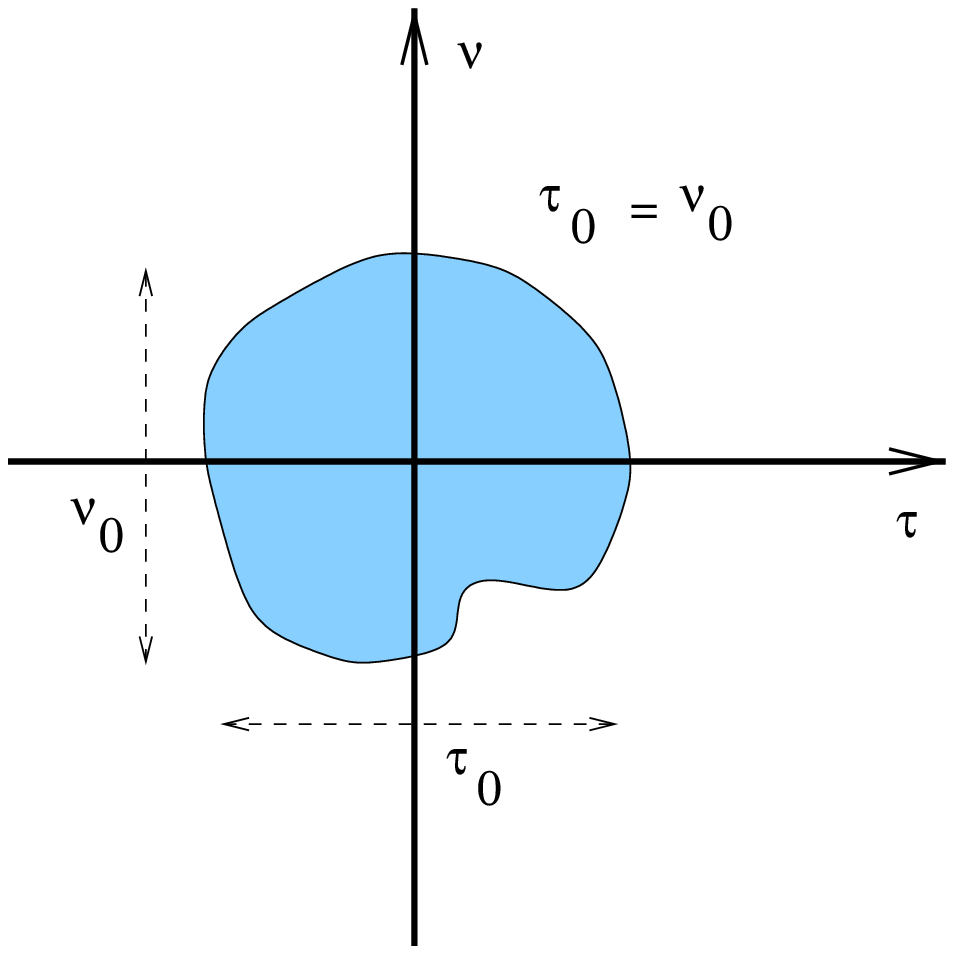,width=40mm}}
\subfigure[]{
\epsfig{file=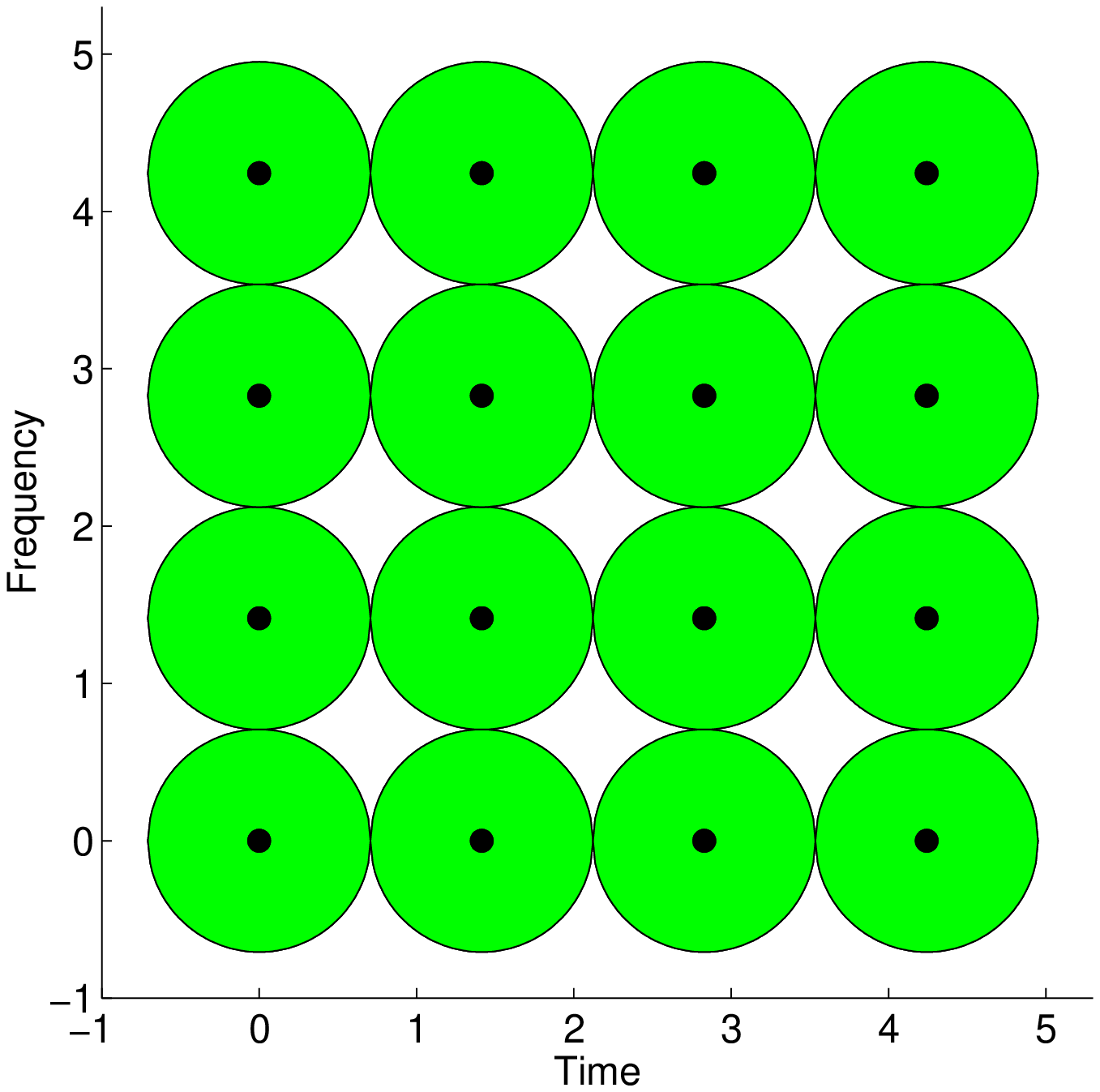,width=55mm}}
\caption{An OFDM system with rectangular lattice in the time-frequency plane. 
The data symbols $c_{k,\l}$ are transmitted at the lattice points $(kT,\l F)$. 
The spheres represent the effective support $\supp_{\eps} A \psikl$ of the 
ambiguity function $A \psikl$. The choice of $\eps$ is such that the 
difference to the LOFDM setup displayed in Figure \ref{fig:hex} becomes obvious.}
\label{fig:rect}
\end{center}
\end{figure}

The ambiguity function of $g$ is rotation-invariant, see \cite{Fol89}, therefore $|A g_{k,\l}|$ is rotation invariant. Hence for any $\eps >0,\,\supp_{\eps}(A g_{k,\l})$ is a sphere with center $(kT,\l F)$. It follows from \eqref{spect} that if the matrix $R$ is well-conditioned then $A\ph \approx Ag$. Hence $\supp_{\eps}(A \psikl)$ is approximately a sphere with center $(kT,\l F)$.  The OFDM system $(\ph,T,F)$ with $T=F$ and $\sigma=1$ can thus be represented
schematically by Figure \ref{fig:rect}. The data symbol $c_{k,\l}$ is transmitted at the lattice point $(kT,\l F)$.  The gray spheres represent $\supp_{\eps}A \psikl$. We have chosen $\eps$ such that the difference between the setup in Figure \ref{fig:rect} and Figure \ref{fig:hex}, see below, becomes obvious.  Certainly the spheres may overlap for smaller $\eps,$ but the $A \psikl$ are mutually orthogonal, since the $\psikl$ are mutually orthogonal. However after passing through the channel $\Hchan$ the orthogonality between the $\psikl$ (and thus between the $A \psikl$) is lost and the support of $A \psikl$ is spread out in time and frequency. And now the overlap between adjacent $A\psikl$'s (or rather between $A \Hchan \psikl$ and $A \Hchan \ph_{k',l'}$) plays a crucial role. Roughly speaking less overlap between the spheres results in less ISI and ICI. \\

Note that since $\psikl(t) = \ord(e^{- |t|})$ we have
\begin{equation}
|\langle \Hchan \ph_{k,\l},\ph_{k \pm 1, l \pm 1} \rangle| \gg
|\langle \Hchan \ph_{k,\l},\ph_{k \pm n, l \pm n}\rangle| , \quad \text{for}
\,\,n>1.
\label{nei}
\end{equation}
Therefore we can concentrate on the interference between $c_{k,\l}$ and direct 
neighbor symbols $c_{k \pm 1, l \pm 1}$, since the interference of
$c_{k,\l}$ with symbols $c_{k \pm n,l \pm n}$ where $ n>1$ will be negligibly 
small by comparison.  Hence the problem reduces to maximizing the distance 
between immediate neighbor lattice points $(k,\l)$ and $(k \pm 1, l \pm 1)$ 
without losing spectral efficiency. The answer to this problem
comes from the theory of sphere packings.\\

The sphere packing problem can be loosely described as follows.
We are given an infinite number of $d$-dimensional spheres of 
equal size. We identify the midpoint of each sphere with a point of
a $d$-dimensional lattice $\La$ in $\R^d$. 
In the sphere packing problem one tries to find the
lattice $\La_{opt}$ that solves
\begin{equation}
\underset{\La}{\max} \,\,\,
\left\{\frac{\text{Volume of a sphere}}{\text{Density of $\La$}}\right\}.
\label{sp1}
\end{equation}
It is easy to see that in the sphere packing problem the radius $r$ of a 
sphere is given by
\begin{equation}
\label{sp2}
r = \frac{1}{2} \Big( \underset{{\footnotesize \begin{matrix}\lambda_k,\lambda_l \in
\La,\\ \lambda_k \neq \lambda_l \end{matrix}}}{\min}
\big\{ \|\lambda_k - \lambda_l\| \big\}\Big).
\end{equation}
In words, the radius is equal to half of the minimal distance between
two lattice points.
We can reformulate \eqref{sp1} as follows:
\begin{equation}
\underset{\La}{\max}\big(\underset{{\footnotesize 
\begin{matrix}
\lambda_k,\lambda_l \in \La,\\ 
\lambda_k \neq \lambda_l 
\end{matrix}}}{\min} \big\{ \|\lambda_k - \lambda_l\| \big\}\big) 
\qquad \text{subject to $\xi(\La)= \xi_0$},
\label{sp3}
\end{equation}
for some arbitrary $\xi_0>0$.  Obviously we find ourselves confronted with the same problem as in OFDM, when we want to maximize the distance between adjacent lattice points in order to minimize ISI/ICI while keeping the spectral efficiency constant.\\

It is well-known that the sphere packing problem in two dimensions is solved by the hexagonal lattice, cf. \cite{CS93b}. Thus for the case $\tau_0 = \nu_0$ an OFDM system constructed from a hexagonal lattice as shown in Figure \ref{fig:hex} will provide optimal protection against ISI/ICI due to the increased distance between adjacent lattice points.  The effect of the increased distance becomes apparent when comparing Figure \ref{fig:rect} with Figure \ref{fig:hex}. In both cases we have chosen the same $\eps$ as the threshold for the effective support of the ambiguity function.\\

\begin{figure}[!ht]
\begin{center}
\subfigure[]{
\epsfig{file=scatcirc.eps,width=40mm}}
\subfigure[]{
\epsfig{file=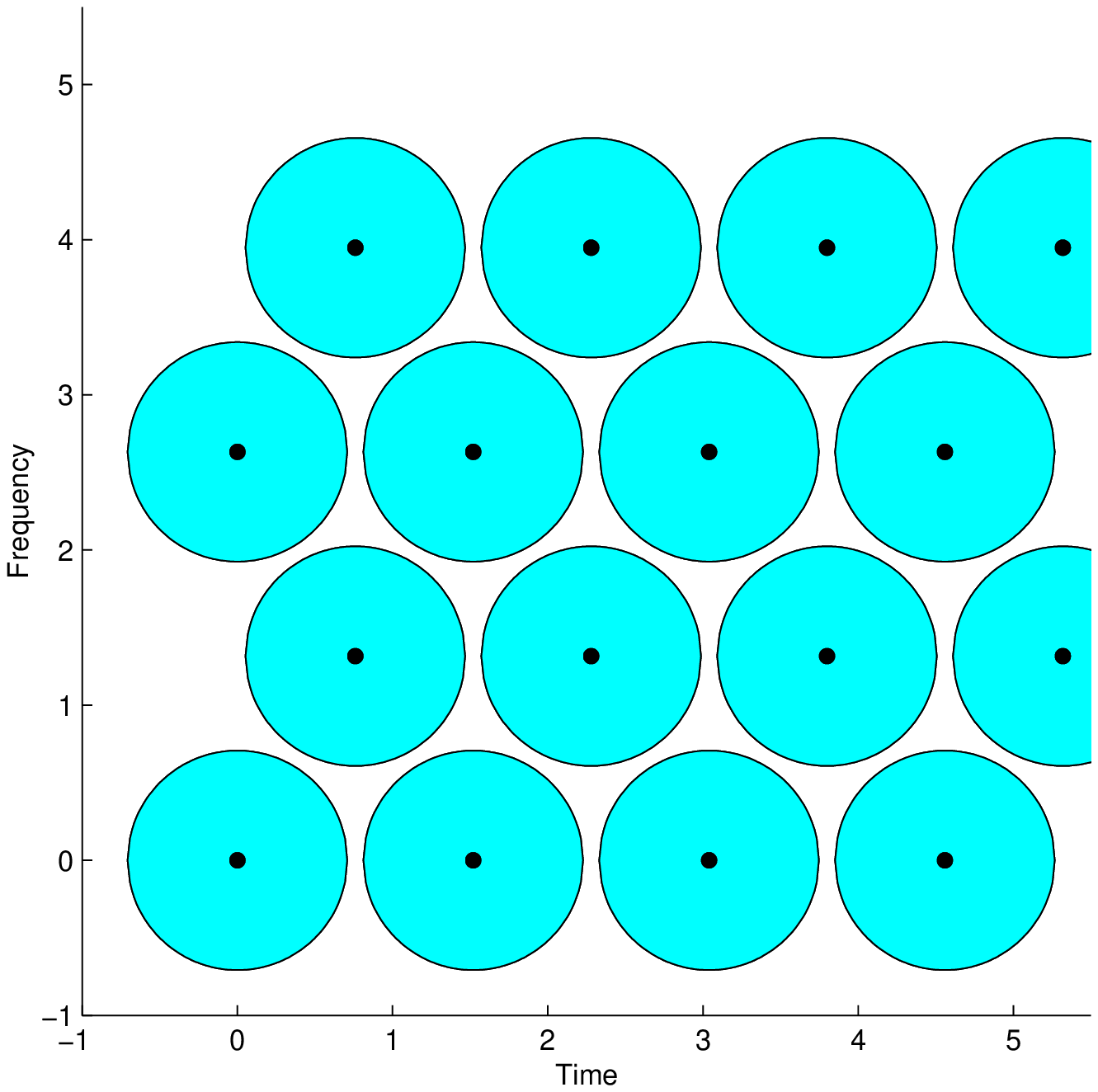,width=55mm}}
\caption{Optimal LOFDM system (Fig. \ref{fig:hex}(b)) for the channel 
scattering shown in Fig. \ref{fig:hex}(a). 
The data symbols $c_{k,\l}$ are transmitted at hexagonal lattice points.
The spheres represent $\supp_{\eps} A \psikl$ with the same
$\eps$ as in Figure \ref{fig:rect}. The use of the hexagonal lattice
clearly increases the distance between adjacent spheres as compared
to the rectangular lattice (cf.~Figure \ref{fig:rect}), thus leading
to reduced ISI/ICI for time-frequency dispersive channels.}
\label{fig:hex}
\end{center}
\end{figure}

Furthermore when we construct an orthogonal set of transmission functions
$(\ph,\La)$ from a non-orthogonal set of Gaussians $(\ggam,\La)$
via $\psigam = L\ddot{o}\,(\ggam,\La)$ the increased distance between adjacent 
lattice points yields a smaller correlation 
$|\langle \psigam_{\la'},\psigam_{\la}\rangle|$ between neighbor functions 
than in the case of a rectangular lattice. This results in a smaller
condition number of the associated Gram matrix $R_{\ggam,\La}$
for $R_{\ggam,T,F}$ at the same spectral efficiency, as can be seen
in Figure \ref{fig:condR}.\\

Of course the assumption $\tau_0= \nu_0$ in our example is very
restrictive. In the next subsection we consider more general channel 
conditions and derive a general LOFDM design rule in order to minimize 
the joint ISI/ICI.\\

We shall from this point forward restrict ourselves to $d=1$ since we are constrained to use only symplectic lattices.  In $\R^2$ any lattice is symplectic, but this is not the case in $\R^d$ with $d>2.$\\

\subsection{Optimal LOFDM design for the minimization of ISI/ICI} \label{ss:gen}
\begin{definition}
We define the dilation operator $\Dil$ by
\begin{equation}
\Dil f(t) = \sqrt{\alpha} f(\alpha t),\qquad \alpha \in \R_+,
\label{dil}
\end{equation}
and multiplication by a chirp is 
\begin{equation}
\Chirp f(t) = f(t) e^{-\pi i \beta t^2}, \qquad \beta \in \R.
\label{chirp}
\end{equation}
$\Dil$ and $\Chirp$ are unitary linear operators on $\LtR,$ and we have
${\cal D}_{\alpha}^{-1} = {\cal D}_{\alpha^{-1}}$ and 
${\cal C}_{\beta}^{-1} = {\cal C}_{-\beta}$.\\
\end{definition}

For $f,h \in \LtR$ we have \cite{Gro01}
\begin{align}
A(\Dil f,\Dil h)(t,\omega)& = A(f,h)(\alpha t, \alpha^{-1} \omega) 
\label{ambdil} \\
A(\Chirp f,\Chirp h)(t,\omega) &= A(f,h)(t, \omega +\beta t).
\label{ambchirp}
\end{align}

\begin{lemma}\label{fnd}
\be
\F\D_{\alpha^{-1}}=\D_\a\F
\end{equation}
\end{lemma}
\begin{myproof}\bea
{\cal F}{\cal D}_{\alpha^{-1}}f(t)&=&\frac{1}{\sqrt{\a}}{\cal F}f(t/\a)\nn\\
&=&\frac{1}{\sqrt{\a}}\int_{\R^d}f(t/\a)e^{-2\pi i\w t}\,dt\nn\\
&=&\sqrt{\a}\int_{\R^d}f(t)e^{-2\pi i\w\a t}\,dt\nn\\
&=&\D_\a\F f(t).\nn
\end{eqnarray}\end{myproof}

Recall that $\tau_0$ and $\nu_0$ describe the support of the channel scattering function.
More precisely, assuming that the support of the scattering function is circumscribed
by a rectangle (or ellipse), then $\tau_0$ gives the length of the rectangle (in the time-direction) and $\nu_0$ corresponds to the width (in the frequency direction).  In the following theorem we extend the result from Section \ref{ss:sphere} to more general (and more realistic) choices for $\tau_0$ and $\nu_0$.  The theorem is essentially only a special case of a deep theorem by Gr\"ochenig \cite{Gro01}, but it shows how we can switch between different lattices by adapting the pulse shape accordingly in order to preserve the properties of the optimal hexagonal case.  In fact our approach will also include the case where the support of the scattering function may be concentrated along a rotated rectangle in the time-frequency domain (like in the case of a chirp). Such a setup has also been investigated in \cite{BT01}. \\

\begin{theorem}\label{th:sn}
 Let $(f_1,\La_1)$ be a non-orthogonal linearly independent system with $\|R_{f_1,\La_1}\|<\infty,$ and denote $\varphi_1=L\ddot{o}\,(f_1,\La_1).$  Assume $\La_2$ is a lattice with density $\xi(\La_2)=\xi(\La_1)$ and let
\be
L_1=\begin{bmatrix}x_1&y_1\\0&z_1\end{bmatrix}\quad\mbox{ and }\quad L_2=\begin{bmatrix}x_2&y_2\\0&z_2\end{bmatrix}
\end{equation}
be the generator matrices of $\La_1$ and $\La_2,$ respectively.  If
\be
\label{f2adapt}
f_2={\mathcal F}{\mathcal D}_{\alpha^{-1}}{\mathcal C}_{-\beta\alpha^2}{\mathcal F}^{-1}f_1
\end{equation}
with
\be
\alpha=\frac{z_2}{z_1},\quad\beta=-\frac{y_2}{z_2} + \frac{y_1}{z_1}\frac{z_1^2}{z_2^2},
\end{equation}
then
\be
\label{eqgram}
R_{f_2,\La_2}=R_{f_1,\La_1}
\end{equation}
and if we denote $\varphi_2:=L\ddot{o}\,(f_2,\La_2),$ then 
\be
\label{orthadapt}
\varphi_2={\mathcal F}{\mathcal D}_{\alpha^{-1}}{\mathcal C}_{-\beta\alpha^2}{\mathcal F}^{-1}\varphi_1
\end{equation}
\end{theorem}
\begin{myproof}
The proof consists of applying Proposition 9.4.4 in \cite{Gro01} twice.  Fix an arbitrary rectangular lattice $\La_R$ with generator matrix
\be
L_R=\begin{bmatrix}T&0\\0&F\end{bmatrix}
\end{equation}
such that $\xi(\La_R)=\xi(\La_1).$  It is known that any ${\mathcal A}\,\in\,$ Sp$(d)$ can be written as a finite product of matrices of the form $J, C_\beta,$ and $D_\alpha.$  It is thus clear with some calculation that $L_1$ and $L_2$ can be expressed as
\bea
L_1&=&JC_{\beta_1}D_{\alpha_1}J^{-1}L_R,\qquad \mbox{ with }\alpha_1=\frac{z_1}{F},\,\beta_1=-\frac{y_1}{z_1}\\
L_2&=&JC_{\beta_2}D_{\alpha_2}J^{-1}L_R,\qquad \mbox{ with }\alpha_2=\frac{z_2}{F},\,\beta_2=-\frac{y_2}{z_2}
\end{eqnarray}
where $J,C_\beta,D_\alpha$ are defined by
\be
\label{JCD}
J=\bbm0&1\\-1&0\end{bmatrix},\quad C_\beta=\bbm1&0\\\beta&1\end{bmatrix},\quad D_\alpha=\bbm\alpha&0\\0&\alpha^{-1}\end{bmatrix},\quad\alpha,\beta\,\in\,\R.
\end{equation}
Let Sp$(1,\R)$ denote the symplectic group\footnote{Any lattice in $\R^2$ has a generator matrix in Sp$(1,\R)$.}, and let $\mu(\cdot)$ denote the symplectic (unitary) operator appearing in the metaplectic representation of Sp$(1,\R)$ (see \eqref{unitarySVN} in Appendix \ref{repSvNHg}). Let $J,C_\beta,D_\alpha$ be defined as in \eqref{JCD}, then
\be
\mu(J)={\mathcal F}^{-1}\quad\mu(C_\beta)={\mathcal C}_{-\beta}\quad\mu(D_\alpha)={\mathcal D}_{\alpha^{-1}}\label{unops}
\end{equation}
Furthermore there holds (despite the fact that $\mu$ is not a homomorphism in general)
\be
\mu(JC_\beta D_\alpha J^{-1})=\mu(J)\mu(C_\beta)\mu(D_\alpha)\mu(J^{-1}).\label{unophomo}
\end{equation}
Define
\be
A_1:=JC_{\beta_1}D_{\alpha_1}J^{-1},\quad A_2:=JC_{\beta_2}D_{\alpha_2}J^{-1}.
\end{equation}
Then
\bea
A_2A_1^{-1}&=&JC_{\b_2}D_{\a_2}J^{-1}JD_{\a_1^{-1}}C_{-\b_1}J^{-1}\nn\\
&=&JC_{\b_2}D_{\a_2}D_{\a_1^{-1}}C_{-\b_1}J^{-1}\nn\\
&=&JC_{\b_2}D_{\frac{\a_2}{\a_1}}C_{-\b_1}J^{-1}\nn\\
&=&JD_{\frac{\a_2}{\a_1}}C_{\left[\b_2\big(\frac{\a_2}{\a_1}\big)^2-\b_1\right]}J^{-1}.\label{a2a1inv}
\end{eqnarray}
where in the last step we used the fact that $C_\b D_\a=D_\a C_{\b\a^2}.$ From \eqref{a2a1inv}, we conclude 
\bea
\m (A_2A_1^{-1})&=&\m (JD_{\frac{\a_2}{\a_1}}C_{\left[\b_2\big(\frac{\a_2}{\a_1}\big)^2-\b_1\right]}J^{-1})\nn\\
&=&\m(J)\m(D_{\frac{\a_2}{\a_1}})\m(C_{\left[\b_2\big(\frac{\a_2}{\a_1}\big)^2-\b_1\right]})\m(J^{-1})\nn\\
&=&{\mathcal{F}}^{-1}\,\,{\mathcal{D}}_{\frac{\a_1}{\a_2}}\,\,{\mathcal{C}}_{\left[\b_1-\b_2\big(\frac{\a_2}{\a_1}\big)^2\right]}\,\,{\mathcal{F}}\nonumber\\
&=&{\mathcal{F}}^{-1}\,\,{\mathcal{D}}_{\a^{-1}}\,\,{\mathcal{C}}_{-\b\a^2}\,\, {\mathcal{F}}.\label{a2a1invop}
\end{eqnarray}

Set $f_R=\mu(A_1)^{-1}f_1$ and $f_2=\mu(A_2)^{-1}f_R.$  From Proposition 9.4.4 in \cite{Gro01} we have that
\be
(f_R,\La_R)=\mu(A_1)^{-1}(f_1,\La_1),
\end{equation}
and
\be
(f_2,\La_2)=\mu(A_2)(f_R,\La_R).
\end{equation}
Since $\mu(A_1)$ and $\mu(A_2)$ are unitary, we get
\be
R_{f_1,\La_1}=R_{f_R,\La_R}=R_{f_2,\La_2}
\end{equation}
where $L_2=A_2A_1^{-1}L_1.$  Equation \eqref{orthadapt} now follows from \eqref{a2a1invop}, \eqref{f2adapt} and \eqref{eqgram}.
\end{myproof}

It is possible to prove Theorem \ref{th:sn} directly using signal analytic
methods. However we feel that the metaplectic representation is simpler
and more elegant. The reader familiar with the mathematical theory
of quantum physics will notice that Theorem \ref{th:sn} is essentially
a consequence of the celebrated Stone-von Neumann Theorem on
unitary representations of the Heisenberg group (see Appendix \ref{repSvNHg}).\\

The main feature of Theorem \ref{th:sn} is that it yields a simple
rule for the design of optimal LOFDM systems for time-frequency
dispersive channels. We illustrate this by several examples.

\subsection{Some examples} \label{ss:ex}
\begin{example}
We compute $\varphi_h,$ the pulse generating a tight frame based on a hexagonal lattice with twofold oversampling, from $\varphi_r,$ the pulse generating a tight frame based on a rectangular lattice with twofold oversampling (these pulses could then be used as generators of OFDM systems on their respective adjoint lattices).  If
\be
L_r=\begin{bmatrix}\frac{\sqrt{2}}{2}&0\\0&\frac{\sqrt{2}}{2}\end{bmatrix},\quad\mbox{and}\quad L_h=\frac{1}{\sqrt[4]{3}}\begin{bmatrix}1&\frac{1}{2}\\0&\frac{\sqrt{3}}{2}\end{bmatrix},\nn
\end{equation}
then $\alpha=\frac{\sqrt[4]{3}}{2}\Big/\frac{\sqrt{2}}{2}=\frac{\sqrt[4]{3}}{\sqrt{2}}$ and  $\beta=-\frac{1}{\sqrt{3}},$ and we simply use \eqref{unops}, \eqref{unophomo}, and \eqref{orthadapt} to get
\be
\ph_h={\mathcal F}{\mathcal D}_{\frac{\sqrt{2}}{\sqrt[4]{3}}}{\mathcal C}_{\frac{1}{2}}{\mathcal F}^{-1} \ph_r.
\end{equation}
\end{example}
In Section \ref{s:pulse} we mentioned that the TFL of the pulse $\ph=L\ddot{o}\,(g,\sqrt{2},\sqrt{2})$ measured in terms of second-order moments is $\Delta \tau_{\ph}\cdot\Delta \nu_{\ph}=1.024\cdot\frac{1}{4\pi},$ which is conjectured in \cite{FAB95} to be the optimum among all pulses creating an OFDM system (or OQAM-OFDM system).  This conjecture is plausible for rectangular lattices, but is untrue if we allow more general lattices. For instance if we use the hexagonal lattice $\Lambda_H$ as defined in \eqref{hexmat} with $T=F=\sqrt{2}$ and set $\ph_H = L\ddot{o}\,\,(g,\Lambda_H)$ we compute numerically the TFL $\Delta \tau_{\ph_H} \cdot \Delta \nu_{\ph_H} = 1.014\cdot\frac{1}{4\pi}$.
\begin{conjecture}
$\ph_H$ yields the optimal TFL among all possible lattices and all possible pulse shapes that create an LOFDM system with spectral efficiency $\r=1/2$ (and similarly for other values of $\r$).
\end{conjecture}

In the following examples we consider LOFDM systems with spectral efficiency
$\rho$.  For the remainder of this section $\La_1$ will denote the hexagonal lattice 
$\La_h$ with
generator matrix
\begin{equation}L_h = L_1 :=
\frac{1}{\sqrt{\r}}\begin{bmatrix}
\frac{\sqrt{2}}{\sqrt[4]{3}} & \frac{\sqrt{2}}{2\sqrt[4]{3}} \\
0 & \frac{\sqrt[4]{3}}{\sqrt{2}}\label{l1}
\end{bmatrix}
\end{equation}
and 
\begin{equation}
\label{l2}
\ph_1:=L\ddot{o}\,(g,\La_h).
\end{equation}
We have seen in Subsection \ref{ss:sphere} that the LOFDM system
$(\ph_1, \La_h)$ is optimally robust against ISI/ICI if
$\tau_0=\nu_0$. Based on Theorem \ref{th:sn} we are able now to
consider more general situations.\\

\begin{example} Fix $\r$ and assume $\tau_0=\frac{1}{\sigma}\nu_0$ for some $\sigma >0.$
The LOFDM system that is optimally adapted to these channel conditions
with respect to minimizing the joint ISI/ICI is given by
$(\psigam, \La_{\sigma})$, where the generator matrix $L_{\sigma}$ of 
$\La_{\sigma}$ can be computed by left-multiplication of $L_1$ by the dilation matrix $D_{1/\sqrt{\s}}$ as \footnote{see \eqref{orthdil} and \eqref{unops}; also note that since it is the {\em column vectors} of the generator matrices which are the generating vectors, left-multiplication is the correct operation} 
\begin{equation}
\label{sigmat}
L_{\sigma}=
\frac{1}{\sqrt{\r}}\begin{bmatrix}
T & T/2 \\
0 & F
\end{bmatrix},
\quad \text{with}\,\,\, T=\frac{\sqrt{2}}{\sqrt[4]{3}\sqrt{\s}},\,\,F=\frac{1}{T}
\end{equation}
and the optimal pulse shape is $\psigam=L\ddot{o}\,(\ggam,\La_{\sigma}).$ 

Consider for instance $\tau_0 = \sqrt{3} \nu_0$ and $\rho=\frac{1}{2}$. By 
using \eqref{sigmat} we obtain 
\begin{equation}
L_2 = 
\begin{bmatrix}
2 & 1\\
0 & 1 
\end{bmatrix}
\label{L2}
\end{equation}
as generator matrix of the optimal lattice $\La_2$. 
Using \eqref{eqgram} and \eqref{f2adapt} with $L_1$ as in \eqref{l1} and
$\ph_1$ as in \eqref{l2} 
we compute $\alpha=\frac{1}{\sqrt[4]{3}}$, $\beta=0$, hence
\begin{equation}
\ph_{\sqrt{3}}(t)={\cal F}{\cal D}_{\alpha^{-1}}{\cal C}_{-\beta\alpha^2}{\cal F}^{-1}
\ph_1(t)={\cal F}{\cal D}_{\sqrt[4]{3}}{\cal F}^{-1}\ph_1(t) = {\cal D}_{\frac{1}{\sqrt[4]{3}}}\ph_1(t),
\label{ex1}
\end{equation}
where we have used Lemma \ref{fnd}.  Thus the optimal LOFDM setup for the given conditions
$\rho=\frac{1}{2}, \tau_0= \sqrt{3}\nu_0$ can be schematically represented as
shown in Figure \ref{fig:fig5}. 

\begin{figure}[!ht]
\begin{center}
\subfigure[]{
\epsfig{file=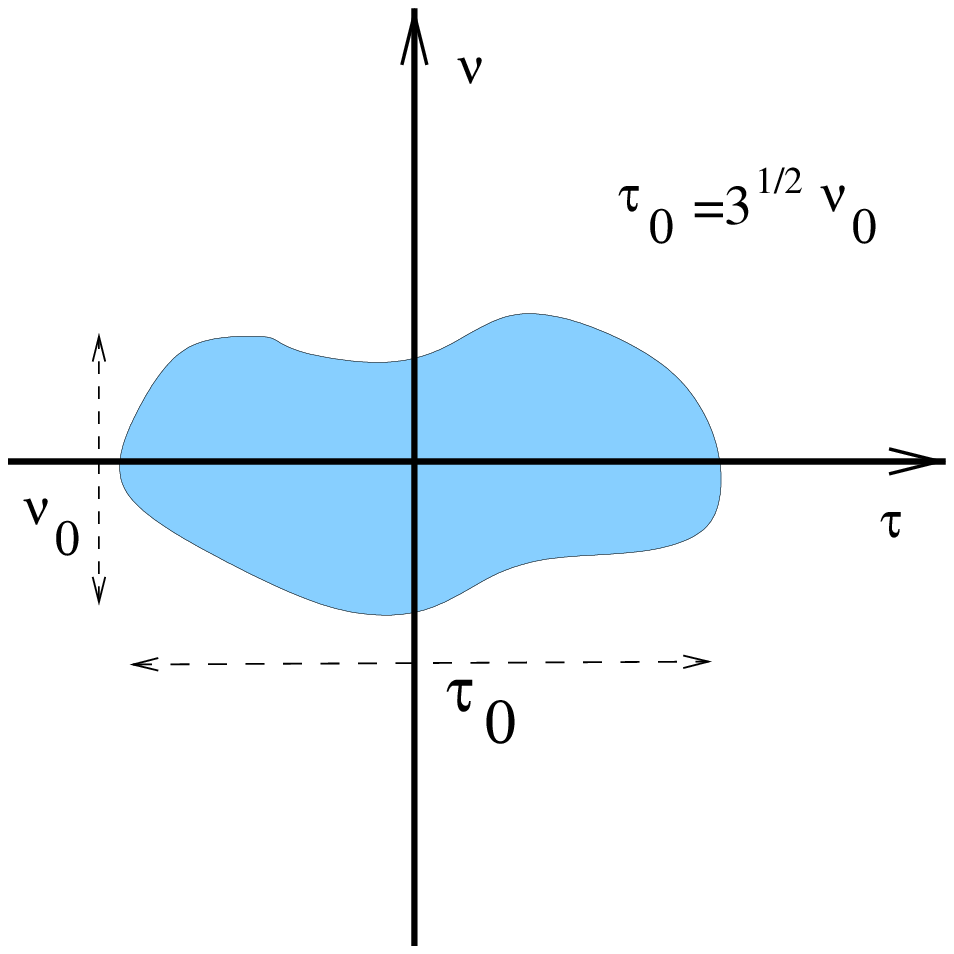,width=40mm}}
\subfigure[]{
\epsfig{file=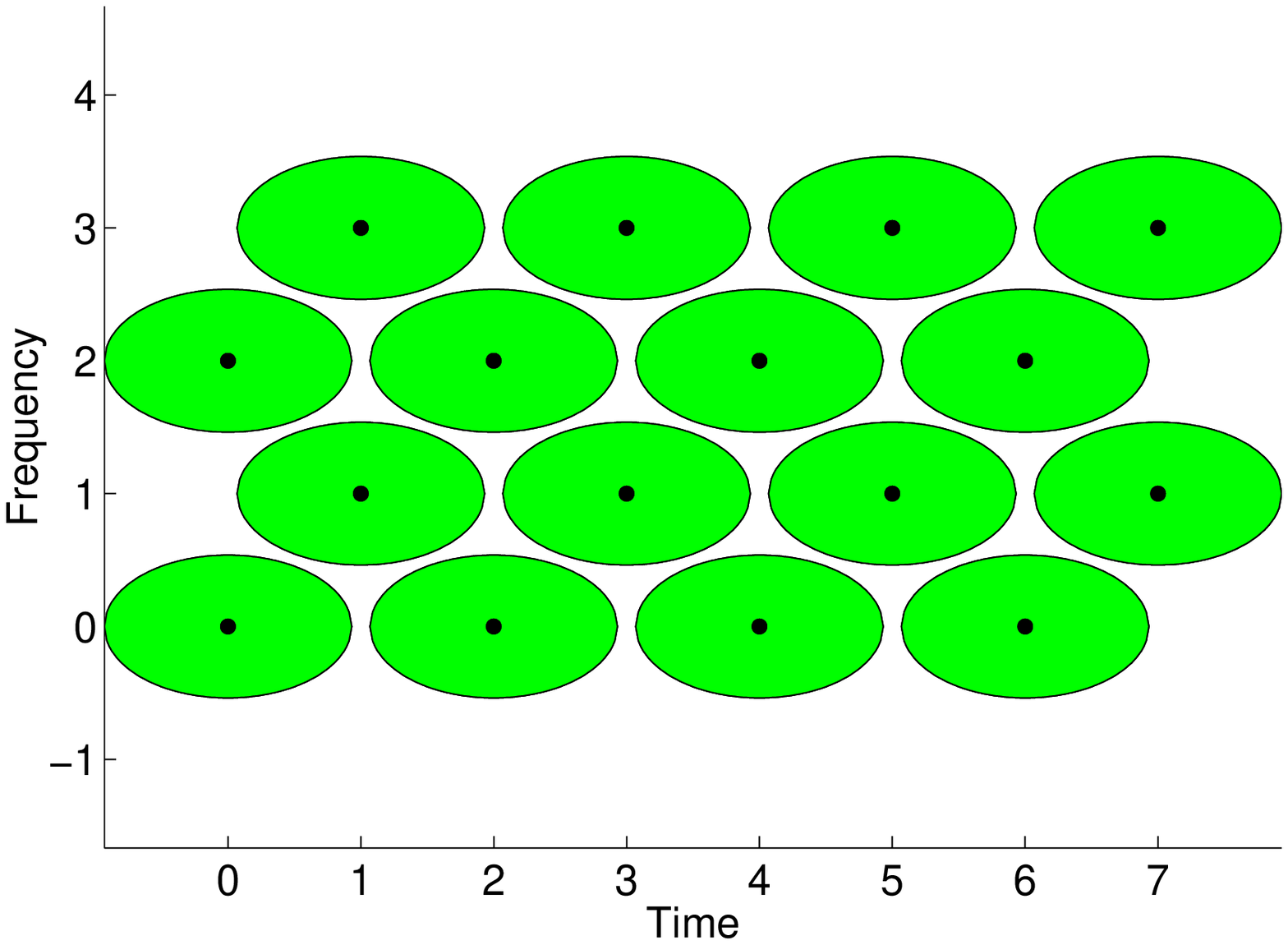,width=55mm}}
\caption{Optimal LOFDM setup (Fig. \ref{fig:fig5}(b)) for a wireless channel 
with $\tau_0 = \sqrt{3} \nu_0$ (scattering function shown in 
Fig. \ref{fig:fig5}(a)). The ellipses represent 
$\supp_{\eps} A (\ph_{\sqrt{3}})_{k,\l}$, where we have 
chosen the same $\eps$ as in Figure \ref{fig:hex}.}
\label{fig:fig5}
\end{center}
\end{figure}

\end{example}

\begin{example} Assume $\rho=\frac{1}{2}$ and $\tau_0=\sqrt{3}\nu_0$, but now let the scattering 
function be rotated clockwise by 45 degrees in the time-frequency plane, 
as shown in Figure \ref{fig:fig6}(a). In this case Theorem \ref{th:sn}
together with the adaptation rule \eqref{adapt} yields 
\begin{equation}
L_2'=Rot_r(-45)\cdot L_1 =\begin{bmatrix}
\frac{\sqrt{2}}{2} & \frac{\sqrt{2}}{2} \\
-\frac{\sqrt{2}}{2} & \frac{\sqrt{2}}{2}
\end{bmatrix}
\begin{bmatrix}
2 & 1 \\
0 & 1
\end{bmatrix}=\begin{bmatrix}
\sqrt{2} & \sqrt{2} \\
-\sqrt{2} & 0
\end{bmatrix}
\label{ex2}
\end{equation}
where $Rot_r(-45)$ is the left-multiplication Givens rotator.  Now $L_2'$ generates the same lattice as 
\be  
L_2=U\cdot L_2'=\begin{bmatrix}
0 & -1 \\
1 & 1
\end{bmatrix}
\begin{bmatrix}
\sqrt{2} & \sqrt{2} \\
-\sqrt{2} & 0
\end{bmatrix}
=\begin{bmatrix}
\sqrt{2} & 0 \\
0 & \sqrt{2}
\end{bmatrix}
\end{equation}
since $U$ is unitary and leaves one lattice-generating vector unchanged, and re-forms the second as a linear combination of the lattice-generating vectors in $L_2'.$  Thus we can take $L_2$ as the generator matrix of the optimal lattice $\La_2$, which is a rectangular lattice in this case.
Proceeding in the same way as in the previous example 
we calculate $\alpha=\frac{\sqrt{2}}{\sqrt[4]{3}}$, $\beta=\frac{1}{2}$, and
\be
\ph'_{\sqrt{3}}(t)= {\cal F}{\cal D}_{\alpha^{-1}}{\cal C}_{-\beta \alpha^2}
{\cal F}^{-1} \ph_1(t)=
{\cal F}{\cal D}_{\frac{\sqrt[4]{3}}{\sqrt{2}}}{\cal C}_{-\frac{1}{\sqrt{3}}}\ph_1(t).
\end{equation}
Since $\ph'_{\sqrt{3}}$ is the Fourier transform of a dilated chirped Gaussian, formulas
\eqref{ambdil}, \eqref{ambchirp}, and \eqref{spect} then show that $\supp_{\eps}A\ph'_{\sqrt{3}}$
is (approximately) an ellipse, with ratio of semi-major vs semi-minor axis equal to $\sqrt{3}$, 
rotated clockwise by 45 degrees. Hence the corresponding LOFDM system
$(\ph'_{\sqrt{3}},\La_R)$ can be schematically
represented as in Figure \ref{fig:fig6}.

\begin{figure}[!ht]
\begin{center}
\subfigure[]{
\epsfig{file=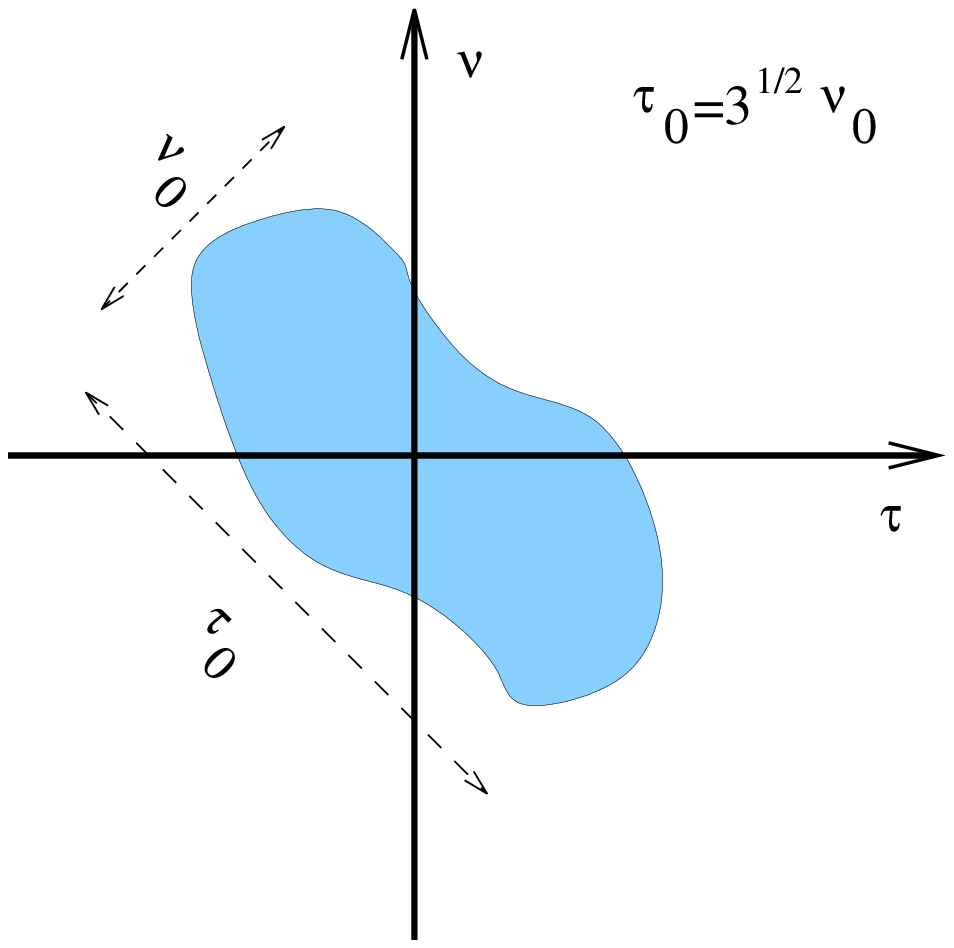,width=40mm}}
\subfigure[]{
\epsfig{file=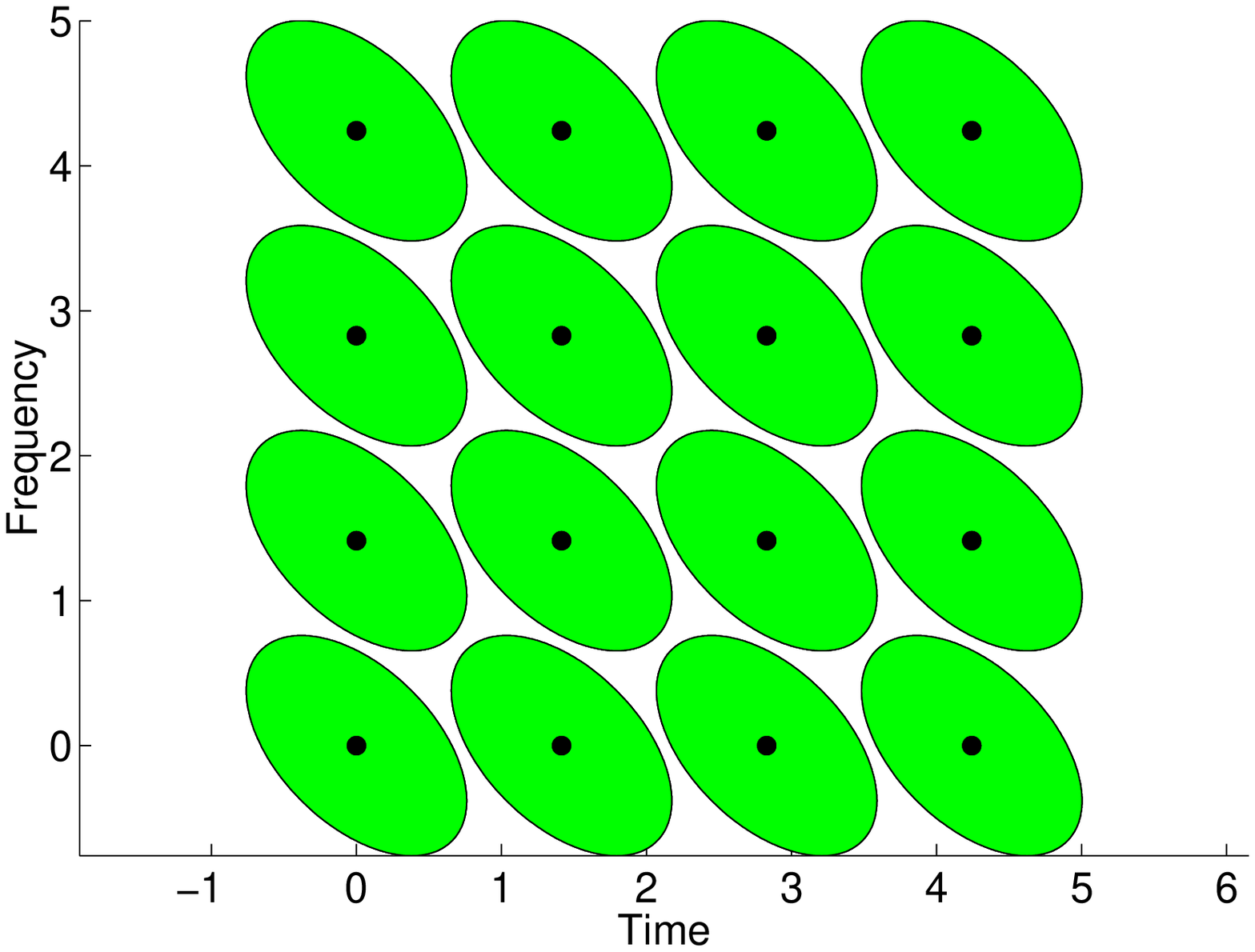,width=55mm}}
\caption{Optimal LOFDM setup (Fig. \ref{fig:fig6}(b)) for a wireless channel 
with 45-degree rotated scattering function satisfying 
$\tau_0 = \sqrt{3} \nu_0$ (see Fig. \ref{fig:fig6}(a)). The ellipses 
represent $\supp_{\eps} A (\psi_2)_{k,\l}$, where 
we have chosen the same $\eps$ as in Figure \ref{fig:hex}.}
\label{fig:fig6}
\end{center}
\end{figure}

\end{example}

\remark The pulse-shaping framework derived in Section \ref{s:pulse} also
holds for general lattices, since Theorem \ref{th:min} 
carries over to general lattices (see Appendix \ref{gensympl}). The proof of the extension of 
Theorem \ref{th:min} can be done by using Theorem \ref{th:sn} and the
results derived in Chapter~9.4 in \cite{Gro01}. \\

\if
Note that Theorem \ref{th:min} does not necessarily mean that 
$A \psigam \approx A \ggam$. The latter is only true if the condition number of
the Gram matrix $R_{\ggam,\La}$ is small. Although $A\ggam$ has always the shape of an ellipsoid (the length and width of which depend only on $\sigma$) this is not the case for $\psigam$, since
$\psigam$ depends on $\ggam$ {\em and} on $\La$, which becomes apparent in Theorem \ref{th:lim}.  When we change $\sigma$, but leave the lattice $\La$ unchanged, the condition number of $R_{\ggam,\La}$ will increase (even if $\det \La \\ll 1$). This follows for rectangular lattices from \cite{TO95} and for general lattices by combining the results in \cite{TO95} with Theorem \ref{th:sn}.  Indeed, as a consequence of Theorem \ref{th:lim} the shape of $\supp_{\eps} A \psigam$ changes from sphere-shaped to rectangular when we let $\sigma \rightarrow 0$ or $\sigma \rightarrow \infty$.  \\

We want to determine the optimal pulse $\psi$ and the optimal lattice $\La$ for given channel conditions and predefined spectral efficiency $\rho$, such that we maximize the robustness of the transmission system with respect to ISI/ICI. We assume that the maximal delay spread $\tau_0$ and the maximal Doppler spread $\nu_0$ are known. We assume that the channel properties are constant during one coherence interval. If we have a channel with feedback -
such as an TDD OFDM system - then we can adapt the setup of the OFDM system to varying channel conditions. All that is required is knowledge about the ratio $\tau_0/\nu_0$. The proposed (L)OFDM adaptation procedure can be divided into the following four steps.

1. Set $T=\frac{1}{\sqrt{\rho}}\sqrt{\tau_0}{\nu_0}$,
$F=\frac{1}{\sqrt{\rho}}\sqrt{\nu_0}{\tau_0}$, and
let $\La_o$ be the lattice with generator matrix
$$L_= \begin{bmatrix} \alpha \frac{\sqrt{2}}{\sqrt[4]{3}}T & \beta \frac{\sqrt{2}}{2\sqrt[4]{3}}T 
\\ 0 & \frac{1}{\alpha} \frac{\sqrt[4]{3}}{\sqrt{2}} F
\end{bmatrix},$$
thus the generator matrix is normalized such that $\det L = 1/\rho$ and 
the lattice $\La_o$ yields the desired spectral efficiency.

2. Set $\sigma\approx \tau_0/\nu_0$ and and compute $\ggam = (2\sigma)^{\frac{1}{4}} e^{-\pi \sigma t^2}$.

3. Define $\alpha = \frac{\sqrt[4]{3}}{\sqrt{2}}\sigma$ and $\beta=-1/2$.

4. Compute the generator matrix $L_o$ for lattice $\La_o$ and calculate $\psi=\text{orth}(\ggam,\La_o)$. Then the LOFDM system given by $(\psi, \La_o)$ is ``essentially optimal'' with respect to minimizing ISI/ICI.\\
\fi

It is interesting to compare the condition number of $R_R:=R_{g,\Lambda_R}$ with that of $R_H:=R_{g,\Lambda_H}$ for $T=F=\sqrt{2}$, where the generator matrices for these lattices are as defined in \eqref{rectmat} and \eqref{hexmat}. We get $\cond (R_R)=\sqrt{2}$ and $\cond(R_H)\approx 1.2599$ which is suspiciously close to $\sqrt[3]{2}$.
\begin{conjecture}
$R_{g,\Lambda}=\sqrt[3]{2}$ is the minimal condition number achievable for any lattice $\La$ with $\xi(\La)=1/2.$
\end{conjecture}

If $\nu_0 \approx 0$ then OFDM systems employing (almost) rectangular pulses such as cyclic prefix OFDM are of course optimal.  The essential support of the ambiguity function of a rectangular pulse has rectangular shape.  It is clear that the optimal lattice for such pulses is a rectangular lattice.  Thus the introduction of general lattices would only decrease the performance in such a situation.

\section{Efficiency considerations in the numerical implementation of LOFDM} \label{ss:eff}
Before we proceed to a performance comparison between OFDM and LOFDM we briefly discuss the numerical implementation of LOFDM. We do not go into detail here, but note that the results in this section can be readily derived by using the link between OFDM and Gabor analysis and the efficient algorithms that have been developed in association with suitable finite-dimensional models for Gabor analysis \cite{Str97a}.\\

The requirement for an efficient implementation places some limits on the possible values of the 
entries $x,y,z$ of the generator matrix of the lattice (see Equation \ref{genmat}). One can easily see that an FFT-based implementation is only possible if $x,y,z\in \Q$. In particular if $x/y \in \N$ then the lattice can be written as sum of $x/y$ rectangular lattices, which differ from each other only by the time-frequency offset. For each of these rectangular sub-lattices we can use the usual FFT-based implementation for pulse-shaping OFDM systems. An example of a ``good'' generator matrix is
\be
L=\begin{bmatrix}
T & T/2 \\
0 & F\nn
\end{bmatrix}
\end{equation}
with $TF \in \Q$. Specifically, if $TF = p/q$ with $p,q \in \N$, the length of the FFT is given by $s^2 pq$, where $s$ is a constant that depends on the ratio of the sampling rate and the decimation factor in the A/D conversion and is often of the form $s=2^n$. 

\section{Comparison of performance of LOFDM vs OFDM} \label{ss:comp}
We compare the proposed LOFDM system to a standard pulse-shaping OFDM system of \eqref{ofdmsystem} associated with the rectangular lattice.  For both systems we use the pulse obtained from the Gaussian via $L\ddot{o}\,(\ggam,\La_R)$ and $L\ddot{o}\,(\ggam,\Lambda_H)$ respectively. For the rectangular lattice the OFDM system obtained in that way  coincides with the setup suggested in \cite{FAB95}. In both cases we optimally adapt the pulse and the lattice parameters to the given delay spread and Doppler spread.\\ 

All parameters are as described in Section \ref{numimplofdm}, except that AWGN on the order of 0.2\% of the signal energy is added during the simulated transmission (the noise level is 9).  We choose this small percentage of noise in our comparison in order to concentrate on the effects caused by ISI and ICI. \\

For both OFDM and LOFDM we compute the SNIR of the received data $d_{k,\l},$ which have been equalized according to \eqref{equal}.  Since both OFDM and LOFDM are orthogonal systems they are equally well (in fact optimally) adapted to AWGN. Hence any statistical difference in the SNIR can be attributed to the difference in robustness against ISI/ICI.\\

We fix the product $\tau_0 \cdot \nu_0$, randomly vary the parameters of the scattering for each experiment and compute the average SNIR for OFDM and LOFDM over 100 such experiments. We repeat this simulation for a variety of values $\tau_0 \cdot \nu_0$ in the range of $10^{-4} - 10^{-2}$ (but we always adapt both systems approximately optimal to the new values of $\tau_0$ and $\nu_0$). The result is shown in Figure \ref{fig:fig3}. As expected from the theoretical analysis LOFDM yields indeed a better SNIR than OFDM. On the average the improvement of LOFDM compared to OFDM is in the range of 1dB. 

\begin{figure}[htb]
\begin{center}
\epsfig{file=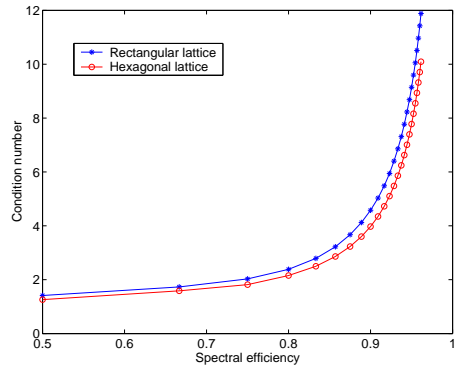,width=60mm}
\caption{Comparison of the condition number of the Gram matrix of the
OFDM-type function system $(\ggam,T,F)$ ($\ggam$ is the Gaussian) 
for the rectangular lattice and the hexagonal 
lattice when the product $(TF)^{-1}$ approaches 1. In both cases the 
condition number is small for $TF \in [\frac{1}{2},\frac{4}{5}]$.
As expected $\cond(R)$ for the hexagonal lattice is definitely smaller
than $\cond(R)$ for the rectangular lattice, thus demonstrating the
potential advantage of designing OFDM systems for general time-frequency
lattices.}
\label{fig:condR}
\end{center}
\end{figure}

\begin{figure}[htb]
\begin{center}
\epsfig{file=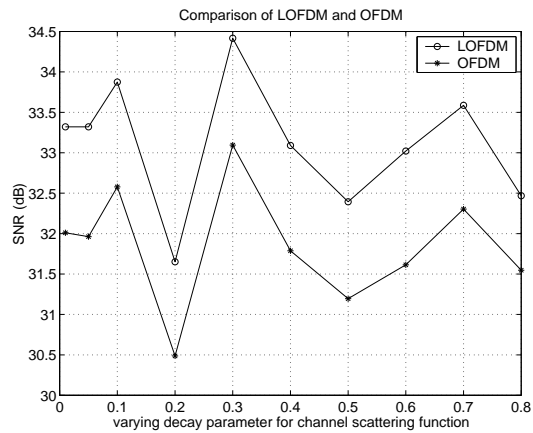,width=70mm}
\caption{Comparison of SNIR for OFDM and LOFDM for $\tau_0 \cdot \nu_0 \in [10^{-4}, 10^{-2}]$.
The proposed LOFDM scheme yields an improvement of about 1dB.}
\label{fig:fig3}
\end{center}
\end{figure}

\np
\pagestyle{myheadings}
\chapter{Future Research}
\thispagestyle{myheadings} 
\markright{  \rm \normalsize CHAPTER 8. \hspace{0.5cm}
 FUTURE RESEARCH}

\section{Weighted (under matrix and Hadamard multiplication) L\"owdin orthogonalization via the SVD}

Let $\Ph\in\C^{m\times n}$ with $m\geq n,$ with rank($\Ph$) $=n.$  We consider the general case of possibly nonsquare, rank-deficient weight matrix $W.$  Let $W\in\C^{n\times k},\,k\geq n$ with rank($W$) $=r\leq n.$  The (reduced) singular value decomposition, or SVD, of $\Ph W$ is 
\bea
\Ph W&=&UDV^*\nonumber\\
&=&[U_1,U_2]\biggl[\begin{tabular}{ll} $D_1$&$0_{k\times(k-r)}$\\$0_{(k-r)\times r}$&$0_{(k-r)\times(k-r)}$\end{tabular}\biggr][V_1,V_2]^*\nonumber\\
&=&U_1D_1V^*_1 
\end{eqnarray}
where $U\in\C^{m\times k},\,V\in\C^{k\times k},\,D\in\R^{k\times k};\,D_1$ is the diagonal $r\times r$ matrix containing the $r$ nonzero singular values of $\Ph W,\,U_1\in\C^{m\times r},$ and $V_1\in\C^{k\times r}.$  The SVD of $W$ is $W=U_WD_WV^*_W=U_{W1}D_{W1}V^*_{W1},$ where $U_W\in\C^{n\times k},D_W\in\R^{k\times k},V_W\in\C^{k\times k},$ and $U_{W1}\in\C^{n\times r},\,D_{W1}\in\R^{r\times r},$ and $V_{W1}\in\C^{k\times r}.$ \\
\begin{lemma}\label{wtdlowdinmat}
Among all unitary matrices $Q,\,\|(\Ph-Q)W\|_F$ is minimized by
\be
Q_{opt}=P:=UV^*,
\end{equation}
where $U$ and $V$ are unitary matrices in $USV^*,$ the reduced, compressed SVD of $\Ph WW^*$ ($U\in\C^{n\times r},\,S\in\R^{r\times r},$ and $V\in\C^{n\times r}).$
\end{lemma}
\begin{myproof}
We have
\bea
\|\Ph W-QW\|_F^2&=&\langle \Ph W-QW, \Ph W-QW\rangle\nonumber\\
&=&\|\Ph W\|_F^2+\|QW\|_F^2-\langle \Ph W,QW\rangle-\langle QW,\Ph W\rangle\nonumber\\
&=&\|\Ph W\|_F^2+\|W\|_F^2-2{\mbox Re}\langle \Ph W,QW\rangle.
\end{eqnarray}
Denote the singular values in $S$ by $\s_i,$ and $V^*Q^*U$ by $\{t_{ij}\}.$  Clearly we want to maximize ${\mbox Re}\langle \Ph W,QW\rangle.$  But
\bea
{\mbox Re}\langle \Ph W,QW\rangle&=&{\mbox Re}\{\mbox{tr} (\Ph WW^*Q^*)\}\nonumber\\
&=&{\mbox Re}\{\mbox{tr} (USV^*Q^*)\}\nonumber\\
&=&{\mbox Re}\{\mbox{tr} (SV^*Q^*U)\nonumber\\
&=&{\mbox Re}\left\{ \displaystyle{\sum_{i=1}^r}\s_it_{ii} \right\}.
\end{eqnarray}\nonumber

This sum is maximized when $t_{ii}=1$ for all $i,$ or when $V^*Q^*U=I_r,$ which yields $Q=UV^*.$
\end{myproof}

\vspace{.1in}

The goal of the above treatment is to enable us to say something about L\"owdin orthogonalization in the case of Gabor frames under a weighted norm.  That is, in the Banach space $L^2_s,$ where $s$ is the weight function, which orthogonalization of the Gabor frame $ (g,\a,\b)$ retains a maximal proportion of the excellent TF-localization of $g$?\\

Obviously, the result above shows that L\"owdin orthogonalization is optimal in the Banach space described by the norm $\|\cdot W\|_F$ where $W$ is any weight factor {\em in the matrix multiplication}.  We now formulate the discretized problem.  Write $g_{nm} :=M_{n\a}T_{m\b}g,$ the continuous translated modulated Gaussian, and write ${\overline g}_{nm} :=M_{n\a}T_{m\b}{\overline g}$ for the discretized frame element, say a row vector for concreteness.  Our linearly independent set $\Ph$ is the matrix whose rows are the ${\overline g}_{mn}$'s.  Note that (for the current treatment, at least) the weight function is constant for a given $m,$ even as $n$ varies; but as $m$ varies, the weight function must translate along with the gaussian.  Thus different translations of the weight function will multiply, point-by point, different row-blocks of $\Ph.$  This is a highly nonlinear operation, and thus the above treatment does not apply in its current form.  We shall denote point-by-point matrix multiplication by $\odot$ and the inverse operation, point-by-point matrix division, by $\oslash.$  Note that both $\odot$ and $\oslash$ conform to the (right or left) distributive rule: $(A+B)\odot C=A\odot C+B\odot C,\,(A+B)\oslash C=A\oslash C+B\oslash C$\\

Let $\Ph$ be as above, but the form of $W$ is now very specific.  Suppose we have a weight function $w:\R\rightarrow\R:t\mapsto w(t)$ with $w(t)>0,\,\forall\, t\in\R;$ a typical weight function is in the standard class of polynomial type $w_s(t)=(1+t^2)^s.$  Denote the discretization of $w$ by ${w^\circ}[d].$  Then the weight matrix $W$ will be of the form \\
\be
W^\circ=\bpm && {w^\circ}[n-\l_1T]&&\\&&\vdots&&\\&& {w^\circ}[d-\l_2T]&&\\&&\vdots&&\\&&\vdots&&\\&& {w^\circ}[d-\l_MT]&&\end{bmatrix}
\end{equation}
where $T$ is the symbol length, $M$ is the total number of translations in the disrete model, and each of the rows ${w^\circ}[d-\l_jT]$ is a particular translation of ${w^\circ}[d]$ and is repeated $N$ times, where $N$ is the finite number of modulations in the model.  Note that $W$ is quite rank-deficient, but more importantly, has no zero entries (in fact for the polynomial weights we have $W_{ij}\geq 1 \,\forall\,\, i,j$).\\

Let $\Phi$ be as above.  The goal is to find, among all unitary matrices $Q,$ the matrix $Q_{min}$ which minimizes $\|(\Ph-Q)\odot W\|_F.$  The method used in the proof of Lemma \ref{wtdlowdinmat} fails in the Hadamard multiplication case.  I will attempt to find a method which yields the optimal $Q$ for this case.

\section{Is $\Ll^1$ a symmetric Banach algebra?}\label{L1symm}
In Section \ref{HBSlemmas}, we noted that it seems as yet unknown whether or not $\Ll^1$ is a symmetric Banach algebra.  
\begin{conjecture}\label{conj1}
Let $h\in\Lone.$  Then there exists $\t_\d(x)=(1+|x|)^\d\,,\,\d>0$ such that $h_L\in\Ll^1_{\t_\d}.$
\end{conjecture}
If Conjecture \ref{conj1} can be proven, we might conclude that $\Ll^1$ is symmetric, provided that we can redefine $\Ll^1$ as the union $\underset{\t_\d\,,\,\d>0}{\bigcup}\Ll^1_{\t_\d}$ and assign a norm to each element of $h\in\Ll^1$ as $\|h\|_{\Ll^1}=\underset{\d}{\max}\,\,\|h\|_{\Ll^1_{\t_\d}}.$  In this case the hypotheses of Barnes' Lemma might be satisfied for $\Ll^1.$

\section{Extending the results regarding Banach algebras to quotient spaces}
It is possible to extend the results in Chapter \ref{decaydistributions} to Banach algebras of distributions containing non-invertible operators, by forming equivalence classes.  Fornasier and Gr\"ochenig demonstrate this concept for pseudo-inverses in \cite{FG04}, but it is, I think, possible to extend their result to the more general setting of quotient spaces, based on the ideas of Barnes in \cite{Bar87}.  Let $\Ll^1_v\circ$ and $\Ll_v\circ$ denote the extensions of $\Ll^1_v$ and $\Ll_v$ containing also non-invertible operators, and let $K_{\Ll^1_v\circ}$ and $K_{\Ll_v\circ}$ denote the kernels of $\Ll^1_v\circ$ and $\Ll_v\circ$ respectively.
\begin{conjecture}
The spaces $\Ll^{1'}_v=\Ll^1_v\circ/K_{\Ll^1_v\circ}$ and $\Ll'_v=\Ll_v\circ/K_{\Ll_v\circ}$ possess the same decay properties as $\Ll^1_v$ and $\Ll_v,$ respectively.
\end{conjecture}

\section{Numerical solution of integral equations and the finite section method}
It is often of interest to solve $Ax=b$ where $A$ is bi-infinite and $x,b$ are infinite.  Of course closed-form solutions are usually impossible, so we turn to the computer for our approximate solution.  For any sequence $y,$ let
\be
P_ny=(\ldots,0,0,y_{-n},y_{-n+1},\ldots,y_{n-1},y_n,0,0,\ldots).
\end{equation}
A typical approach to finding an approximate solution to $Ax=b$, termed the finite section method (see e.g., \cite{StrChr03}), is to truncate $A$ via the projection $P_nAP_n,$ which yields a finite section of $A$ padded with zeros outside the section defined by $n,$ and similarly with $P_nb.$  We then solve the system $P_nAP_n\tilde{x}=P_nb,$ where in general $\tilde{x}\neq P_nx.$ \\

The rate of convergence of $\tilde{x}\rightarrow x$ is of interest, whether in the discrete setting of matrix equations or the continuous setting of integral equations.  The results of Chapter \ref{decaydistributions} will have applicability to this question, since the kernels in $\Ll^1_v$ and $\Ll_v$ are defined by weights which by construction allow us to quantify convergence rates of numerical methods.
 
\section{Dependence of EBFDM SNIR on $\a\,,\,\b$ and $p$}
The dependence of EBFDM SNIR on $\a\,,\,\b$ and $p$ is not obvious (at least to the author), and warrants further investigation, as does the transition from the cusp-like maximum at $p=1/2$ for no dispersion to the smooth maximum for dispersive channels.

\appendix
\chapter{Integral Operators and Kernel Decay}
\markright{  \rm \normalsize APPENDIX A - INTEGRAL OPERATORS AND KERNEL DECAY}
\pagestyle{headings} 
\thispagestyle{headings} 

\section{Associativity of $\star$}\label{starassoc}

Let $f\star g=d, \,\,g\star h=p,\,\, (f\star g)\star h=q, \,\,f\star(g\star h)=r.$  The task is to show that $q=r.$  We have 
\be
d(x,z)=\int_{\R^d}f(t,z)g(x,t)\,dt\quad\mbox{ and }\quad p(y,\w)=\int_{\R^d}g(t,\w)\,h(y,t)\,dt\nn
\end{equation}
Thus
\bea
q(y,z)&=&\int_{\R^d}d(x,z)\,h(y,x)\,dx\nn\allowdisplaybreaks\\
&=&\int_{\R^d}\left(\int_{\R^d}f(t,z)g(x,t)\,dt\right)\,h(y,x)\,dx\allowdisplaybreaks\nn\\
&=&\int_{\R^d}\int_{\R^d}f(t,z)g(x,t)\,h(y,x)\,dt\,dx\allowdisplaybreaks\nn\\
&=&\int_{\R^d}f(t,z)\left(\int_{\R^d}g(x,t)\,h(y,x)\,dt\right)\,dx\allowdisplaybreaks\nn\\
&=&\int_{\R^d}f(t,z)\,p(y,t)\,dx\allowdisplaybreaks\nn\\
&=&r(y,z)\nn
\end{eqnarray}
and thus
\be
(f\star g)\star h=f\star (g\star h).
\end{equation}

\section{Schur's test}\label{schurstest}

\begin{lemma} Suppose that $f$ is a measurable function on $\R^{2d}$ which satifies
\bea
&&\underset{y\,\in\R}{\sup}\,\int_\R |f(t,y)|\,dt\leq C_1,\mbox{ and }\nn\\
&&\underset{t\,\in\R}{\sup}\,\int_\R |f(t,y)|\,dy\leq C_2,\label{schur}
\end{eqnarray}
Then the integral operator $f$ defined by $fh(y)=\int_\R f(t,y)h(t)\,dt$ is bounded from $\Lp$ to $\Lp$ for $1\leq p\leq\infty,$ and the operator norm of $f$ is bounded by
\be
\displaystyle{ \|f\|_{L^p(\R)\rightarrow L^p(\R)}\leq C_1^{1/p'}C_2^{1/p} }
\end{equation}
where $p'$ is the H\"older conjugate of $p.$
\end{lemma}

\begin{myproof} Since $1/p'+1/p=1,$ we can apply H\"older's inequality to $$\displaystyle{\int_{\R^d} f^{1/p'}(t,y)f^{1/p}(t,y) h(t)\,dt}$$
to get
\bea
|fh(y)|&\leq&\biggl(\int_{\R^d}|f^{1/p'}(t,y)|^{p'}\,dt\biggr)^{\frac{1}{p'}}\biggl(\int_{\R}|f^{1/p}(t,y)|^p|h(t)|^p\,dt\biggr)^{\frac{1}{p}}\nn\\
&\leq&\biggl(\int_{\R^d}|f(t,y)|\,dt\biggr)^{\frac{1}{p'}}\biggl(\int_{\R}|f(t,y)||h(t)|^p\,dt\biggr)^{\frac{1}{p}}\nn
\end{eqnarray}
which yields
\bea
\|fh\|^p_p&\leq&\int_\R\biggl(\int_{\R^d}|f(t,y)|\,dt\biggr)^{\frac{p}{p'}}\biggl(\int_{\R}|f(t,y)|\,|h(t)|^p\,dt\biggr)\,dy\nn\\
&\leq&C_1^{p/p'}\int_{\R^d}\int_{\R^d}|f(t,y)||h(t)|^p\,dt\,dy\nn\\
&\leq&C_1^{p/p'}\int_{\R^d}\int_{\R^d}|f(t,y)||h(t)|^p\,dy\,dt\qquad \mbox{(Fubini)}\nn\\
&\leq&C_1^{p/p'}\int_{\R^d}\biggl(\int_{\R^d}|f(t,y)|\,dy\biggr)|h(t)|^p\,dt\nn\\
&\leq&C_1^{p/p'}C_2\|h\|^p_p
\end{eqnarray}
and by taking $p^{th}$ roots, we are done.
\end{myproof}

\section{Properties of admissible weights}\label{weights1}

\subsection{Submultiplicativity of concave $\k$ with $\k(0)=0$}
Recall that $\k$ denotes the exponent in $v(x)=e^{\k(|x|)}.$  Since $\k(0)=0,$ by a simple corollary to Lemma 5.16 in \cite{Royd88} we have
$\frac{\k(b)}{b}\leq\frac{\k(a)}{a}\Rightarrow \frac{a\k(b)}{b}\leq\k(a),$ and
\bea
&&\frac{\k(a+b)}{a+b}\leq\frac{\k(b)}{b}\nn\\
&\Rightarrow&\k(a+b)\leq (1+\frac{a}{b})\k(b)\nn\\
&=&\k(b)+\frac{a}{b}\k(b)\leq\k(b)+\k(a),
\end{eqnarray}
and thus \eqref{submult} holds.\\

The proof that concave $\k$ with $\k\geq 0$ and $\k(0)=0$ is nondecreasing is menial but instructive.  We ignore the trivial case $\k\equiv 0.$  Since $\k(0)=0$ and $\k\geq 0,$ if $\k$ is decreasing at any point $y>0, \k(z)$ must be decreasing $\forall z>y,$ since if not, concavity would obviously be violated.  So suppose $\,\exists\,y\,\ni\,\k(z)$ is decreasing $\forall z>y.$  Then we must have $\k(z)\rightarrow 0$ as $z\rightarrow\infty.$  Let $y_1>y.$  Then since $\k$ is decreasing on $[y,\infty)$ we must have $\k(y_1)<\k(y).$  Write $m_1=\frac{\k(y)-\k(y_1)}{y-y_1}.$  Let $y_2>y_1-\k(y_1)/m>y_1$ (recall $m<0$).  Then since $\k(y_2)>0,\,\, m_2:=\frac{\k(y_1)-\k(y_2)}{y_1-y_2}<m_1,$ and thus $\k(ty+(1-t)y_2)\geq t\k(y)+(1-t)\k(y_2),$ so $\k$ is convex from $y$ to $y_2.$  Thus we conclude that if $\k(\x)$ is concave with $\k(0)\!=\!0,\, \k$ must be nondecreasing.\\ 

\subsection{$\k$ concave $\Rightarrow\,\,v$ continuous}
It suffices to prove that $\k$ is continuous, in light of the fact that the exponential function is continuous, and the composition of continuous functions is continuous.  Suppose $s<x<y<t.$  Now let $\la_1$ be such that $\la_1s+(1-\la_1)\,y=x,$ and $\la_2$ be such that $\la_2\,x+(1-\la_2)\,t=y.$  Obviously as $y\rightarrow x, \la_1\downarrow 0,$ and $\la_2\uparrow 1.$ We have
\bea
\k(x)&\geq&\la_1\,\k(s)+(1-\la_1)\k(y)\nn\\
\k(y)&\geq&\la_2\,\k(x)+(1-\la_2)\k(t).\nn
\end{eqnarray}
Since $s$ and $t$ are fixed, we have, as $y\rightarrow x,$
\bea
\k(x)&\geq&\underset{\la_1\downarrow 0}{\lim}\,\left\{\la_1\k(s)+(1-\la_1)\k(y)\right\}=\k(y)\nn\\
\k(x)&\leq&\underset{\la_2\uparrow 1}{\lim}\,\left\{\frac{1}{\la_2}\,\k(y)+\frac{(1-\la_2)}{\la_2}\k(t)\right\}=\k(y).\nn
\end{eqnarray}
Thus as $y\rightarrow x,$ we must have $\k(y)\rightarrow\k(x).$

\subsection{$v(x)=e^{\k(|x|)}$ is neither log-concave nor log-convex}  
Certainly $v$ restricted to certain subsets of its domain is log-concave, but $\forall\,\,x,$ we have
\bea
\log(v)(\frac{1}{2}x+(1-\frac{1}{2})(-x))=\k(0)=0\nn\\
\leq \k(\|x\|)=\frac{1}{2}\,\log(v)(x)+(1-\frac{1}{2})\log(v)(-x)\nn
\end{eqnarray}
violating log-concavity in general.

\begin{figure}[htb]
\begin{center}
\epsfig{file=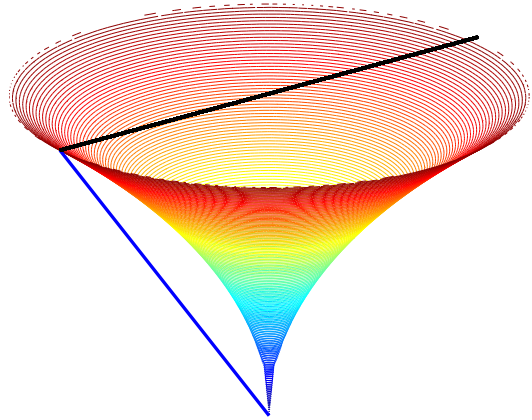,width=70mm}
\caption{{For admissible $v,\,\k(|x|)$ is neither concave nor convex}}
\label{fig:notconcconv}
\end{center}
\end{figure}

\section{The GRS-condition}\label{grsappendix}
To illuminate further the GRS-condition, note that $v(x)=e^{ ^{\sqrt{|x|}}}$ satisfies \eqref{grs}:
\be
\left(e^{ ^{\sqrt{|nx|}}}\right)^{1/n}=e^{ ^{\sqrt{\frac{|x|}{n}}}}\overset{n\rightarrow\infty}{\longrightarrow} 1\quad\forall x\in\R^d\nn
\end{equation}
but $v_\circ(x)=e^{\scriptscriptstyle{|}x\scriptscriptstyle{|}^2}$ does not:
\be
\left(e^{\scriptscriptstyle{|}nx\scriptscriptstyle{|}^2}\right)^{1/n}=e^{n\scriptscriptstyle{|}x\scriptscriptstyle{|}^2}\overset{n\rightarrow\infty}{\longrightarrow}\infty\quad\forall x\in\R^d.\nn
\end{equation}
Note, however, that condition \eqref{subexpdef} does not imply condition \eqref{grs}.  Let us restrict ourselves to $d=1$ and treat only $x\geq 0,$ to simplify the computations.  Consider the example $v=e^{\k(|x|)}$ where $\k(t)=t^{\,r(t)}$ and $r$ is asymptotic to $1$ from above:
\be
v_1(x)=e^{(x+1)^{\frac{x+2}{x+1}}-1}.\nn
\end{equation}
We denote $\k_1(x)=(x+1)^{\frac{x+2}{x+1}}-1.$  By construction $v(0)=e^{\k(0)}=1.$  The concavity of $\k_1$ can be shown by demonstrating that the second derivative $\k_1''$ is negative where $x>0$\footnote{I confess I have only shown this indirectly with Mathematica; the inequalities involved are unsolvable in closed-form.  The constrained maximum of $\k_1''(y)$ with $y\geq 0$ appears to be at $\k_1(0.295403)=-0.419112.$}.\\

But in checking the GRS-condition, we have 
\bea
\log\bigl(v_1(nx)^{1/n}\bigr)&=&\frac{1}{n}\left[(nx+1)^{\left(\frac{nx+2}{(nx+1)}\right)}-1\right]\nn\\
&=&\frac{1}{n}(nx+1)\cdot(nx+1)^{\frac{1}{(nx+1)}}-\frac{1}{n}\nn\\
&\geq&x\cdot(nx+1)^{\frac{1}{(nx+1)}}\nn\\
&\geq&x,\label{blowup}
\end{eqnarray}
since $y^{1/y}\geq 1$ if $y\geq 1.$  Thus
\be
\underset{n\rightarrow\infty}{\lim}v_1(nx)^{1/n}=\underset{n\rightarrow\infty}{\lim}e^{ ^{\scriptstyle{[}(nx+1)^{\scriptstyle{(}1+\frac{1}{(nx+1)}\scriptstyle{)} }-1\scriptstyle{]}}}\geq e^x >1,\label{nx}
\end{equation}
and so $v_1(x)$ does not satisfy the GRS-condition.

\section{An example of a subconvolutive weight}\label{subcon}

The admissible weight $\t_s(x)=(1+|x|)^s,\,s\geq 2d,$ is subconvolutive.  To show this, fix $x$ and divide $\R$ into two regions: 
\bea
t\in N_x\quad&\Rightarrow&\quad|t-x|\leq\frac{|x|}{2}\nn\\
t\in N^c_x\quad&\Rightarrow&\quad|t-x|>\frac{|x|}{2}.\nn
\end{eqnarray}
For $t\in N_x,$ we have 
\bea
&&-\frac{x}{2}\leq t-x\leq\frac{x}{2}\quad\Rightarrow\quad\frac{x}{2}\leq t\leq\frac{3x}{2}\nn\\
&\Rightarrow&|t|\geq\frac{|x|}{2}\quad\Rightarrow\quad(1+|t|)^{-s}\leq(1+\frac{|x|}{2})^{-s}.\nn
\end{eqnarray}
We conclude that
\be
(1+|t|)^{-s}\leq\frac{2^s}{(2+|x|)^s}\leq 2^s(1+|x|)^{-s}.
\end{equation}
Thus
\bea
\int_{N_x}(1+|t|)^{-s}(1+|x-t|)^{-s}\,dt&\leq&2^s(1+|x|)^{-s}\,\int_{N_x}(1+|x-t|)^{-s}\,dt\nn\\
&\leq&2^s(1+|x|)^{-s}\,\int_\R (1+|t|)^{-s}\,dt.
\end{eqnarray}
On $N^c_x,$ we have $1+|t-x|>1+\frac{|x|}{2},$ so $(1+|t-x|)^{-s}<(1+\frac{|x|}{2})^{-s}.$  Thus we have
\bea
\int_{N^c_x}(1+|t|)^{-s}(1+|x-t|)^{-s}\,dt&<&(1+\frac{|x|}{2})^{-s}\int_{N^c_x}(1+|t|)^{-s}\,dt\nn\\
&<&(1+\frac{|x|}{2})^{-s}\int_\R(1+|t|)^{-s}\,dt.\nn
\end{eqnarray}
If $s\geq 2,$ the integrals converge, and we have $\t_s^{-1}*\t_s^{-1}\leq C\t_s^{-1}$ as was required.

\section{Non-separability of $\Ll^1_v$ and $\Ll_v$}\label{nonsep}

Neither ${\mathcal L}_{v}$ nor ${\mathcal L}_{v}^1$ is separable.  To see this, consider the interval $I=[0,1].$  Represent each $z\in I$ by its binary equivalent $.d_1d_2d_3\hdots.$  Form the $\l^{\infty}(\Z)$ sequence $y_z$ via $(y_x)_k=d_k.$  ``Unfold'' $y_z$ by defining $\tilde{z}$ as follows:
$$\tilde{z}_k=\begin{cases}(y_z)_{-2k+1} & j \in\Bbb{Z},\,\, \mbox{$j \leq 0$}  ,\\
                               (y_z)_{2k}    & j \in\Bbb{Z},\,\, \mbox{$j>0$}.
                  \end{cases}$$
Then form a function $f_z=\underset{k\in\Z}{\sum}\chi_{[k,k+1)}(x)\chi_{[k,k+1)}(y)\tilde{z}_k.$  Thus in each square centered diagonally on the line $x=y,$ $f_z$ assumes only either of the constant values 1 or 0, and off the diagonal squares $f_z$ assumes only zero.\\

Denote the collection of such functions by ${\mathcal L}_{\{0,1\}}.$  It is clear that ${\mathcal L}_{\{0,1\}}\subset{\mathcal L}_{v}\cap {\mathcal L}_{v}^1.$  We have $\|f_{z_1}-f_{z_2}\|_{{\mathcal L}_v}>1$ and $\|f_{z_1}-f_{z_2}\|_{{\mathcal L}_v^1}\,>\, 1 \,\,\forall\,\,f_{z_1},f_{z_2}\in {\mathcal L}_{\{0,1\}}$ with  $z_1\neq z_2$ since $v(0)=1$ is a global minimum.     So, in the topology generated by either norm, $B_{\frac{1}{3}}(f_{z_1})\cap B_{\frac{1}{3}}(f_{z_2})=\emptyset.$  Since $[0,1]$ is uncountable, there exist an uncountable number of balls of radius $\frac{1}{3},$ containing only singletons, around the elements of ${\mathcal L}_{\{0,1\}}$ (and ${\mathcal L}_{\{0,1\}}$ is not even dense in ${\mathcal L}_{v}\cap {\mathcal L}_{v}^1$). Thus ${\mathcal L}_{\{0,1\}}$ is not separable, and hence neither is ${\mathcal L}_v$ nor ${\mathcal L}_{v}^1.$  

\np

\pagestyle{myheadings} 
\chapter{Gabor analysis, LOFDM, and the Stone-von Neumann theorem}
\markright{  \rm \normalsize APPENDIX B - GABOR ANALYSIS, LOFDM, AND THE SVN THEOREM}
\thispagestyle{myheadings} 

\section{Powers of the Gabor frame operator commute with TF-shifts}\label{spowercommute}
For all $t\in\R\backslash\{1\},$ let $\Upsilon$ be the commutator
\be
\label{allcommute}
\Upsilon:=S^tT_xM_\w - T_xM_\w S^t.
\end{equation}
Now $S^{1-t}\Upsilon S^{1-t}=0,$ and since the frame operator $S$ is a positive operator, the operators $S^\nu,\,\nu\neq 0,$ are well-defined and positive by the spectral theorem and thus invertible. We conclude that $\Upsilon=0.$   In particular $S^{-\frac{1}{2}}$ commutes with the time-frequency-shifts as well, and we have
\be
\label{gabtgtgab}
S^{-\frac{1}{2}} (g,\La)= (S^{-\frac{1}{2}}g,\La).
\end{equation}

\section{Simultaneous pulse and lattice dilation commutes with L\"owdin orthogonalization of the Gaussian}\label{dilorth}

We prove \eqref{orthdil}.  First note that the Gram matrix of $(\ggam,1/\sqrt{\sigma \rho},\sqrt{\sigma/\rho})$ equals the Gram  matrix for $(g,1/\sqrt{\p},1/\sqrt{\p}).$  Consider an entry in the Gram matrix of dilated Gaussians on the dilated lattice:
\be
\langle g^\s_{k'\frac{\a}{\sqrt{\s}},\l'\b\sqrt{\s}}\,, \,g^\s_{k\frac{\a}{\sqrt{\s}},\l\b\sqrt{\s}}\rangle=\int_\R(2\s)^{\frac{1}{4}}e^{-\pi\s (t-k'\frac{\a}{\sqrt{\s}})^2}e^{2\pi i\l'\b\sqrt{\s}t} (2\s)^{\frac{1}{4}}e^{-\pi\s (t-k\frac{\a}{\sqrt{\s}})^2}e^{2\pi i\l\b\sqrt{\s}t}\,dt.
\end{equation}
By changing variables to $\t=\sqrt{\s}t,$ we get
\be
\int_\R 2^{\frac{1}{4}}e^{-\pi(\t-k'\a)^2}e^{2\pi i\l'\b\t}\,\,2^{\frac{1}{4}}e^{-\pi(\t-k\a)^2}e^{2\pi i\l\b\t}\,d\t=\langle g_{k'\a,\l'\b}\,, \,g_{k\a,\l\b}\rangle.
\end{equation}
From this we can conclude that for the Gram matrices derived from dilated and undilated Gaussians, we have
\be
R^{-1/2}_{k\frac{\a}{\sqrt{\s}},\l\b\sqrt{\s},0,0}=R^{-1/2}_{k\a,\l\b,0,0}
\end{equation}
leading to
\be
\sum_{k,\l}R^{-1/2}_{k\frac{\a}{\sqrt{\s}},\l\b\sqrt{\s},0,0}\ggam_{k,\l}=
\sum_{k,\l}R^{-1/2}_{k\a,\l\b,0,0}\ggam_{k,\l}.
\end{equation}
Now since the dilation operator $\D$ is a bounded linear operator on $\LtR,$ we have from, e.g., Lemma 5.3.2 in \cite{Gro01} that
\be
\psigam=\sum_{k,\l}R^{-1/2}_{k\a,\l\b,0,0}\,\,\ggam_{k,\l}=\D_{\sqrt{\s}}\sum_{k,\l}R^{-1/2}_{k\a,\l\b,0,0}\,\,g_{k,\l}=\D_{\sqrt{\s}}\,\,\ph.
\end{equation}

\section{Representations, The Stone-von Neumann theorem, and the Heisenberg group}\label{repSvNHg}
A {\it representation} $\z$ of a group $G$ on a vector space $V$ is defined as a homomorphism
\be
\z : G \longrightarrow GL(V),
\end{equation}
with $dim\,\z = dim\,V,$ and where $GL(V)$ denotes the group of linear operators on $V.$   The rotation groups on $\R^3$ are examples of representations on a vector space without regard to a choice of basis.  If we choose a basis, this defines an isomorphism
\bea
GL(V)&\longrightarrow& GL_n(F),\\
T &\longrightarrow& (\mbox{matrix of } T),
\end{eqnarray}
where $GL_n(V)$ denotes the group of invertible $n\times n$ matrices and $F$ is the field over which the chosen basis lies.  In this way we arrive at $\zeta\longrightarrow GL(V)\longrightarrow GL_n(F)$ which is a {\it matrix representation} of $G.$  A matrix representation can be viewed as a representation of $G$ on the space of column vectors $F^n.$\\

Let $\z$ be a representation of $G.$  Denote the image of $g\in G$ under $\z$ as $\z_g.$  Then $\z_g$ is a linear operator on $V,$ and for $g,h\in G$ the associative law yields $\zeta_{gh}=\z_g\z_h.$\\

\begin{definition}(Unitary Representations):  A {\it Hermitian} vector space is defined as a complex vector space together with a positive-definite Hermitian form $\langle\cdot ,\cdot\rangle.$  Recall that an operator $T$ on $V$ is unitary $\Leftrightarrow \langle x,y\rangle = \langle Tx,Ty\rangle ,\, \forall x,y \in V.$  It is thus natural to define {\em unitary representations} via $\langle u,w\rangle = \langle \zeta_gu,\zeta_gw\rangle ,\,\, \forall u,w \in V,\, \forall g\in G.$\\
\end{definition}

Consider the dihedral group $D_4$ acting on $\R^3$; there are two proper invariant subspaces - the plane containing the square $\Delta$ on which $D_4$ is defined, and the axis of rotation perpendicular to $\Delta.$ \\

A {\em proper} $G$-{\em invariant subspace} of $V$ under $\z$ is a subspace of $V,$ not equal to $V,$ which is invariant under the action of the representation $\z.$\\

\begin{definition}(Irreducible Representations):  If a representation $\z$ of a group $G$ has no proper $G-$invariant subspace, then $\z$ is called {\it irreducible}. \\
\end{definition}

\subsection{Schr\"odinger's Representation of the Heisenberg Group}

\noindent Define on $\H'=(\R\times\R\times\T)$ a group multiplication by
\be
(x,y,\tau)\cdot (u,v,\eta ) = (x+u,y+v,e^{2\pi i(\tau+\eta)}e^{i \pi (uy-xv)}).
\end{equation}
This is the {\it reduced Heisenberg group}.  Further, let $\H=(\R\times\R\times\R)$ and define a multiplication
\be
(x,y,\tau)\cdot (u,v,\eta ) = (x+u,y+v,\tau+\eta+\frac{1}{2}(uy-xv)),
\end{equation}
and the result is the {\it full} Heisenberg group.  Both the full and reduced Heisenberg groups have the same distinguished unitary representation, the Schr\"odinger representation.  The Schr\"odinger representation $\pi$ from $\H$ into the unitary operators on $\LtR$ is given by
\be
\pi (x,y,\tau)\,f(t)= \tau e^{\pi ixy} e^{2\pi iyt} f(t+x),\quad x,y,t\in \R, \, \tau\in \T,
\end{equation}
or
\be
\pi (x,y,\tau)=\tau e^{\pi ixy} M_yT_{-x}.\label{Schro1}
\end{equation}
The Schr\"odinger representation is clearly unitary, and is irreducible since if we suppose ${\mathcal K}\neq 0$ is an invariant subspace of ${\mathcal H}$ under $\pi$ and fix $g\neq 0$ in ${\mathcal K}$ and let $f\in {\mathcal K}^\perp,$ then the {\it representation coefficients} equal $0$ for all ${\mathbf h}=(x,y,\tau) \in \H$:
\be
0=|\langle f,\pi({\mathbf h})g\rangle |=|V_gf(x,y)|,\quad\forall\,(x,y,\tau) \in \H
\end{equation}
and since the STFT is injective, $f\equiv 0,$ so ${\mathcal K}^\perp=0$ and ${\mathcal K}={\mathcal H}.$\\

For convenience we will often write $z=(x,y)$ for $x,y\in\R^2$ so $\pi=\pi (z,\tau).$  Every lattice in $\R^2$ is of the form $\La=L(\Z^2)$ for some $L\in GL_2(\R).$  Since we use the tools of representation theory of the Heisenberg group, we must first determine which automorphisms $z\mapsto Lz$ of $\R^2$ extend to automorphisms $i_L(z,\tau)=(Lz,\tau)$ of $\H.$  If $L$ leaves the symplectic form $[(x,y),(u,v)]=yu-xv$ on $\R^2$ invariant, $L$ defines an automorphism of $\H$ in the following way:
\be
i_L (z,\tau)=(Lz,\tau),\quad \forall\, z\in\R^2,\,\tau\in\T,
\end{equation}
yielding the sequence of mappings
\be
\H\overset{i_L}{\longrightarrow}\H \overset{\pi}{\longrightarrow}{\mathcal U}_{\mathcal H}:{\mathcal H}\overset{U_L}{\longrightarrow}{\mathcal H}
\end{equation}
Such a matrix $L$ is characterized by\footnote{In higher dimensions, this characterization defines the symplectic group.}
\bea
&&det\,\left( L\cdot\begin{bmatrix}y&v\\x&u\end{bmatrix}\right)=\pm\,det\,\begin{bmatrix}y&v\\x&u\end{bmatrix}\\
&&\Rightarrow det\,L=\pm 1\\
&&\Rightarrow L\in SL_2(\R).
\end{eqnarray}
\begin{theorem}{\em (Stone-von Neumann)} \label{svnthm}
Any irreducible unitary representation $\zeta$ of $\H$ on a Hilbert space ${\mathcal H}$ satisfying $\zeta (0,\tau)f=\tau f,\,\tau\in\T,\, f\in {\mathcal H},$ is unitarily equivalent to the Schr\"odinger representation $\pi.$  Thus there exists unitary $U : {\mathcal H} \rightarrow \LtR$ such that $\zeta (x,y,\tau ) = U^{-1} \pi (z,\tau) U,\,\,\forall (z,\tau)\in\H.$
\end{theorem}

For $L\in SL_2(\R)$ we compose the Schr\"odinger representation with the automorphism $i_L$ to yield a new representation $\pi_L:=\pi\circ i_L$ given by $\pi_L(z,\tau)=\pi(Lz,\tau).$
Clearly $\pi_L$ is unitary, and is irreducible since $L$ is a unitary change of coordinates, so (see the proof irreducibility for $\pi$ above) the appeal to the STFT's injectivity still holds.   Since $\pi_L(0,\tau)=\pi(0,\tau),\,\forall\,\tau\,\in\,\T,$ then by Theorem \ref{svnthm}, there exists unitary $\mu(L),$ called a {\em symplectic operator}, (defined only up to a phase factor) such that 
\be
\pi_L(z,\tau)=\mu(L)\pi(z,\tau)\mu(L)^{-1}\quad\forall (z,\tau)\in\H.\label{unitarySVN}
\end{equation}
The phase factors of the symplectic operators $\mu(\cdot)$ can be defined so that either $\mu(\A\B)=\mu(\A)\mu(\B)$ or $\mu(\A\B)=-\mu(\A)\mu(\B)$ (but neither case can hold universally), by invoking $Mp(d),$ the two-fold cover of $Sp(d).$  The following theorem is due to Shale and Weil (\cite{Sh62},\cite{We64}).
\begin{theorem}
Let $Mp(d)$ be the two-fold cover of $Sp(d),$ and let $\tilde{\A}\rightarrow\A$ be the projection from $Mp(d)$ onto $Sp(d)\cong Mp(d)/\{\pm I\}.$  Then there exists a unitary representation $\mu$ of $Mp(d)$ on $\Ltwo$ such that
\be
\mu(\tilde{\A})\pi(z,\t)\mu(\tilde{\A})^{-1}=\pi(\A z,\t).
\end{equation}
\end{theorem}

We call $Mp(d)$ the {\em metaplectic group} and $\mu$ the {\em metaplectic representation}.

\subsection{Examples relevant to Sections \ref{ss:gen} and \ref{ss:ex}}
In this section, we drop the reference to $\t$ in \eqref{unitarySVN}, since it is of no consequence to the calculations and its inclusion would only be cumbersome.  Note that each of the following matrices $A,B,$ and $C$ are members of $Sp(1).$\\

\begin{example}
Let $A=\begin{bmatrix}0&1\\-1&0\end{bmatrix};$ then $\mu(A)={\mathcal F}.$  Indeed, by equations  \eqref{commrels}, \eqref{transfourier}, \eqref{modfourier}, and \eqref{Schro1}, we have
\bea
{\mathcal F}\pi(x,\omega)&=&e^{\pi ix\omega}{\mathcal F}T_xM_\omega=e^{\pi ix\omega}M_{-x}T_\omega{\mathcal F}\nn\\
&=&e^{-\pi ix\omega}T_\omega M_{-x}{\mathcal F}=\pi(A(x,\omega)){\mathcal F}.
\end{eqnarray}
\end{example}

\begin{example}
Let $B=\begin{bmatrix}\alpha&0\\0&\alpha^{-1}\end{bmatrix};$ and set ${\mathcal D}_\alpha f(x) = \alpha^{1/2}f(\alpha x).$  Then we have
\bea
{\mathcal D}_{\alpha^{-1}}T_xM_\omega f(t)&=&\alpha^{-1/2}e^{2\pi i\omega(\alpha^{-1}t-x)}f(\alpha^{-1}t-x)\nn\\
&=&\alpha^{-1/2}e^{2\pi i(\alpha^{-1}\omega)(t-\alpha x)}f(\alpha^{-1}(t-\alpha x))\nn\\
&=&T_{\alpha x}M_{\alpha^{-1}\omega}{\mathcal D}_{\alpha^{-1}}f(t).
\end{eqnarray}
Thus $\mu(B)={\mathcal D}_{\alpha^{-1}}.$
\end{example}

\begin{example}
Let $C=\begin{bmatrix}1&0\\\b&1\end{bmatrix}$ with $\b=\bar{\b},$ and set $\Chirp=e^{-\pi i\b t^2}f(t).$  Then $\mu(C)=\Chirp^{-1}$ results from
\bea
{\Chirp}^{-1}\pi(x,\omega){\Chirp}f(t)&=&e^{\pi i\b t^2}e^{\pi ix\omega}e^{2\pi i\omega(t-x)}e^{-\pi i\b(t-x)^2}f(t-x)\nn\\
&=&e^{\pi ix(\omega+\b x)}e^{2\pi i(\omega +\b x)(t-x)}f(t-x)\nn\\
&=&\pi(x,\omega +\b x)f(t)=\pi(C(x,\omega))f(t).
\end{eqnarray}
\end{example}

\subsection{Extension to general symplectic lattices}\label{gensympl}
The following is from \cite{Gro01}.  In the context of frame theory, we have, for a window $g,$ arbitrary $f\in\LtR,$ and an arbitrary lattice $\La$
\be
S_{g,\La}f=\underset{\la\in\La}{\sum}\langle f,\pi(\la,\tau)g\rangle\,\pi(\la,\tau)g.
\end{equation}
\begin{lemma}
Suppose $\La=M\Z^2,\,g\in S_0(\R)$ (functions whose STFT is integrable).  Writing $M=\alpha L,\,\textup{det}\,L=1,\,\La_1=(\alpha\Z)^2$ we have
\be
S_{g,\La} = \s(L)S_{\s(L)^{-1}g,\La_1}\s(L)^{-1},
\end{equation}
where $\s$ is the (unitary) metaplectic representation.\\

\begin{myproof}
Write $S=S_{g,\La}.$  Then 
\bea 
S(f)&=&\underset{\la\in\La}{\sum}\langle f ,\pi (\la)g\rangle\,\pi (\la)g\\
&=&\underset{k\in\Z^2}{\sum}\langle f ,\pi (L\alpha k)g\rangle\,\pi (L\alpha k)g\\
&=&\underset{k\in\Z^2}{\sum}\langle f ,\s(L)\pi (\alpha k)\s(L)^{-1}g\rangle\,\s(L)\pi (\alpha k)\s (L)^{-1}g\\
&=&\s(L)\underset{k\in\Z^2}{\sum}\left\langle \s(L)^{-1}f ,\pi (\alpha k)\s(L)^{-1}g\right\rangle\pi (\alpha k)\s (L)^{-1}g\\
&=&[\s(L)\,S_{\s(L)^{-1}g,\La_1}\,\s(L)^{-1}](f).
\end{eqnarray}
\end{myproof}
\end{lemma}

Thus the frame operator on a general lattice is unitarily related to the frame operator on a separable lattice.  This fact permits us to conclude that all results on Gabor frames extend to the general lattice case.\\

\np

\pagestyle{myheadings} 
\chapter{Glossary of notation}
\markright{  \rm \normalsize APPENDIX C - GLOSSARY OF NOTATION}
\thispagestyle{myheadings} 
\renewcommand{\baselinestretch}{1.6}\large
In this appendix we provide a referral for the notation commonly used in this paper.
\section{Basics}
\begin{tabbing}
$\mathcal{H}\hspace{.9in}$\= arbitrary Hilbert space\\
${\mathcal B}({\mathcal H})$ \> the space of bounded linear operators on ${\mathcal H}$\\
$\N$ \> the nonnegative integers\\
$\Z$ \> the integers\\
$\Z_+$ \> the positive integers\\
$\Z^d$ \> the set of $d$-tuples of integers\\
$\Q$ \> the rational numbers\\
$\R$ \> the real numbers\\
$\R_+$ \> the positive real numbers\\
$\R^d$ \> the set of $d$-tuples of real numbers\\
$i$ \> $\sqrt{-1}$\\
$\C$ \> the complex numbers\\
$\overline{y}$ \> the complex conjugate of $y=a+bi$: $\overline{y}=a-bi$\\
$\exists$ \> ``there exists''\\
$\forall$ \> ``for all''\\
$\emptyset$ \> the empty set\\
$\wt{\O}$ \hspace{.9in}\= the space of all functions on $\R^{2d}$ \\
$C^\infty(\R^d)$ \> the space of all infinitely differentiable functions on $\R^d$\\
$C_0$ \> $\{f\,|\,f\,\text{ continuous and }\,f(x)\rightarrow 0 \,\text{ as }\,|x|\rightarrow\infty\}$\\
$|x|$ \> the finite-dimensional 2-norm of $x$\\
$\a$ \> $(\a_1,\ldots,\a_d)\in\Z^d_+;$ a multi-index \\
$x^\a$ \> $\prod^d_{j=1}x^{\a_j}_j$\\[1mm]
$D^\a$ \> $\left(\frac{\partial}{\partial x_1}\right)^{\a_1}\cdots\left(\frac{\partial}{\partial x_d}\right)^{\a_d}$\\
$\S$ \> the Schwartz space; all $\ph\in\C^\infty(\R^d)$ such that\\[-2mm] 
          \> $\underset{x\in\R^d}{\sup}\left|x^\a D^{\,\b}\ph(x)\right|<\infty$\\[1mm]
$\S^*$ \> the space of tempered distributions, the dual of $\S$\\
$\d$ \> the Dirac delta function (point measure)\allowdisplaybreaks\\
$\F$ \> the Fourier transform operator\\
$\hat{f}(\w)$ \> $(\F f)(\w)$\\
STFT \> the short-time Fourier transform\\
$Af$ \> the ambiguity function of $f$\\
$A(f,h)$ \> the cross-ambiguity function of $f$ and $h$\\
$f\star g$ \> composition of integral operators with kernels $f$ and $g,$ also \\[-1.5mm]
    \> frequently denoted $fg$\\
$f*g$ \> the convolution of $f$ with $g$\\
$W(\R^d)$ \> the Wiener algebra : $f\in W(\R^d)\Leftrightarrow f,\hat{f}\in\Lone$\\
$\int f(x)\,dx$ \> the Lebesgue integral of the function $f$\\
$\|f\|$ \> the norm of $f$\\
$\|f\|_{op}$ \> the operator norm of $f$ given by $\underset{\|x\|=1}{\sup}\|f(x)\|$ \\
$\|f\|_p$ \> the $p$-norm of $f,$ given by $\|f\|_p = \left(\int_{\R^d}|f(x)|^p\,dx\right)^\frac{1}{p}$ \\[6mm]
$\Lp$ \> the space of $p^{th}$-power Lebesgue-integrable functions on $\R^d,$\\[-2.5mm]
      \> $1\leq p\leq\infty$\\
$\underset{x\in\R^d}{\mbox{ess }\sup}\,f(x)$ \> the essential supremum of $f,$ given by\\[-3.5mm]                     \> $\underset{x\in\R^d}{\mbox{ess }\sup}\,f(x)=\inf\{M\,|\,m\{x\,|\,f(x)>M\}=0\}$\\
$\|f\|_\infty$ \> the infinity-norm of $f,$ given by $\|f\|_\infty :=\underset{x\in\R^d}{\mbox{ess }\sup}\,|f(x)|$ \\
$\|A\|_F$ \> the Frobenius norm of the matrix A, given by $\sum_{i=1}^m\sum_{j=1}^n |a_{ij}|^2$\\
$A^c$ \> the complement of $A$ in some superspace
\end{tabbing}

\section{Banach algebras and off-diagonal decay}
\begin{tabbing}
$\d(x-y)$ \hspace{.5in}\= the kernel of the identity integral operator on $f:\R^d\rightarrow\R;$  \\[-2.5mm]
    \> also denoted $e$ \\
$\O$  \> the space of all functions on $\R^{2d}$ appended with $\d(x-y)$\\
$\O_1$  \> a superspace of $\O$ admitting point measures, of which no \\[-2.5mm]
    \>hyperplane can contain an infinite number\\
$v$ \> an admissible weight, which satisfies the GRS-condition and \\[-2.5mm]
    \> whose exponent $\k$ satisfies concavity on $|x|$ \\
$\k$ \> the function in the exponent of $v\,:\,v(x)=e^{\k(|x|)}$\\
$f\odot g$ \> the Hadamard (or pointwise) multiplication of $f$ and $g$\\
$\Ll^1_v$ \> the space of equivalence classes of all operators in $\O_1$ which are\\[-2.5mm]
        \>integrable against the weight $v$ in both coordinates\\
$\Ll^1$ \> $\Ll^1_v$ with $v\equiv 1$\\
$\Ll_v$ \> the space of equivalence classes of all operators in $\O_1$ which are \\[-2.5mm]
          \> ess-sup-bounded pointwise against $v$ in both coordinates\allowdisplaybreaks\\
$\t_s$ \> the polynomial weight $(1+|x|)^s$ with $s>0$\\
$\B_{u,s}$ \> derived from $\Ll_v$ : $v=u\t_s$ with $u$ admissible;                                 $\B_{u,s}=\Ll^1_u\cap\Ll_v$\allowdisplaybreaks\\
$f_L$  \> the convolution (Laurent) operator generated from $f:\R^d\rightarrow\R$\allowdisplaybreaks\\
$\tilde{f}$ \> the non-distributional part of the kernel $f$ of an integral operator\\
$\s_A(f)$ \> the algebraic spectrum of the operator $f$\\ 
$\s(f)$ \> the spectrum (as an operator on $\mathcal{H}$) of the operator $f$\\
$\r_A(f)$ \> the algebraic spectral radius of the operator $f$\\ 
$\r(f)$ \> the spectral radius (as an operator on $\mathcal{H}$) of the operator $f$
\end{tabbing}

\section{Pseudodifferential operators}
\begin{tabbing}
$\Hchan$ \hspace{.9in} \= time-varying communications channel operator\\
$K_\s$ \> a pseudodifferential operator\\
$\s$ \> the symbol of $K_\s$\\
$k_\s$ \> the kernel of $K_\s$\\
$\F_2$ \> the partial Fourier transform in the second variable\\
$A$ \> a linear time-varying system
\end{tabbing}

\section{OFDM and LOFDM}
\begin{tabbing}
OFDM \hspace{.5in} \= orthogonal frequency-division multiplexing\\
BFDM \> biorthogonal frequency-division multiplexing\\
EBFDM \> extended biorthogonal frequency-division multiplexing\\
SNIR \> signal-to-noise-plus-interference ratio\\
AWGN \> additive white Gaussian noise\\
TFL \> time-frequency localization\\
$T$ \> the symbol length for a rectangular OFDM system\\
$F$ \> the carrier separation for a rectangular OFDM system\\
$\r$ \> the spectral efficiency of a communications system \\
$C_\Hchan$ \> scattering function of the channel $\Hchan$\\
$S_\Hchan$ \> spreading function of $\Hchan$ (randomized version of $C_\Hchan$)\\
$h(t,\t)$ \> time-varying impulse response of $\Hchan$\\
$\t_0$ \> maximum delay due to multipath propagation\\
$\nu_0$ \> maximum Doppler spread\\
$\Delta\t_f$ \> second-order time-localization of the function $f$\\
$\Delta\nu_f$ \> second-order frequency-localization of the function $f$\\
$\psi$ \> the atom, or generating pulse-shape, of an OFDM system\\
$f_{k,\l}$ \> $f(t-kT)e^{2\pi i\l Ft};$ notation used only for rectangular OFDM systems\\
$T_x$ \> the translation operator: $T_xf(t)=f(t-x)$\\
$M_\w$ \> the modulation operator: $M_\w f(t)=f(t)e^{2\pi i\w\cdot t}$\\
$z(\la)$ \> $T_xM_\w : \la=(x,\w)\in\La$\\
$\La$ \> an arbitrary symplectic lattice\\
$L$ \> the generator matrix for the lattice $\La$\\
$\xi(\La)$ \> the density of the lattice $\La$\\
$\La_\circ$ \> the adjoint lattice of $\La$ given by $\xi^{-1/d}\La$\\
$J$ \> an arbitrary index set\\
$e_j$ \> a frame element\\
$S$ \> the general frame operator: $Sf = \sum_{j\in J}\langle f,e_j\rangle e_j$\\
$S^\nu(g,\La)$ \> the action of the $\nu^{th}$-power of the Gabor frame operator\\[-2.5mm]
    \> $S$ on the entire Gabor frame $(g,\La)$\\
$S_{f,h,\La}$ \> the frame operator for the frame generated by the dual functions\\[-2.5mm]
    \> $f$ and $h$ on the lattice $\La$\\
$(f,\La)$ \> Gabor system generated by $f$ on the lattice $\La;$ also used to\\[-2.5mm]
  \> denote an LOFDM system\\
$L\ddot{o}(f,\La)$ \> The L\"owdin-orthogonalized pulse derived from $f$ on an \\[-2.5mm]
  \> arbitrary time-frequency lattice $\La$ \\
$(f,T,F)$ \> same as $(f,\La)$ but with $\La$ rectangular with parameters $T$ and $F$\\
$L\ddot{o}(f,T,F)$ \> $L\ddot{o}(f,\La)$ where $\La$ is rectangular with parameters $T$ and $F$ \\
$\phi$ \> $L\ddot{o}(f,\La)$\\
$\ph_\s$ \> $L\ddot{o}(g_\s,\La)$ where $g_\s$ is a dilated Gaussian\\
$\ph$ \> $L\ddot{o}(g_\s,\La)$ where $\s=1$\\
$\supp_{\e}f$ \> the effective support of $f,$ where $|f|>\e$\\
$R_{f,\La}$ \> the Gram matrix of $(f,\La)$\\
$R_{f,T,F}$ \> the Gram matrix of $(f,T,F)$\\
$\Dil f(t)$ \> the dilation operator $\sqrt{\alpha} f(\alpha t),\,\alpha \in \R_+$\\
$\Chirp f(t)$ \> multiplication by a chirp: $f(t) e^{-\pi i \beta t^2},\,\beta \in \R$\\
{\em SL}$(2,\R)$ \> the special linear group; $2\times 2$ real matrices $A$ with $|$det$(A)|=1$ 
\end{tabbing}
 
\section{Future research and appendices}
\begin{tabbing}
$W$ \hspace{.9in}\=a weight matrix under the usual matrix multiplication\\
$W^\circ$ \> a weight matrix under Hadamard matrix multiplication\\
$B_\e(f)$ \> the ball of radius $\e$ about a function $f$\\
$[a,b]$ \> the set of real numbers $x$ such that $a\leq x\leq b$\\
$V$ \> a vector space\\
$\z$ \> a reperesentation of a group $G$ on $V$\\
$\z_g$ \> the image of $g\in G$ in the vector space $V$ under the representation $\z$\\
$\T$ \> the torus: $\{(x,y)\in\R^2\,|\,0\leq x,y\leq 1\}$ with $(0,y)$ and $(1,y)$ identified\\[-2.5mm]
    \> and $(x,0)$ and $(x,1)$ identified\\
$\H'$ \> $(\R\times\R\times\T)$\\
$\H$ \> $(\R\times\R\times\R)$\\
$\pi$ \> the Schr\"odinger representation $\pi (x,y,\tau)=\tau e^{\pi ixy} M_yT_{-x},$\\[-2.5mm]
      \> $x,y,t\in \R,\,\tau\in\T$ from either $\H'$ or $\H$ to the unitary operators\\[-2.5mm]
      \> on $\Ltwo$\\
$\pi_L(z,\tau)$ \> the representation $\pi(Lz,\tau)\,\,\forall\,\,(z,\tau)\in\H$ for $L\in SL_2(\R)$ \\
{\em Sp}$(d,\R)$ \> the symplectic group, also denoted {\em Sp}$(d)$\\
{\em Mp}$(d,\R)$ \> the metaplectic group, the twofold cover of {\em Sp}$(d)$\\
$\mu(\cdot)$ \> the symplectic (unitary) operator appearing in the metaplectic\\[-2.5mm]
    \> representation of {\em Mp}$(d,\R)$
\end{tabbing}


\np

\pagestyle{myheadings} 
\markright{  \rm \normalsize BIBLIOGRAPHY}

\end{document}